\documentclass[11pt]{article}
\usepackage{amsmath}
\usepackage{amsfonts}
\usepackage{amsthm}
\usepackage{pictex}
\usepackage{graphicx}
\usepackage{comment}
\usepackage{a4wide}

\def\Ai{\mathop{\mathrm{Ai}}\nolimits}

\newtheorem{theorem}{Theorem}[section]
\newtheorem{proposition}[theorem]{Proposition}
\newtheorem{lemma}[theorem]{Lemma}
\newtheorem{corollary}[theorem]{Corollary}
\newtheorem{Definition}[theorem]{Definition}
\newenvironment{definition}{\begin{Definition}\rm}{\end{Definition}}
\newtheorem{Remark}[theorem]{Remark}
\newenvironment{remark}{\begin{Remark}\rm}{\end{Remark}}
\numberwithin{equation}{section}
\newcommand{\C}{\mathbb{C}}

\renewcommand{\O}{\mathcal{O}}
\renewcommand{\Re}{{\rm Re} \,}
\renewcommand{\Im}{{\rm Im} \,}

\newcommand{\ds}{\displaystyle}

\begin{document}
\title{Quadratic Hermite--Pad\'e approximation to the exponential function:
a Riemann--Hilbert approach\thanks{This work is supported by INTAS 2000-0272,
 projects G.0176.02 and G.0184.02 of FWO-Vlaanderen}}
\author{A. B. J. Kuijlaars, W. Van Assche,  \\
  Katholieke Universiteit Leuven, Belgium \\[10pt]  and  \\[10pt]
  F. Wielonsky \\
  Universit\'e des Sciences et Technologies Lille 1, France, and \\
  INRIA Sophia Antipolis, France}
\maketitle

\begin{abstract}
We investigate the asymptotic behavior of the polynomials
$p$, $q$, $r$ of degrees $n$ in type I Hermite-Pad\'e approximation
to the exponential function, defined by
$p(z)e^{-z} +q(z) + r(z) e^{z} = \O(z^{3n+2})$ as $z \to 0$.
These polynomials are characterized by a Riemann--Hilbert problem
for a $3\times 3$ matrix valued function. We use the
Deift-Zhou steepest descent method for Riemann--Hilbert problems
to obtain strong uniform asymptotics
for the scaled polynomials $p(3nz)$, $q(3nz)$, and $r(3nz)$
in every domain in the complex plane.
An important role is played by a three-sheeted Riemann surface
and certain measures and functions derived from it.
Our work complements recent results of Herbert Stahl.
\end{abstract}

\tableofcontents
\section{Introduction}

The study of Hermite--Pad\'e approximation  for the exponential function
was initiated by Hermite \cite{hermite},
in connection with his proof of the transcendency of $e$.
A special case of this is the quadratic Hermite--Pad\'e approximation
of type I, where for any given three integers $n_1,n_2,n_3 \geq 0$,
one asks for polynomials $p_{n_1,n_2,n_3}$, $q_{n_1,n_2,n_3}$, and
$r_{n_1,n_2,n_3}$ of degrees $n_1,n_2,n_3$, respectively, such
that
\begin{equation} \label{eq:HPdef}
  p_{n_1,n_2,n_3}(z)  + q_{n_1,n_2,n_3}(z) e^{z}  + r_{n_1,n_2,n_3}(z)e^{2z}
   = \O(z^{n_1+n_2+n_3+2}), \qquad z \to 0.
\end{equation}
These polynomials exist and are unique up to a common multiplicative
factor. We put
$$
    e_{n_1,n_2,n_3}(z) =
    p_{n_1,n_2,n_3}(z)e^{-z} + q_{n_1,n_2,n_3}(z) +
    r_{n_1,n_2,n_3}(z)e^{z},
$$
so that
$$
   e_{n_1,n_2,n_3}(z)  = \O(z^{n_1+n_2+n_3+2}), \qquad z \to 0.
$$
Putting the left-hand side of (\ref{eq:HPdef}) equal to zero, we obtain
the quadratic algebraic Hermite--Pad\'e approximants
\[ \frac{-q_{n_1,n_2,n_3}(z) \pm \sqrt{q_{n_1,n_2,n_3}(z)^2 -
        4p_{n_1,n_2,n_3}(z)r_{n_1,n_2,n_3}(z)}}{2r_{n_1,n_2,n_3}(z)}
\]
to $e^{z}$.

The analytic theory of the Hermite--Pad\'e approximation to the
exponential function was investigated, among others, by
Mahler \cite{MAH}, Aptekarev \cite{APT}, Chudnovsky \cite{CHU},
Borwein \cite{BOR}, Driver \cite{DRI}, Driver and Temme \cite{DRITEM},
and Wielonsky \cite{wiel1,wiel2}. The paper \cite{NUT} by Nuttall deserves a
special mention since it contains many deep ideas and conjectures about the
asymptotics of Hermite--Pad\'e approximants to general functions. Several
surveys on the subject also exist, see
De Bruin \cite{BRU}, Baker and Lubinsky \cite{BAK-LUB},
and Aptekarev and Stahl \cite{APT-STA}.

\begin{figure}[tb]
\begin{center}
\includegraphics[width=10cm]{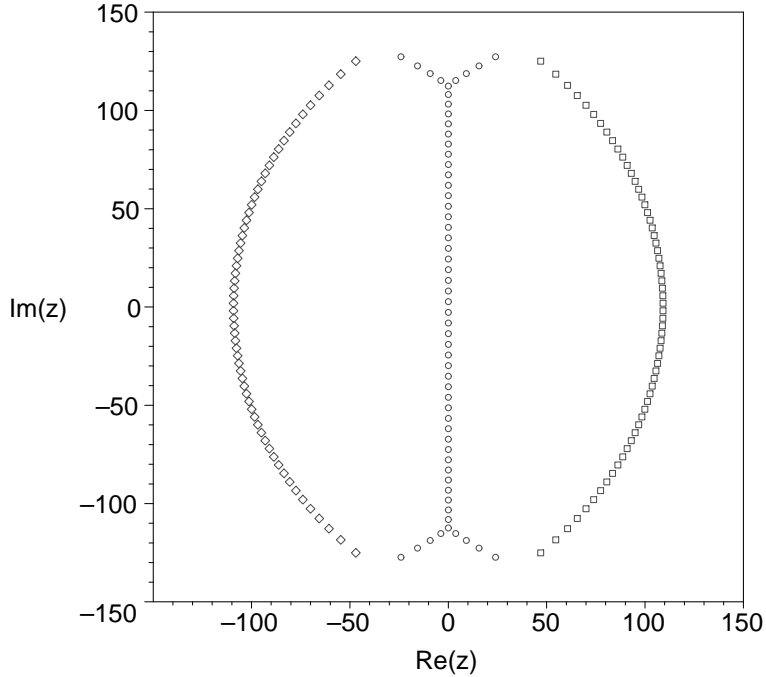}
\caption{Zeros of the diagonal quadratic Hermite--Pad\'e polynomials
$p_{60,60,60}$ (the diamonds on the left) $q_{60,60,60}$ (the circles in
  the middle),
and $r_{60,60,60}$ (the boxes on the right)}
\label{zeroHP60}
\end{center}
\end{figure}

Very recently, Herbert Stahl \cite{stahl1,stahl2,stahl3} gave
detailed asymptotic results for the scaled diagonal polynomials
\begin{equation}
\label{PQR}
P_{n}(z) =  p_{n,n,n}(3nz), \qquad
   Q_{n}(z) =  q_{n,n,n}(3nz), \qquad
   R_{n}(z) =  r_{n,n,n}(3nz),
\end{equation}
and for the remainder term
\begin{equation}
\label{En}
E_n(z)= P_n(z) e^{-3nz} + Q_n(z) + R_n(z) e^{3nz}
\end{equation}
in the approximation.
From these asymptotics, the limit distributions of the zeros of
$P_{n}, Q_{n}, R_{n}$ and $E_n$ follow. The zeros accumulate on
specific arcs in the complex plane. In Figure \ref{zeroHP60}, the zeros of
$p_{n,n,n}$, $q_{n,n,n}$, and $r_{n,n,n}$ have been plotted for the
value $n=60$.
It shows their remarkable distribution, especially the zeros of $q_{n,n,n}$
distribute themselves
in a very particular way. These results on the distribution of the
zeros of Hermite--Pad\'e approximants to exponentials may be seen as a
continuation of the study initiated by Szeg\H o in \cite{SZE}
concerning the distribution of the zeros of Taylor sections of the
series for $e^z$, subsequently generalized by Saff and Varga in
\cite{SAF-VAR1, SAF-VAR2, SAF-VAR3} to the zeros of the Pad\'e
approximants to $e^z$, see also \cite{KuMc,VAR-CAR}.

Stahl \cite{stahl1,stahl2,stahl3} obtained his results from a saddle point analysis for
the explicit integral formulas for $P_n$, $Q_n$, and $R_n$, which are
\begin{eqnarray}
\label{rep-intP}
P_n(z) & = & \frac{C e^{3nz}}{2\pi i}\oint_{C_{-1}}\frac{e^{3nzw}dw}{[w(w^2-1)]^{n+1}}, \\[10pt]
\label{rep-intQ}
Q_n(z) & = & \frac{C}{2\pi i}\oint_{C_0} \frac{e^{3nzw}dw}{[w(w^2-1)]^{n+1}}, \\[10pt]
\label{rep-intR}
R_n(z) & = & \frac{C e^{-3nz}}{2\pi i}\oint_{C_1} \frac{e^{3nzw}dw}{[w(w^2-1)]^{n+1}}.
\end{eqnarray}
Here $C_j$ is a closed contour in the complex plane encircling $j$ in the positive
direction, which does not encircle the other points in $\{ -1, 0, 1\}$. The number
$C$ in (\ref{rep-intP})--(\ref{rep-intR}) is a normalization constant. Stahl used
$C = n! 2^{n+1}(3n)^{-n}$ in order to make $P_n$ a monic polynomial. For our
analysis below, we found it convenient to make $Q_n$ monic, and we choose
\begin{equation} \label{eq:constant}
    C = \frac{n! (-1)^{n+1}}{(3n)^n}.
\end{equation}

In this paper we propose a different approach to the asymptotic analysis.
Our approach is based on a Riemann--Hilbert formulation for the polynomials
$P_n$, $Q_n$, and $R_n$, combined with a steepest descent analysis for Riemann--Hilbert
problems. This technique originated with Deift and Zhou \cite{deiftzhou}
and was applied to the asymptotics of orthogonal polynomials by Deift et al.\
\cite{deift,deiftal1,deiftal2}. See also
\cite{bleher,bleher2,BKMM,KrMc,Kuij,KuMc,KuMc2,KMVV} for recent developments.
The Riemann--Hilbert problem for multiple orthogonal polynomials
was given by Van Assche et al.\ \cite{VAGeKu}.
In the present situation it gives rise to a Riemann--Hilbert problem for a
$3 \times 3$-matrix valued function.
As far as we are aware, this is the first time that the steepest descent method
for Riemann--Hilbert problems is applied to a $3 \times 3$ problem.

The Riemann--Hilbert problem is to find a $3\times 3$ matrix valued function
$Y: \mathbb{C}\setminus \Gamma \to \mathbb{C}^{3\times 3}$
where $\Gamma$ is  a closed contour in the complex plane encircling the
origin once in the positive direction, such that
\begin{enumerate}
\item $Y$ is analytic in $\mathbb{C} \setminus \Gamma$.
\item $Y$ satisfies the jump condition
\begin{equation}  \label{eq:Yjump}
   Y_+(z) = Y_-(z) \begin{pmatrix}
                   1 & z^{-3n-2} e^{-3nz} & 0 \\
                   0 & 1 & 0 \\
                   0 & z^{-3n-2}e^{3nz} & 1
                   \end{pmatrix}, \qquad z \in \Gamma,
\end{equation}
where $Y_+(z)$ and $Y_-(z)$ denote the limiting values of $Y(z')$ as $z'$
approaches $z \in \Gamma$ from the inside and outside of $\Gamma$, respectively.
\item For large $z$
\begin{equation}  \label{eq:Yasym}
    Y(z) = \left( I + \O\left(\frac1z\right) \right)
           \begin{pmatrix}
            z^{n+1} & 0 & 0 \\
            0 & z^{-2n-2} & 0 \\
            0 & 0 & z^{n+1}
            \end{pmatrix}, \qquad z \to \infty.
\end{equation}
\end{enumerate}
We will show in Section 5.1 that the Riemann--Hilbert problem has
a unique solution and that $Y_{21}(z) = P_n(z)$,
$Y_{22}(z) =  z^{-3n-2} Q_n(z)$ for $z$ outside $\Gamma$,
$Y_{22}(z) = z^{-3n-2} E_n(z)$ for $z$ inside $\Gamma$,
and $Y_{23}(z) = R_n(z)$.
\medskip

The steepest descent analysis consists of a number of transformations.
To make the transformations work, we make heavy use of the works of Stahl
\cite{stahl1,stahl2,stahl3}. A crucial role is played  by the Riemann
surface defined by
\begin{equation} \label{z=zw}
    z= \frac{w^2- \frac{1}{3}}{w(w^2-1)},
\end{equation}
which is considered as a three sheeted surface with cuts along
two arcs $\Gamma_P$ and $\Gamma_R$.
The jumps of the different inverse mappings of (\ref{z=zw}) across
the arcs determine probability measures $\mu_P$ and $\mu_R$
supported on $\Gamma_P$ and $\Gamma_R$, respectively.
These measures turn out to be limiting distributions of the normalized
zero counting measures of the polynomials $P_n$ and $R_n$, respectively.

We choose the closed contour $\Gamma$ in the Riemann--Hilbert problem
for $Y$ so that it contains the cuts $\Gamma_P$ and $\Gamma_R$.
The measures $\mu_P$ and $\mu_R$ and their $g$-transforms
\[ g_P(z) = \int \log(z-s) d\mu_P(s),
    \qquad g_R(z) = \int \log(z-s) d\mu_R(s) \]
are used to make the first transformation of the Riemann--Hilbert problem,
which has the effect of normalizing the problem at infinity.
Then we follow the general scheme, as presented in
\cite{deiftal2} or \cite{deift}, for the asymptotic analysis of
Riemann--Hilbert problems. It leads to a final
Riemann--Hilbert problem whose solution has an explicit asymptotic
behavior for $n \to \infty$, see Theorem \ref{thm:T} in Section 6.3.
Tracing our steps back to the original Riemann--Hilbert
problem, we obtain strong asymptotic formulas for the
Hermite--Pad\'e approximants in every region of the complex plane,
including the sets where the zeros are and their endpoints.

In Section 2 we state the asymptotic results for the
polynomials $P_n$, $Q_n$, and $R_n$, and for the remainder $E_n$,
as well as for the quadratic Hermite--Pad\'e
approximants. The asymptotic results make use of the functions
obtained from the Riemann surface. The Riemann surface and the measures
and functions derived from it are also described in Section 2.
In Section 3 we prove the statements about the Riemann surface and
other geometrical objects involved in the problem.
To prepare for the transformations of the Riemann--Hilbert problem
we need a large number of relations between the various functions
involved, such as the inverse mappings of (\ref{z=zw}) and the
functions $g_P$ and $g_R$. These properties are established in
terms in Section 4. Sections 5 and 6 contain the transformations
of the Riemann--Hilbert problem and all asserted asymptotic
results are finally proven in Section 7.

\section{Statement of results}

\subsection{The Riemann surface}
In order to state our results we first introduce an
appropriate Riemann surface. The Riemann surface is
chosen so that it parameterizes the critical points of the function
\begin{equation} \label{def-integrand}
    w \mapsto 3zw-\log \left[w(w^2-1) \right].
\end{equation}
Note that the integrals in formulas (\ref{rep-intP})--(\ref{rep-intR})
have the form
\begin{equation} \label{integral}
    \oint_{C_j} \frac{1}{w(w^2-1)} e^{n\left(3zw-\log \left[w(w^2-1)\right]\right)} dw
\end{equation}
and that by the classical saddle point analysis for the asymptotic evaluation
of integrals, the main contribution to the integral (\ref{integral}) comes from
a critical point of (\ref{def-integrand}).
So we define $\mathcal{R}$ as the Riemann surface for the function
\begin{equation}  \label{eq:zw}
    z = z(w)
     = \frac13 \left( \frac1w + \frac{1}{w-1}+\frac{1}{w+1} \right)
     = \frac{w^2-\frac13}{w(w^2-1)}.
\end{equation}
Note that we obtain (\ref{eq:zw}) if we set the derivative of (\ref{def-integrand})
equal to zero and solve for $z$.

\unitlength=1.00mm
\begin{figure}[ht]
\centering
\begin{picture}(140.00,80.00)(0,60)
\put(30.00,140.00){\line(1,0){110.00}}
\put(140.00,140.00){\line(-3,-2){30.00}}
\put(110.00,120.00){\line(-1,0){110.00}}
\put(0.00,120.00){\line(3,2){30.00}}
\put(30.00,110.00){\line(1,0){110.00}}
\put(140.00,110.00){\line(-3,-2){30.00}}
\put(110.00,90.00){\line(-1,0){110.00}}
\put(0.00,90.00){\line(3,2){30.00}}
\put(30.00,80.00){\line(1,0){110.00}}
\put(140.00,80.00){\line(-3,-2){30.00}}
\put(110.00,60.00){\line(-1,0){110.00}}
\put(0.00,60.00){\line(3,2){30.00}}
\bezier{60}(80.00,65.00)(87.00,69.00)(90.00,75.00)
\bezier{60}(50.00,125.00)(54.00,132.00)(60.00,135.00)
\bezier{60}(50.00,95.00)(54.00,102.00)(60.00,105.00)
\bezier{60}(80.00,95.00)(87.00,99.00)(90.00,105.00)
\put(50.00,125.00){\line(0,-1){3.00}}
\put(50.00,120.00){\line(0,-1){3.00}}
\put(50.00,115.00){\line(0,-1){3.00}}
\put(50.00,110.00){\line(0,-1){3.00}}
\put(50.00,105.00){\line(0,-1){3.00}}
\put(50.00,100.00){\line(0,-1){3.00}}
\put(50.00,125.00){\circle*{2.00}}
\put(50.00,95.00){\circle*{2.00}}
\put(60.00,135.00){\line(0,-1){3.00}}
\put(60.00,130.00){\line(0,-1){3.00}}
\put(60.00,125.00){\line(0,-1){3.00}}
\put(60.00,120.00){\line(0,-1){3.00}}
\put(60.00,115.00){\line(0,-1){3.00}}
\put(60.00,110.00){\line(0,-1){3.00}}
\put(60.00,135.00){\circle*{2.00}}
\put(60.00,105.00){\circle*{2.00}}
\put(80.00,95.00){\line(0,-1){3.00}}
\put(80.00,90.00){\line(0,-1){3.00}}
\put(80.00,85.00){\line(0,-1){3.00}}
\put(80.00,80.00){\line(0,-1){3.00}}
\put(80.00,75.00){\line(0,-1){3.00}}
\put(80.00,70.00){\line(0,-1){3.00}}
\put(80.00,95.00){\circle*{2.00}}
\put(80.00,65.00){\circle*{2.00}}
\put(90.00,105.00){\line(0,-1){3.00}}
\put(90.00,100.00){\line(0,-1){3.00}}
\put(90.00,95.00){\line(0,-1){3.00}}
\put(90.00,90.00){\line(0,-1){3.00}}
\put(90.00,85.00){\line(0,-1){3.00}}
\put(90.00,80.00){\line(0,-1){3.00}}
\put(90.00,105.00){\circle*{2.00}}
\put(90.00,75.00){\circle*{2.00}}
\put(56.00,138.00){\makebox(0,0)[cc]{$z_1$}}
\put(46.00,128.00){\makebox(0,0)[cc]{$z_2$}}
\put(84.00,63.00){\makebox(0,0)[cc]{$z_3$}}
\put(94.00,73.00){\makebox(0,0)[cc]{$z_4$}}
\put(84.00,73.00){\makebox(0,0)[cc]{$\Gamma_R$}}
\put(55.00,127.00){\makebox(0,0)[cc]{$\Gamma_P$}}
\put(127.00,126.00){\makebox(0,0)[cc]{$\mathcal{R}_P$}}
\put(127.00,96.00){\makebox(0,0)[cc]{$\mathcal{R}_Q$}}
\put(127.00,66.00){\makebox(0,0)[cc]{$\mathcal{R}_R$}}
\end{picture}
\caption{The Riemann surface $\mathcal{R}$}
\label{fig:surface}
\end{figure}
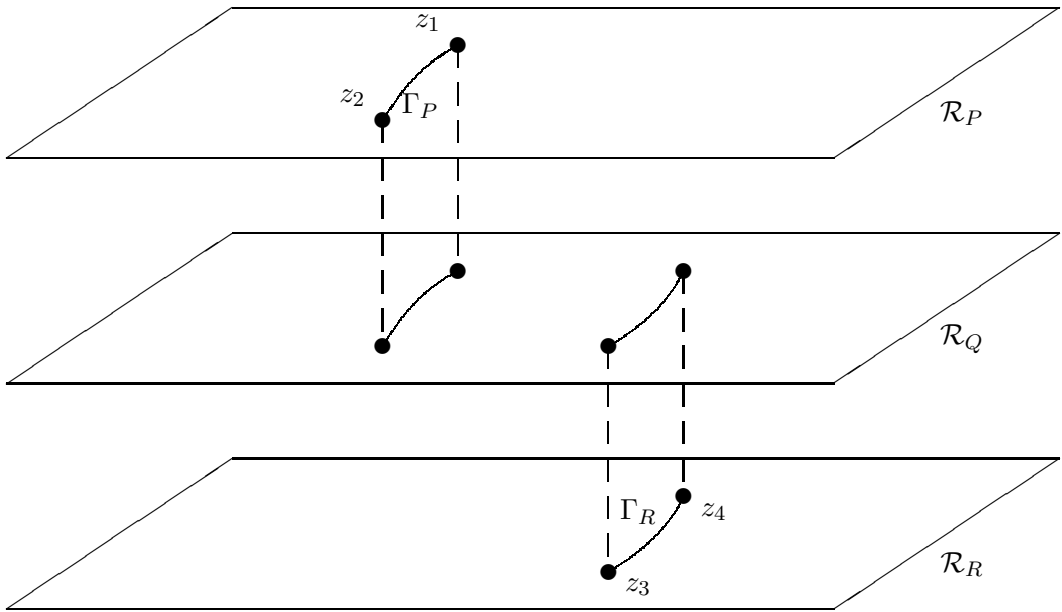

The rational function (\ref{eq:zw}) has three inverse mappings. These are
the three solutions of the cubic equation
\begin{equation}  \label{eq:cubic}
   zw^3-w^2-zw+\frac{1}{3} = 0.
\end{equation}
The Riemann surface $\mathcal{R}$ consists of three sheets $\mathcal{R}_P$,
$\mathcal{R}_Q$, and $\mathcal{R}_R$ as shown in Figure \ref{fig:surface},
see Proposition \ref{prop:GammaPR}. The
bijective mapping
$\psi : \mathcal{R} \to \overline{\mathbb{C}}$ is the inverse of (\ref{eq:zw}).
We denote its restriction to the three sheets by $\psi_P$, $\psi_Q$,
and $\psi_R$, respectively. So $\psi_P(z)$, $\psi_Q(z)$, and $\psi_R(z)$
are the three solutions of (\ref{eq:cubic}). Typically we will identify
the three sheets with copies of the complex plane, and so $\psi_P$, $\psi_Q$,
and $\psi_R$ are defined on $\mathbb C$ with appropriate cuts.
The sheets $\mathcal{R}_P$ and $\mathcal{R}_Q$ are glued together along a cut $\Gamma_P$
connecting two branch points $z_1$ and $z_2$, and the sheets $\mathcal{R}_Q$ and
$\mathcal{R}_R$ are glued together along a cut $\Gamma_R$ connecting the
other two branch points $z_3$ and $z_4$.

The Riemann surface has four branch points
$z_1=z(w_1)$, $z_2=z(w_2)$, $z_3=z(w_3)$, $z_4=z(w_4)$, which
are related to the points $w_1$, $w_2$, $w_3$, $w_4$ for which $z'(w) = 0$.
Simple calculations give
\begin{equation}  \label{eq:wk}
  w_k = 3^{-1/4} \omega_8^{-2k-1}, \qquad k=1,2,3,4,
\end{equation}
where $\omega_8 = e^{2\pi i/8}$  is the primitive 8th root of
unity. The corresponding values of $z_k = z(w_k)$ are
\begin{equation} \label{eq:zk}
z_1 = 3^{-1/4} \omega_{24}^{7},\quad
z_2 = 3^{-1/4} \omega_{24}^{17},\quad
z_3 = 3^{-1/4} \omega_{24}^{19},\quad
z_4 = 3^{-1/4} \omega_{24}^{5},
\end{equation}
where $\omega_{24} = e^{2\pi i/24}$ is the primitive 24th root of unity.
The precise sheet structure of $\mathcal{R}$ is given in the following proposition.

\begin{proposition} \label{prop:GammaPR}
There is an analytic curve $\Gamma_P$ from $z_1$ to $z_2$ lying in the left half-plane,
and an analytic curve $\Gamma_R$ from $z_3$ to $z_4$ lying in the right half-plane, such
that the following hold.
\begin{enumerate}
\item[\rm (a)] Three inverse mappings $\psi_P$, $\psi_Q$, and $\psi_R$ of
{\rm(\ref{eq:zw})}
exist so that $\psi_P$ is defined and analytic on $\C \setminus \Gamma_P$,
$\psi_Q$ is defined and analytic on $\C \setminus (\Gamma_P \cup \Gamma_R)$, and
$\psi_R$ is defined and analytic on $\C \setminus \Gamma_R$.
\item[\rm (b)] At infinity, we have the values
$\psi_P(\infty) = -1$, $\psi_Q(\infty) = 0$, and $\psi_R(\infty) = 1$.
\item[\rm (c)] For $z \in \Gamma_P$, we have
\begin{equation} \label{eq:GammaP}
    \frac{3}{2\pi i} \int_{z_1}^z (\psi_Q-\psi_P)_+(s) ds \in \mathbb R,
\end{equation}
with integration along the $+$side of $\Gamma_P$. [This is the side of $\Gamma_P$
that is on the left while going from $z_1$ to $z_2$ along $\Gamma_P$.]
\item[\rm (d)] For $z \in \Gamma_R$, we have
\begin{equation}  \label{eq:GammaR}
  \frac{3}{2\pi i} \int_{z_3}^z (\psi_Q-\psi_R)_+(s)\,ds \in \mathbb{R},
\end{equation}
with integration along the $+$side of $\Gamma_R$.
\end{enumerate}
\end{proposition}

\begin{figure}[htb]
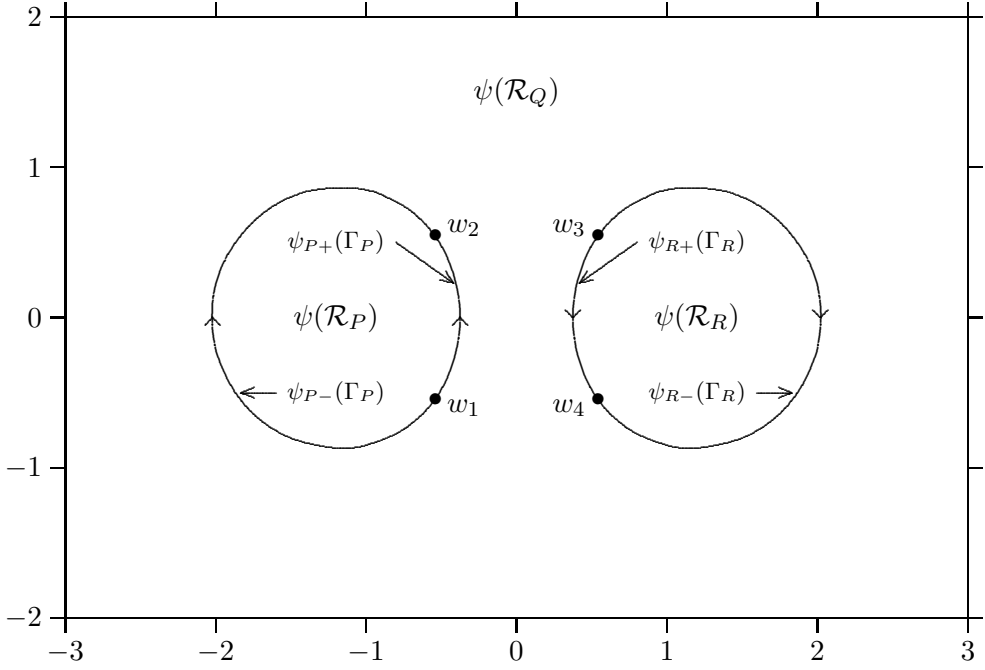

\hfil
\input traject8
\hfil
\caption{$\psi$-image of the Riemann surface $\mathcal{R}$}  \label{fig:psi}
\end{figure}

Because of symmetry, $\Gamma_R$ is the mirror image of $\Gamma_P$
under reflection with respect to the imaginary axis.
We take the sheet structure on the Riemann surface $\mathcal R$ as in
Figure \ref{fig:surface}.
The functions $\psi_P$, $\psi_Q$, $\psi_R$ are defined on the $P$, $Q$,
and $R$ sheet of $\mathcal R$, respectively. Together they constitute
a conformal map from $\mathcal R$ onto the Riemann sphere. The images of
the different sheets are shown in Figure \ref{fig:psi}.

\begin{figure}[htb]
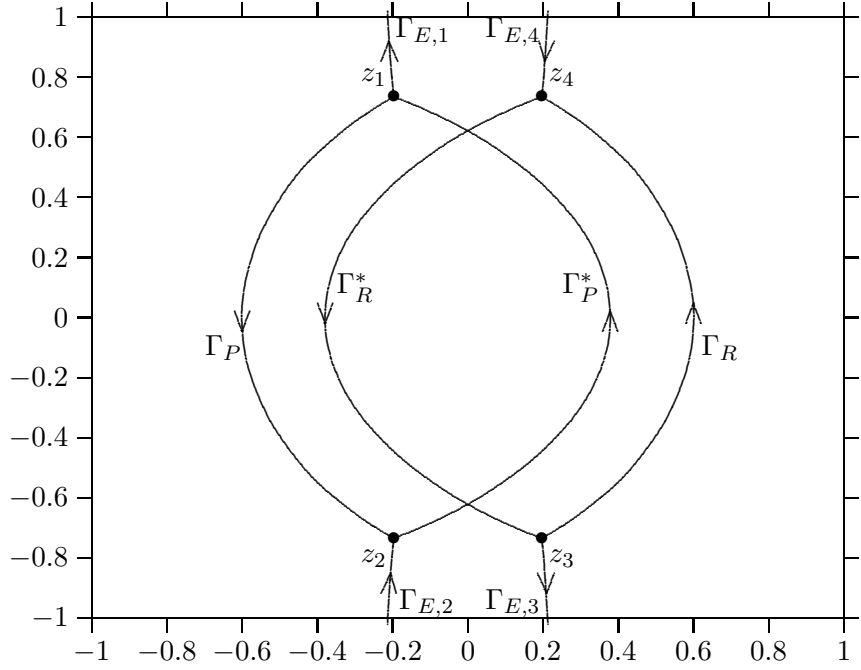

\hfil
\input traject
\hfil
\caption{Curves for which either $\frac{3}{2\pi i} \int_{z_1}^z (\psi_Q-\psi_P)(s)\,ds$
or $\frac{3}{2\pi i} \int_{z_3}^z (\psi_Q-\psi_R)(s)\,ds$ is real.}
 \label{fig:traject}
\end{figure}

Other curves of interest to our problem are defined by the property that either
$\frac{3}{2\pi i} \int_{z_1}^z (\psi_Q-\psi_P)(s) ds$
or $\frac{3}{2\pi i} \int_{z_3}^z (\psi_Q-\psi_R)(s) ds$
is real. These curves are described by the following proposition,
see Figure \ref{fig:traject}.

\begin{proposition} \label{prop:GammaP*R*}
\begin{enumerate}
\item[\rm (a)] There are four analytic curves where
$\frac{3}{2\pi i} \int_{z_1}^z (\psi_Q-\psi_P)(s)\,ds$ is real.
One of them is $\Gamma_P$, a second one joins $z_2$ with $z_1$ and
intersects the positive real axis. We call this curve $\Gamma_P^*$.
The other two curves are unbounded and lie in the left half-plane.
They join $z_1$ and $z_2$ with infinity, and we call
them $\Gamma_{E,1}$ and $\Gamma_{E,2}$, respectively.
\item[\rm (b)] There are four analytic curves where
$\frac{3}{2\pi i} \int_{z_3}^z (\psi_Q - \psi_R)(s)\,ds$ is real.
One of them is $\Gamma_R$, a second one joins $z_4$ with $z_3$
and intersects the negative real axis. We call this curve $\Gamma_R^*$.
The other two curves are unbounded and lie in the right half-plane.
They join $z_3$ and $z_4$
with infinity and we call them $\Gamma_{E,3}$ and $\Gamma_{E,4}$, respectively.
\item[\rm (c)] The curves $\Gamma_P^*$ and $\Gamma_R$ do not intersect.
\item[\rm (d)] The curves $\Gamma_R^*$ and $\Gamma_P$ do not intersect.
\item[\rm (e)] The curves $\Gamma_P^*$ and $\Gamma_R^*$ intersect in the points $\pm i y^*$
on the imaginary axis. The value of $y^*$ is approximately $y^* = 0.621391\cdots$.
\end{enumerate}
\end{proposition}

By symmetry, $\Gamma_P^*$ and $\Gamma_R^*$ are mirror images with respect
to the imaginary axis. Similarly, $\Gamma_{E,1}$ and $\Gamma_{E,4}$ are mirror
images of each other, as well as $\Gamma_{E,2}$ and $\Gamma_{E,3}$.
All contours are oriented as shown in Figure \ref{fig:traject}. The orientation
induces a $+$side and a $-$side on each contour, where the $+$side is on the
left and the $-$side on the right while traversing the contour according to its orientation.
Propositions \ref{prop:GammaPR} and \ref{prop:GammaP*R*} are proved in Section 3.
\medskip

\begin{figure}[ht]
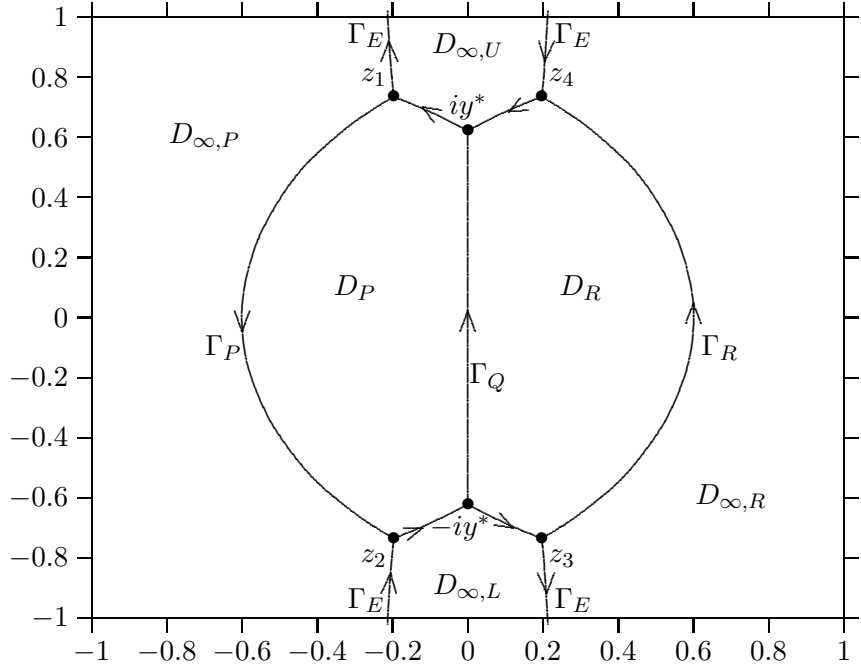

\hfil
\input trajecta
\hfil
\caption{Curves $\Gamma_P$, $\Gamma_Q$, $\Gamma_R$, and $\Gamma_E$, and domains
$D_P$, $D_R$, $D_{\infty,P}$,$D_{\infty,R}$,$D_{\infty,U}$, and $D_{\infty,L}$.
Note that $\Gamma_Q$ does not only consist of the interval $[-iy^*,iy^*]$ but also
of the small arcs from $z_2$ via $-iy^*$ to $z_3$ and from $z_4$ via $iy^*$ to $z_1$.
The domain $D_{\infty}$ is the union of $\Gamma_E$ with the unbounded domains $D_{\infty,P}$, $D_{\infty,R}$,
$D_{\infty,U}$, and $D_{\infty,L}$.}
 \label{fig:trajecta}
\end{figure}

Now we can also define the contours $\Gamma_Q$ and $\Gamma_E$.
We define
\begin{equation} \label{eq:GammaE}
     \Gamma_E = \bigcup_{k=1}^4 \Gamma_{E,k}.
\end{equation}
where $\Gamma_{E,k}$, $k=1,2,3,4$ are as in Proposition \ref{prop:GammaP*R*},
and
\begin{equation} \label{eq:GammaQ}
    \Gamma_Q = [-iy^*, iy^*] \cup \left(\Gamma_P^* \cap \{ \Re z < 0 \} \right)
    \cup \left( \Gamma_R^* \cap \{ \Re z > 0 \} \right),
\end{equation}
where $\Gamma_P^*$, $\Gamma_R^*$ and $iy^*$ are as in Proposition \ref{prop:GammaP*R*}.
We choose an orientation on $\Gamma_Q$ by orienting the interval $[-iy^*,iy^*]$
from $-iy^*$ to $iy^*$ and by orienting $\Gamma_P^*$ from $z_2$ to $z_1$
and $\Gamma_R^*$ from $z_4$ to $z_3$, see Figure \ref{fig:trajecta}.

The contours $\Gamma_P$, $\Gamma_Q$, $\Gamma_R$, and $\Gamma_E$ divide the complex plane
into six domains. We denote the unbounded domains by $D_{\infty,P}$, $D_{\infty,R}$,
$D_{\infty,U}$ and $D_{\infty,L}$, as shown in Figure \ref{fig:trajecta}. The
bounded domains are denoted by $D_P$ and $D_R$, where $D_P$ is in the left half-plane,
and $D_R$ is in the right half-plane, see also Figure \ref{fig:trajecta}.
We also put
\begin{equation} \label{eq:defDinfty}
    D_{\infty} = \left(D_{\infty,P} \cup D_{\infty,R} \cup D_{\infty,U} \cup D_{\infty,L} \right)
    \cup \left(\Gamma_E \setminus \{z_1,z_2,z_3,z_4\}\right).
\end{equation}
This is the unbounded domain bounded by $\Gamma_P$, $\Gamma_Q$ and $\Gamma_R$.
Of further interest is the bounded domain $D_P^*$ bounded by $\Gamma_P$ and $\Gamma_P^*$,
see Figure \ref{fig:Re}, and its mirror image $D_R^*$ with respect to the imaginary axis.

This concludes the description of the Riemann surface.

\subsection{The measures $\mu_P$, $\mu_Q$, $\mu_R$, and $\mu_E$}
We now define a measure on each of the curves $\Gamma_P$, $\Gamma_Q$, $\Gamma_R$,
and $\Gamma_E$. The complex line element $ds$ is defined according to
the orientation of these curves given in Figure \ref{fig:trajecta}.

\begin{definition}
We define a measure $\mu_P$ on $\Gamma_P$ by
\begin{equation} \label{eq:muP}
    d\mu_P(s) = \frac{3}{2\pi i} (\psi_Q-\psi_P)_+(s)\, ds \qquad s \in \Gamma_P,
\end{equation}
a measure $\mu_R$ on $\Gamma_R$ by
\begin{equation} \label{eq:muR}
  d\mu_R(s) = \frac{3}{2\pi i} (\psi_Q-\psi_R)_+(s)\, ds \qquad s \in \Gamma_R,
\end{equation}
and a measure $\mu_Q$ on $\Gamma_Q$ by
\begin{equation} \label{eq:muQ}
    d\mu_Q(s) = \left\{ \begin{array}{l}
        \frac{3}{2\pi i} (\psi_Q-\psi_P)(s)\, ds \qquad s \in \Gamma_P^* \cap \{ \Re z < 0 \}, \\[10pt]
        \frac{3}{2\pi i} (\psi_Q-\psi_R)(s)\, ds \qquad s \in \Gamma_R^* \cap \{ \Re z > 0 \}, \\[10pt]
        \frac{3}{2\pi i} (\psi_R-\psi_P)(s)\, ds \qquad s \in [-iy^*, iy^*].
        \end{array} \right.
\end{equation}
The measure $\mu_E$ on $\Gamma_E$ is defined as
\begin{equation}  \label{eq:muE}
  d\mu_E(s) = \left\{ \begin{array}{l}
       \frac{3}{2\pi i} (\psi_Q-\psi_P)(s)\, ds \qquad s \in \Gamma_{E,1} \cup \Gamma_{E,2}, \\[10pt]
       \frac{3}{2\pi i} (\psi_Q-\psi_R)(s)\, ds \qquad s \in \Gamma_{E,3} \cup \Gamma_{E,4}.
          \end{array} \right.
\end{equation}
\end{definition}
A priori these are complex measures. Our first result is that $\mu_P$, $\mu_Q$,
and $\mu_R$ are in fact probability measures and $\mu_E$ is a positive measure.

\begin{theorem}
\label{mes-proba}
We have that $\mu_P$ is a probability measure on $\Gamma_P$, $\mu_Q$ is
a probability measure on $\Gamma_Q$ and $\mu_R$ is a probability measure
on $\Gamma_R$. The measure $\mu_E$ is a positive measure on $\Gamma_E$.
\end{theorem}
Theorem \ref{mes-proba} follows from Propositions \ref{prop:muPmuR} and
\ref{prop:muQ}.
\medskip

The relevance of these measures is shown by the following theorem. For
every polynomial $p$ of exact degree $n$, we denote by $\nu_{p}$ the normalized
zero counting measure. Thus
\[ \nu_{p} = \frac{1}{n} \sum_{p(z) = 0} \delta_{z} \]
where each zero is counted according to its multiplicity. We also define a zero
counting measure for the remainder function $E_n$, namely
\[ \nu_{E_n} = \frac{1}{n} \sum_{\stackrel{E_n(z) = 0}{z\neq 0}} \delta_{z}, \]
where
the normalization by $n$ now corresponds to the degree of approximation
and the $3n+2$ interpolatory zeros of $E_n$ at the origin have been excluded.
\begin{theorem}
\label{lim-zeros}
We have
\begin{equation} \label{eq:convnuPQRn}
    \nu_{P_n} \stackrel{*}{\to} \mu_P, \qquad
   \nu_{Q_n} \stackrel{*}{\to} \mu_Q, \qquad
   \nu_{R_n} \stackrel{*}{\to} \mu_R,
\end{equation}
where the convergence is in the sense of weak$^*$ convergence of measures,
i.e., $\mu_n \stackrel{*}{\to} \mu$ if $\int f d\mu_n \to \int fd\mu$ for
every bounded continuous function $f$.
Furthermore, we have
\begin{equation} \label{eq:convnuEn}
    \nu_{E_n} \to \mu_{E},
\end{equation}
in the sense that
\[ \lim_{n\to\infty} \int f(s) d\nu_{E_n}(s) = \int f(s) d\mu_E(s) \]
for every continuous function $f$ such that $f(s) = \O(s^{-2})$
as $s \to \infty$.
\end{theorem}

The convergence of the zero counting measures is due to Stahl \cite{stahl2}.
The theorem shows that the measures $\mu_P$, $\mu_Q$, $\mu_R$, and $\mu_E$ agree
with the measures that Stahl introduced in a different way.

In contrast to the measures $\nu_{P_n}$, $\nu_{Q_n}$, and $\nu_{R_n}$ which are
probability measures, the measures $\nu_{E_n}$ have infinite mass. They also have
unbounded support. As a result, the proof of the limit (\ref{eq:convnuEn}) is more
involved than that of (\ref{eq:convnuPQRn}).

\subsection{The $g$-functions}
For the strong asymptotic results we need the log-transforms
(or complex logarithmic potentials) of the measures $\mu_P$,
$\mu_Q$, and $\mu_R$.
\begin{definition}
We introduce three functions
\begin{eqnarray}
  g_P(z) & = & \int_{\Gamma_P} \log (z-s) \, d\mu_P(s),
  \qquad z \in \mathbb C \setminus \Gamma_P, \label{eq:gP} \\
  g_Q(z) & = & \int_{\Gamma_Q} \log(z-s) \, d\mu_Q(s),
  \qquad z \in \mathbb C \setminus \Gamma_Q, \label{eq:gQ} \\
  g_R(z) & = & \int_{\Gamma_R} \log (z-s) \, d\mu_R(s),
  \qquad z \in \mathbb C \setminus \Gamma_R, \label{eq:gR}
\end{eqnarray}
which are defined modulo $2\pi i$.
\end{definition}
Thus $g_P$, $g_Q$, and $g_R$ are multivalued functions,
depending on the specific choice of the branches of the logarithmic functions.
Our results will involve expressions like $e^{n g_P}$ and $e^{n g_R}$, and
then the multivaluedness will play no role.

\subsection{Functions $\varphi_P$ and $\varphi_R$}
Two other important functions are the functions $\varphi_P$ and $\varphi_R$ given by
\begin{equation}  \label{eq:varphiP}
   \varphi_P(z) = \frac{3}{2} \int_{z_1}^z (\psi_Q-\psi_P)(s)\, ds,
\end{equation}
and
\begin{equation} \label{eq:varphiR}
   \varphi_R(z) = \frac{3}{2} \int_{z_3}^z (\psi_Q-\psi_R)(s)\, ds.
\end{equation}
The paths of integration in (\ref{eq:varphiP}) and (\ref{eq:varphiR})
are in $\mathbb C \setminus \left(\Gamma_P \cup \Gamma_R \cup \{0\}\right)$.
The functions $\varphi_P$ and $\varphi_R$ are multivalued but the real parts
are well-defined.
From Proposition \ref{prop:GammaP*R*} we know that $\Re \varphi_P = 0$
on the curves $\Gamma_P$, $\Gamma_P^*$, $\Gamma_{E,1}$, and $\Gamma_{E,2}$,
and that $\Re \varphi_R = 0$ on $\Gamma_R$, $\Gamma_R^*$, $\Gamma_{E,1}$,
and $\Gamma_{E,2}$.
We collect the main properties of $\varphi_P$ and $\varphi_R$
in the following lemma.
\begin{lemma} \label{lemma:phiPR}
\begin{enumerate}
\item[\rm (a)] The real part of $\varphi_P$ is zero exactly on $\Gamma_P$, $\Gamma_P^*$,
$\Gamma_{E,1}$, and $\Gamma_{E,2}$.

The real part of $\varphi_P$ is negative in $D_{\infty,P} \cup
D_P^*$, and it is positive in the remaining part of the plane.
\item[\rm (b)] The real part of $\varphi_R$ is zero exactly on $\Gamma_R$, $\Gamma_R^*$,
$\Gamma_{E,3}$, and $\Gamma_{E,4}$.

The real part of $\varphi_R$ is negative in $D_{\infty,R} \cup D_R^*$, and it is
positive in the remaining part of the plane.
\item[\rm (c)]
On the imaginary axis, we have $\Re \varphi_P = \Re \varphi_R$. We
have $\Re \varphi_P < \Re \varphi_R$ in the left half-plane, and $\Re \varphi_P > \Re \varphi_R$
in the right half-plane.
\end{enumerate}
\end{lemma}
Lemma \ref{lemma:phiPR} is proved in Section \ref{subsec55}.

\begin{figure}[tb]
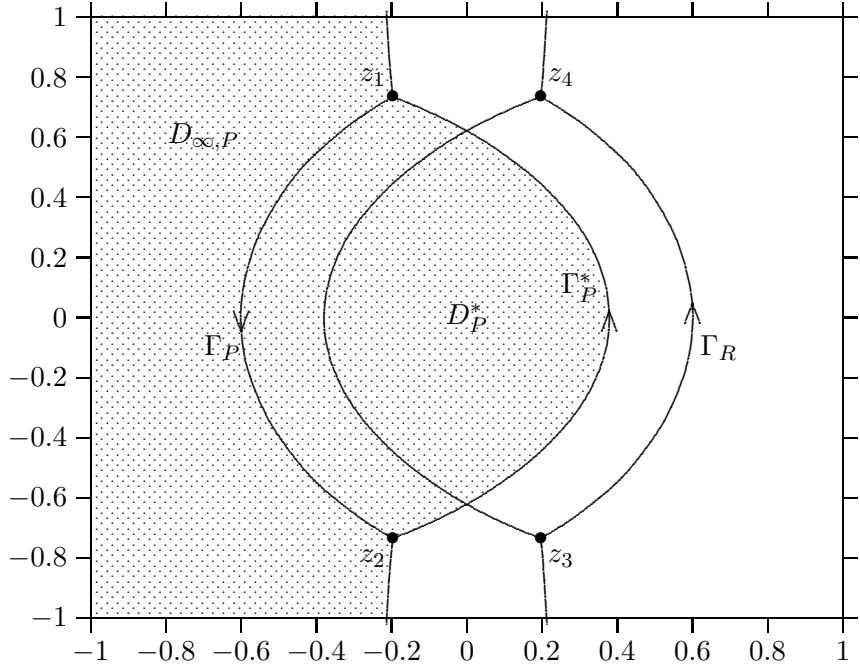

\hfil
\input traject1
\hfil \caption{The shaded region is where $\Re \varphi_P$ is
negative. It consists of the two parts $D_{\infty, P}$ and
$D_P^*$, where $D_P^*$ is bounded by $\Gamma_P$ and $\Gamma_P^*$.}
\label{fig:Re}
\end{figure}

\subsection{Strong asymptotics away from the zeros}

Now we can state the strong asymptotic results for the polynomials $P_n$,
$Q_n$, $R_n$, and the remainder function $E_n$. Recall that the three
polynomials are given by (\ref{rep-intP}),
(\ref{rep-intQ}) and (\ref{rep-intR}) with the constant $C$
as in (\ref{eq:constant}) so that $Q_n$ is a monic polynomial. Then $P_n$
and $R_n$ both have leading coefficients $(-\frac{1}{2})^{n+1}$, and we will
state our asymptotic results for the monic polynomials $(-2)^{n+1} P_n(z)$
and $(-2)^{n+1} R_n(z)$.
Throughout the rest of the paper, we will use the function
$\sqrt{3w^4+1}$ which branches at the four points $w_k$ given in (\ref{eq:wk}).
We choose as cuts for
this function the two curves $\psi_{P+}(\Gamma_P)$ and
$\psi_{R+}(\Gamma_R)$ (see Figure \ref{fig:psi}), and assume that it
is positive for large positive $w$. So, in particular we have
that $\sqrt{3w^4+1} = -1$ for $w=0$.

The following theorem gives the strong asymptotics of the polynomials
$P_n$, $Q_n$, $R_n$ and the remainder term $E_n$ away from their zeros.
These results are due to Stahl \cite{stahl2}, but we will give independent
proofs below.
\begin{theorem}
\label{asymp}
With the functions defined above, we have
\begin{equation} \label{asympP}
(-2)^{n+1} P_n(z) = \frac{2 e^{n g_P(z)}}{\sqrt{3\psi_P^4(z)+1}}
    \left(1+ {\O}\left(\frac{1}{n}\right)\right)
\end{equation}
uniformly for $z$ in compact subsets of $\mathbb C \setminus \Gamma_P$,
\begin{equation} \label{asympR}
(-2)^{n+1} R_n(z) = \frac{2 e^{n g_R(z)}}{\sqrt{3\psi_R^4(z)+1}}
    \left(1+ {\O}\left(\frac{1}{n}\right)\right)
\end{equation}
uniformly for $z$ in compact subsets of $\mathbb C \setminus \Gamma_R$,
and
\begin{equation} \label{asympQ}
Q_n(z) = \left\{ \begin{array}{ll} \ds
    -\frac{e^{n g_Q(z)}}{\sqrt{3\psi_Q^4(z) + 1}}   \left(1 + {\O}\left(\frac{1}{n}\right)\right)
    & \mbox{ for } z \in D_{\infty} \\[10pt] \ds
    \frac{e^{n g_Q(z)}}{\sqrt{3\psi_P^4(z) + 1}}   \left(1 + {\O}\left(\frac{1}{n}\right)\right)
    & \mbox{ for } z \in D_P \\[10pt] \ds
    \frac{e^{n g_Q(z)}}{\sqrt{3\psi_R^4(z) + 1}}  \left(1 + {\O}\left(\frac{1}{n}\right)\right)
    & \mbox{ for } z \in D_R,
    \end{array} \right.
\end{equation}
uniformly for $z$ in compact subsets of $\mathbb C \setminus \Gamma_Q$. Furthermore we have
\begin{equation} \label{asympE}
 E_n(z) = \left\{ \begin{array}{ll} \ds
    - \left( -\frac{1}{2} \right)^n \frac{e^{n (g_R(z)+3z)}}{\sqrt{3\psi_R^4(z) + 1}}   \left(1 + {\O}\left(\frac{1}{n}\right)\right)
    & \mbox{ for } z \in D_{\infty,R} \\[10pt] \ds
    -\left( -\frac{1}{2} \right)^n \frac{e^{n (g_P(z)-3z)}}{\sqrt{3\psi_P^4(z) + 1}}   \left(1 + {\O}\left(\frac{1}{n}\right)\right)
    & \mbox{ for } z \in D_{\infty,P} \\[10pt] \ds
    - \frac{z^{3n} e^{-n(g_P(z)+g_R(z))}}{\sqrt{3\psi_Q^4(z) + 1}}  \left(1 + {\O}\left(\frac{1}{n}\right)\right)
    & \mbox{ for } z \in D_P\cup D_R \cup D_{\infty,U} \cup D_{\infty,L} \\[-10pt]
    & \qquad \quad \cup \; \Gamma_Q \setminus \{z_1, z_2, z_3, z_4\},
    \end{array} \right.
\end{equation}
uniformly for $z$ in compact subsets of $\mathbb{C} \setminus \Gamma_E$.
\end{theorem}

\subsection{Asymptotics near the curves $\Gamma_P$, $\Gamma_Q$,
  $\Gamma_R$ and $\Gamma_E$}

We will also obtain uniform asymptotics near the curves $\Gamma_P$,
$\Gamma_Q$,  $\Gamma_R$ and
$\Gamma_E$, as well as in neighborhoods of the branch points $z_k$.
We start with the polynomials $Q_n$ and the remainder function $E_n$.
\begin{theorem}
\label{asym-Q}
Uniformly for $z$ in compact subsets of the region
\[ \left(\Gamma_Q \setminus \{z_1,z_2, z_3, z_4\} \right) \cup D_P \cup D_R
    \cup D_{\infty,U} \cup D_{\infty, L}, \]
see Figure \ref{fig:trajecta},
we have
\begin{eqnarray} \nonumber
Q_n(z) & = & z^{3n} e^{-n(g_P(z) + g_R(z))}
    \left[
    \frac{e^{-2n \varphi_P(z)}}
{\sqrt{3\psi_P^4(z)+1}} \left(1 +
    \O\left(\frac{1}{n}\right)\right) \right. \\
& &  \label{asympQ2} \quad
    \left.
    -\frac{1}{\sqrt{3\psi_Q^4(z)+1}} \left(1 + \mathcal{ O} \left(\frac{1}{n}\right)\right)
    +\frac{e^{-2n \varphi_R(z)} }
{\sqrt{3\psi_R^4(z)+1}}\left(1 +
    \O \left(\frac{1}{n}\right)\right)
    \right].
\end{eqnarray}

Uniformly for $z$ in compact subsets of $D_{\infty}$
(which includes the arcs $\Gamma_E$ but not the branch points $z_1,z_2,z_3,z_4$),
we have
\begin{eqnarray} \nonumber
E_n(z) & = & -z^{3n} e^{-n(g_P(z) + g_R(z))}
    \left[
    \frac{e^{-2n \varphi_P(z)}}
{\sqrt{3\psi_P^4(z)+1}} \left(1 +
    \O\left(\frac{1}{n}\right)\right) \right. \\
& &  \label{asympE2} \quad
    \left.
    +\ \frac{1}{\sqrt{3\psi_Q^4(z)+1}} \left(1 + \mathcal{ O} \left(\frac{1}{n}\right)\right)
    +\frac{e^{-2n \varphi_R(z)} }
{\sqrt{3\psi_R^4(z)+1}}\left(1 +
    \O \left(\frac{1}{n}\right)\right)
    \right].
\end{eqnarray}
\end{theorem}
Note that the asymptotic formula (\ref{asympQ2}) (resp. (\ref{asympE2}))
holds in particular on
$\Gamma_Q$ (resp. $\Gamma_E$), away from the branch points $z_1$, $z_2$, $z_3$ and $z_4$,
that is, on the curve where the zeros of $Q_n$ (resp. $E_n$) accumulate.

It may be checked that the three terms in (\ref{asympQ2})
are analytic continuations of the different asymptotic formulas we have
in (\ref{asympQ}). For $z \in D_{\infty}$, we have $g_Q(z) = -g_P(z) - g_R(z) + 3 \log z$
by Corollary \ref{corexplicit},
and then (\ref{asympQ2}) reduces to the first formula in (\ref{asympQ})
since $\Re \varphi_P(z) > 0$ and $\Re \varphi_R(z) > 0$ in $D_{\infty,U}$ and $D_{\infty,L}$.
For $z \in D_P$, we have $g_Q(z) = -g_P(z)-g_R(z) + 3 \log z - 2 \varphi_P(z)$
and we obtain the second formula in (\ref{asympQ}). Finally, for $z \in D_R$,
we have $g_Q(z) = - g_P(z) - g_R(z) + 3 \log z - 2 \varphi_R(z)$, and we
get the third formula in (\ref{asympQ}).

On $\Gamma_Q$ two of the terms in (\ref{asympQ2}) have comparable absolute values.
This causes the zeros of $Q_n$ to be close to $\Gamma_Q$. At the
points $\pm i y^*$,
all three terms have comparable absolute values. Similarly, on $\Gamma_E$ two of the terms
in (\ref{asympE2}) have comparable absolute values, which
causes the zeros of $E_n$ in $D_\infty$ to be close to $\Gamma_E$.

For the polynomials $P_n$ and $R_n$ we have an asymptotic formula on $\Gamma_P$
and $\Gamma_R$, respectively, away from the branch points, which now consists
of two terms.

\begin{theorem}
\label{asym-PR} Uniformly for $z$ in compact subsets of the region
\[ \left(\Gamma_P \setminus \{z_1, z_2 \} \right) \cup
    D_{\infty,P} \cup D_P^* \]
(this is the shaded region in Figure \ref{fig:Re}),
we have
\begin{equation} \label{asympP1}
(-2)^{n+1} P_n(z)= e^{n g_P(z)}
    \left[ \frac{2}{\sqrt{3\psi_{P}^4(z)+1}} \left(1 + \O\left(\frac{1}{n}\right)\right)
        \pm \frac{2e^{2n \varphi_P(z)}}{\sqrt{3\psi_{Q}^4(z)+1}} \left(1 + \O\left(\frac{1}{n}\right)\right) \right]
\end{equation}
where the $+$sign holds in $D_P \cup D_R$ and the $-$sign holds in $D_{\infty}$.

Uniformly for $z$ in compact subsets of the domain
\[ \left(\Gamma_R \setminus \{z_3, z_4 \} \right) \cup \{ z\in \mathbb C \mid \Re \varphi_R(z) < 0 \} \]
we have
\begin{equation} \label{asympR1}
(-2)^{n+1} R_n(z)= e^{n g_R(z)}
    \left[ \frac{2}{\sqrt{3\psi_{R}^4(z)+1}} \left(1 + \O\left(\frac{1}{n}\right)\right)
        \pm \frac{2e^{2n \varphi_R(z)}}{\sqrt{3\psi_{Q}^4(z)+1}} \left(1 + \O\left(\frac{1}{n}\right)\right) \right]
\end{equation}
where the $+$sign holds in $D_R \cup D_P$ and the $-$sign holds in $D_{\infty}$.
\end{theorem}

For $z$ away from $\Gamma_P$ so that $\Re \varphi_P(z) < 0$, the asymptotic formula
(\ref{asympP1}) reduces to (\ref{asympP}). On $\Gamma_P$ we have $\Re \varphi_P(z) =0$,
and then the two terms in (\ref{asympP1}) are of comparable magnitude. For $z \in \Gamma_P$,
(\ref{asympP1}) can be re-written as
\[ (-2)^{n+1} P_n(z) = \frac{2e^{n g_{P-}(z)}}{\sqrt{3\psi_{P-}^4(z)+1}} \left(1 + \O\left(\frac{1}{n}\right)\right)
    + \frac{2e^{n g_{P+}(z)}}{\sqrt{3\psi_{P+}^4(z)+1}} \left(1 + \O\left(\frac{1}{n}\right)\right). \]
Similar remarks hold for (\ref{asympR1}).

\subsection{Asymptotics near the branch points}
Near the branch points, the asymptotic formulas involve the Airy
function $\Ai$, which is the solution of the differential equation
$$ y''(z)=z y(z),$$
satisfying
$$\Ai(z) = \frac{1}{2\sqrt{\pi}} z^{-1/4} e^{-\frac{2}{3}z^{3/2}}
    \left(1 + \O\left(\frac{1}{z^{3/2}}\right)\right)
$$
as $z \to \infty$ with $| \arg z | < \pi$, see, e.g., \cite{OLV} for more
details. We only deal with the asymptotic behavior of $P_n$, $Q_n$ and
$E_n$ near
$z_1$. Similar results can be given for the behavior near the other branch
points.
Here we will take the branch of the function $\varphi_P(z)$ which is $0$
at $z=z_1$. So it behaves like
\[ \varphi_P(z) = c(z-z_1)^{3/2} + \O\left((z-z_1)^{5/2}\right) \]
as $z \to z_1$, with $c \neq 0$.
Then we define the function
\begin{equation} \label{deffz}
    f_1(z) = \left[\frac{3}{2} \varphi_P(z) \right]^{2/3}
\end{equation}
which is analytic for $z$ in a neighborhood of $z_1$.
We take that $2/3$rd power so that $f_1(z)$ is real and negative for $z \in \Gamma_P$.
Then
\begin{equation} \label{def-f}
    f_1(z) = c_1 (z-z_1) + \O\left((z-z_1)^2\right)
\end{equation}
as $z \to z_1$ with some constant $c_1 \neq 0$. Explicit calculations show that
\begin{equation} \label{c_1}
    c_1 = f_1'(z_1) = 2^{1/3} 3^{5/12} e^{-\frac{7}{36} \pi i}.
\end{equation}
\begin{theorem}
\label{asym-branch}
There is a $\delta > 0$, such that we have, uniformly for $|z-z_1| < \delta$,
\begin{multline} \label{asympP3}
(-2)^{n+1} P_n(z)  =  \sqrt{\pi} e^{(n+1) (g_P(z) + \varphi_P(z))} \\
 \left[ n^{1/6} h_1(z) \Ai\left( (n+1)^{2/3} f_1(z)\right) \left(1 + \O\left(\frac{1}{n}\right)\right)
  \right. \\   \left.
    + n^{-1/6} h_2(z) \Ai'\left((n+1)^{2/3} f_1(z)\right) \left(1+\O\left(\frac{1}{n} \right)\right)
\right],
\end{multline}
\begin{multline} \label{asympQ3}
(-2)^{n+1} Q_n(z) =  \sqrt{\pi} e^{3z} e^{(n+1) (g_P(z) + \varphi_P(z) -3z)} \\
\left[ n^{1/6} h_1(z) e^{-2 \pi i/3} \Ai\left(e^{-2\pi i/3} (n+1)^{2/3} f_1(z)\right)
    \left(1 + \O\left(\frac{1}{n}\right)\right)
    \right. \\ \left.
    + n^{-1/6} h_2(z) e^{2\pi i/3} \Ai'\left(e^{-2\pi i/3} (n+1)^{2/3} f_1(z)\right) \left(1+\O\left(\frac{1}{n} \right)\right)
\right],
\end{multline}
and
\begin{multline} \label{asympE3}
(-2)^{n+1} E_n(z) =  -\sqrt{\pi} e^{3z} e^{(n+1) (g_P(z) + \varphi_P(z) - 3z)} \\
\left[ n^{1/6} h_1(z) e^{2 \pi i/3} \Ai\left(e^{2\pi i/3} (n+1)^{2/3} f_1(z)\right)
    \left(1 + \O\left(\frac{1}{n}\right)\right)
    \right. \\ \left.
    + n^{-1/6} h_2(z) e^{-2\pi i/3} \Ai'\left(e^{2\pi i/3} (n+1)^{2/3} f_1(z)\right)
    \left(1+\O\left(\frac{1}{n} \right)\right)
\right],
\end{multline}
where $h_1$ and $h_2$ are two analytic functions without zeros
in $|z-z_1| < \delta$, which have explicit expressions
\begin{equation} \label{eq:defh1}
    h_1(z) = \left(N_{21}(z) + i z^{-1} e^{-3z} N_{22}(z) \right) f_1(z)^{1/4},
\end{equation}
with the branch of the fourth root in $f_1(z)^{1/4}$ taken with a cut along $\Gamma_P$, and
\begin{equation} \label{eq:defh2}
    h_2(z) = \left(-N_{21}(z) + iz^{-1} e^{-3z} N_{22}(z) \right) f_1(z)^{-1/4}.
\end{equation}
Here
\begin{equation} \label{eq:defN21N22}
N_{21}(z) = \frac{2 e^{-g_P(z)}}{\sqrt{3 \psi_P^4(z) + 1}}, \qquad
N_{22}(z) = - \frac{e^{g_P(z) + g_R(z)}}{z^2 \sqrt{3 \psi_Q^4(z) + 1}}.
\end{equation}
\end{theorem}
We use $N_{21}$ and $N_{22}$ to denote these functions, since in what follows
they will appear as entries in a matrix $N$.
Note that the function $g_P + \varphi_P$ is analytic near $z=z_1$.

From the asymptotics near the branch points, one can deduce the
behavior of the extreme zeros of $P_n$, $Q_n$, $R_n$ and $E_n$ near the
branch points. We only
state the result for the zeros of $P_n$, $Q_n$ and $E_n$ near $z_1$.
Recall that the Airy function $\Ai$ has only negative real zeros,
which we denote by $0 > -\iota_1 > -\iota_2 > \cdots > -\iota_{\nu} > \cdots$.
\begin{corollary} \label{asym-zero}
Let $z^P_{\nu,n}$, $\nu=1,\ldots,n$, be the zeros of $P_n$, ordered
by increasing distance to $z_1$. Then for every $\nu \in \mathbb N$,
we have
\begin{eqnarray} \nonumber
    z^P_{\nu,n} & = & z_1 -  \frac{\iota_\nu}{f_1'(z_1)} n^{-2/3}  +
\O\left(\frac{1}{n}\right),
\\
    & = & \label{zerosPn}
     z_1  +  2^{-1/3} 3^{-5/12} e^{-\frac{29}{36}\pi i}
    \iota_{\nu} n^{-2/3} + \O\left(\frac{1}{n}\right),
\end{eqnarray}
as $n \to \infty$.
Let $z^Q_{\nu,n}$, $\nu =1,\ldots, n$, be the zeros of $Q_n$, ordered
by increasing distance to $z_1$. Then for every $\nu \in \mathbb N$,
\begin{eqnarray} \nonumber
    z^Q_{\nu,n} & = & z_1 - e^{2\pi i/3} \frac{\iota_\nu}{f'_1(z_1)}
    n^{-2/3} + \O\left(\frac{1}{n}\right),  \\
    & = & \label{zerosQn}
    z_1 + 2^{-1/3} 3^{-5/12} e^{-\frac{5}{36} \pi i}
    \iota_{\nu} n^{-2/3} + \O\left(\frac{1}{n}\right),
\end{eqnarray}
as $n \to \infty$.
Let $z^E_{\nu,n}$, $\nu =1,\ldots, n$, be the zeros of $E_n$, ordered
by increasing distance to $z_1$. Then for every $\nu \in \mathbb N$,
\begin{eqnarray} \nonumber
    z^E_{\nu,n} & = & z_1 - e^{-2\pi i/3} \frac{\iota_\nu}{f'_1(z_1)}
    n^{-2/3} + \O\left(\frac{1}{n}\right),  \\
    & = & \label{zerosEn}
    z_1 + 2^{-1/3} 3^{-5/12} e^{\frac{19}{36} \pi i}
    \iota_{\nu} n^{-2/3} + \O\left(\frac{1}{n}\right),
\end{eqnarray}
as $n \to \infty$.
\end{corollary}

\subsection{Convergence of the Hermite--Pad\'e approximants}
The polynomials $P_{n}$, $Q_{n}$, $R_{n}$
satisfy
\[
P_n(z) + Q_n(z)e^{3nz}  + R_n(z) e^{6nz}=\O\left(z^{3n+2}\right),
\]
as $z\to 0$. Replacing the $\O$ term by $0$, one might expect that
the quadratic equation
\begin{equation} \label{quad}
    P_{n}(z)+Q_{n}(z)X_n(z)+R_{n}(z)X_n^2(z)=0
\end{equation}
has a solution
\begin{equation} \label{formulaXn}
    X_{n}(z)=\frac{-Q_n(z)\pm\sqrt{Q_n^2(z)-4P_n(z)R_n(z)}}{2R_n(z)}
\end{equation}
which is close to the exponential $e^{3nz}$.
This is indeed true for $z$ in the domains $D_P$ and $D_R$.
\begin{theorem}
\label{algebra}
For $z\in D_P\cup D_R$, there is a solution to {\rm(\ref{quad})} such that
\begin{equation} \label{conv-Xn}
    X_{n}(z)= e^{3nz} \left(1+ \O\left(\frac{1}{n}\right) \right).
\end{equation}
The $\O$-term is uniform on compact subsets of $D_P\cup D_R$.
\end{theorem}

\begin{remark} The domain of convergence does {\it not} include the
imaginary axis. On the imaginary axis one has roots of $Q_n^2-4P_nR_n$
which are branch points of $X_n$, accumulating along the interval
$[-i y^*,i y^*]$. Hence, convergence does not take place in a
neighborhood of 0 even though the approximants come from a
Hermite--Pad\'e expansion around 0.
\end{remark}
\begin{remark}
It turns out that $Q_n^2$ dominates $P_nR_n$ in the region
$D_P\cup D_R$. If we agree to take the square root in (\ref{formulaXn})
so that
\[ \sqrt{Q_n^2(z)-4P_n(z)R_n(z)} \simeq Q_n(z), \]
then in (\ref{formulaXn}) we have to take the $+$ sign in
$D_P$ and the $-$ sign in $D_R$. So
\begin{eqnarray*}
    X_n(z) &= & \frac{-Q_n(z)+\sqrt{Q_n^2(z)-4P_n(z)R_n(z)}}{2R_n(z)} \\
    & = & \frac{-2P_n(z)}{Q_n(z)+\sqrt{Q_n^2(z)-4P_n(z)R_n(z)}}
    \simeq -\frac{P_n(z)}{Q_n(z)} \qquad \mbox{for } z \in D_P,
\end{eqnarray*}
and
\begin{eqnarray*}
    X_n(z) &= & \frac{-Q_n(z)-\sqrt{Q_n^2(z)-4P_n(z)R_n(z)}}{2R_n(z)}\simeq
    -\frac{Q_n(z)}{R_n(z)} \qquad\mbox{for } z \in D_R.
\end{eqnarray*}
\end{remark}
\begin{remark}
It is possible to show that in $D_{\infty}$, there is no relative
convergence to $e^{3nz}$ of any of the two roots $X_n(z)$. Absolute
convergence holds in $D_P$ and in a subdomain of
$D_R$ bounded by the curve $\Re(3z+2\varphi_R(z))=0$.
\end{remark}

\section{Geometry of the problem}

In this section we will prove Propositions  \ref{prop:GammaPR} and \ref{prop:GammaP*R*}
together.

\subsection{Proof of Propositions \ref{prop:GammaPR} and \ref{prop:GammaP*R*}}

We recall that $\psi_P$, $\psi_Q$, and $\psi_R$ will be the three inverse mappings
of
\[ z = z(w) = \frac{w^2-\frac{1}{3}}{w(w^2-1)}. \]
It is easy to see that $z(w)$ is real if $w$ is real, and that $z(w)$ is purely imaginary
if $w \in i \mathbb R$. In fact, for $z \in i\mathbb R$, there is one solution
of the cubic equation (\ref{eq:cubic}) on the imaginary axis, one in the left half-plane,
and one in the right half-plane. We denote these by $\psi_Q(z)$, $\psi_P(z)$, and $\psi_R(z)$,
respectively.
Then part (b) of Proposition \ref{prop:GammaPR} holds. We also note that $\psi_Q(0) = \infty$.

These functions have analytic continuations, first to a vertical strip
$-\varepsilon < \Re z < \varepsilon$. around the imaginary axis.
Then $\psi_Q$ has a pole at $0$ and from (\ref{eq:zw}) it follows that
\begin{equation} \label{eq:poleQ}
    \psi_Q(z) = \frac{1}{z} + \O(z) \qquad \mbox{as } z \to 0.
\end{equation}
Near infinity, we have from analyzing (\ref{eq:zw})
\begin{eqnarray}
   \psi_P(z) & = & -1 + \frac{1}{3z} + \O\left(\frac1{z^2}\right), \label{eq:Pinf} \\
   \psi_Q(z) & = &  \frac{1}{3z} + \O\left(\frac1{z^2}\right), \label{eq:Qinf} \\
   \psi_R(z) & = & 1 + \frac{1}{3z} + \O\left(\frac1{z^2}\right),  \label{eq:Rinf}
\end{eqnarray}
as $z \to \infty$.

For $0 < \Re z < \varepsilon$, we then have that $\Re \psi_Q(z) > 0$, $\Re \psi_R(z) > 0$,
while $\Re \psi_P(z) < 0$.
Since for $z \notin i\mathbb R$, there can be no solution of (\ref{eq:cubic})
on the imaginary axis, it then follows by continuity that for every $z$ in the right half-plane,
there are two solutions of the cubic equation (\ref{eq:cubic})
with positive real part, and one solution with negative real part.
For $z$ with $\Re z > 0$, we then define $\psi_P(z)$ as the unique solution
of (\ref{eq:cubic}) with negative real part. This defines $\psi_P$
as an analytic function for $\Re z > - \varepsilon$.
Similarly, for $\Re z < 0$, we define $\psi_R(z)$ as the unique solution
of (\ref{eq:cubic}) with positive real part, and we have $\psi_R$ as an
analytic function for $\Re z < \varepsilon$.

For the moment, we let $\Gamma_P$ be an arbitrary simple curve that connects the
branch points $z_1$ and $z_2$, and lies entirely in the left half-plane.
We let $\Gamma_R$ be the mirror image of $\Gamma_P$ with respect to the imaginary axis
and we assume that these cuts determine the sheet structure of
the Riemann surface $\mathcal R$.
Then we can extend $\psi_P$, $\psi_Q$, and $\psi_R$ to analytic functions
on $\C \setminus \Gamma_P$, $\C \setminus (\Gamma_P \cup \Gamma_R)$,
and $\C \setminus \Gamma_R$, respectively. So we have part (a) of Proposition
\ref{prop:GammaPR}.

Next we want to study curves so that
$\frac{3}{2\pi i} \int_{z_1}^z (\psi_Q-\psi_P)(s) ds$ is real. If $z= z(t)$
is arclength parametrization of such a curve, then
\[ \Re \left[ \int_{z_1}^{z(t)} (\psi_Q-\psi_P)(s) ds \right] = 0. \]
which upon differentiating leads to $\Re \left[ z'(t) (\psi_Q-\psi_P)(z(t)) \right]  = 0$.
Since $z'(t) \neq 0$, and $(\psi_Q-\psi_P)(z) \neq 0$ (except at the branch
points $z_1$ and $z_2$), we find that
$- z'(t)^2 (\psi_Q-\psi_P)(z(t))^2 > 0$. Thus the curve is what
is known in geometric function theory as a trajectory of the quadratic differential
\begin{equation} \label{eq:quaddiff}
    - (\psi_Q-\psi_P)^2(z) dz^2,
\end{equation}
see \cite{Pom,Strebel}.
Now, $(\psi_Q-\psi_P)^2(z)$ is well-defined in the left half-plane, irrespective of
the exact choice for $\Gamma_P$. It is analytic with simple zeros at the branch
points $z_1$ and $z_2$, and a double pole at $0$.

Trajectories of the quadratic differential (\ref{eq:quaddiff}) which start from or end at
$z_1$ or $z_2$ are called critical trajectories. From the local structure of
trajectories of quadratic differentials, it is known that three trajectories
emanate from any simple zero. So three
critical trajectories emanate from $z_1$, and we would like to label
them by $\Gamma_P$, $\Gamma_P^*$, and $\Gamma_{E,1}$. The labeling will
depend on the global behavior of these trajectories. We consider
the global behavior  in the second quadrant, which we call $G$.

It is known that any trajectory in $G$ must begin and end either at $z_1$,
or at infinity, or on the boundary of $G$. It is also known that there can
be no closed Jordan curve consisting
of trajectories, since there are no poles of the quadratic differential in $G$.
For these properties of
trajectories of quadratic differentials, see \cite[Chapter 8]{Pom} and also \cite{BDMMZ}.
So the three critical trajectories that start at $z_1$, cannot end at $z_1$, and they
should end in $G$ in one of the following ways
\begin{enumerate}
\item[(1)] on the negative real axis,
\item[(2)] on the positive imaginary axis, or
\item[(3)] at infinity.
\end{enumerate}

If a critical trajectory meets the real axis, then by symmetry its mirror image
with respect to the real axis will be its continuation to $z_2$. If there were two critical trajectories
that come to the negative real axis, then they both would continue to $z_2$. Then there
would be a closed Jordan curve consisting of trajectories, which would not enclose a pole.
This is impossible. Hence there can be at most one critical trajectory that meets the negative real axis.

On the critical trajectories we have that $\varphi_P(z)$, see (\ref{eq:varphiP}), is purely imaginary.
Its derivative is $\varphi_P'(z) = \frac{3}{2\pi i} (\psi_Q-\psi_P)(z)$, and
\begin{equation} \label{eq:RevarphiP}
    \Re \varphi_P'(z) = \frac{3}{2\pi} \Im (\psi_Q-\psi_P)(z).
\end{equation}
We claim that $\Im \psi_Q(z) = \Im \psi_P(z)$ cannot happen for $z$ on the positive
imaginary axis.
To prove this claim, let us assume that for $z=iy$ with $y>0$, we have $\Im \psi_Q(z) = \Im \psi_P(z)$.
Note that $\psi_Q$ maps the positive imaginary axis to the negative imaginary axis.
So we have $\psi_Q(iy) = - iv$ with $v > 0$ and
\begin{equation} \label{eq:yv}
    y = \frac{v^2+\frac{1}{3}}{v(v^2+1)}.
\end{equation}
Recall that $\psi_P(z)$ is always in the left half-plane.
So we have $\psi_P(iy) = -u -iv$ with $u > 0$, and by symmetry, we
then have $\psi_R(iy) = u-iv$.
Since $\psi_P$, $\psi_Q$, and $\psi_R$ are the three solutions of (\ref{eq:cubic}),
we have
\[ \psi_P(z) + \psi_Q(z) + \psi_R(z) = \frac{1}{z}, \qquad z \in \C. \]
Then for $z = iy$, this leads to $-3iv = \frac{1}{iy}$, so that $y = \frac{1}{3v}$.
However, this is incompatible with (\ref{eq:yv}) if $v \neq 0$. This contradiction
proves the claim that $\Im \psi_Q(z) = \Im \psi_P(z)$ does not happen for $z$
on the positive imaginary axis.
As a result we see from (\ref{eq:RevarphiP}) that the real part of $\varphi_P(z)$ is strictly
monotone as $z$ varies over the positive imaginary axis. Thus there can be at most one  value
$y^* > 0$ so that $\Re \varphi_P(iy^*) = 0$. Since $\Re \varphi_P = 0$ on the critical trajectories,
we see that there can be at most one critical trajectory that meets the positive imaginary axis.

Now consider trajectories that end at infinity. At infinity, we have
$-(\psi_Q-\psi_P)^2(\infty) = -1$. This means that the quadratic differential has a pole
of order $4$ at infinity, see \cite{Strebel}, and all trajectories that extend to infinity
arrive there with a vertical tangent. Assume that two critical trajectories extend to infinity
in $G$. Then these trajectories are the boundary of  a region $G'$ in $G$. Any trajectory
in $G'$ begins and ends at infinity with a vertical tangent. However,
according to \cite[Theorem 7.4]{Strebel} a pole of order $4$ has a neighborhood,
so that any closed trajectory lying entirely in that neighborhood begins and
ends at the pole, but from opposite directions. This is a contradiction, since there
will be trajectories in $G'$ that are arbitrarily close to infinity. This contradiction
shows that there can be at most one critical trajectory that extends to infinity.

So now we proved that each of the possible ways (1), (2), and (3) that a critical
trajectory can end in $G$, happens at most once. Since there are three critical trajectories,
it follows that each possibility happens exactly once. So there is a unique
critical trajectory that meets the negative real axis. This trajectory continues
to $z_2$, and we call it $\Gamma_P$. We then have that
part (c) of Proposition \ref{prop:GammaPR} is satisfied.

The critical trajectory that meets the positive imaginary axis is called $\Gamma_P^*$,
and the critical trajectory that extends to infinity is called $\Gamma_{E,1}$.
By symmetry, we also find the trajectory $\Gamma_{E,2}$ that emanates
from $z_2$ and extends to infinity.

By symmetry with respect to the imaginary axis, we also find the contours $\Gamma_R$,
$\Gamma_{E,3}$, $\Gamma_{E,4}$, and $\Gamma_R^*$ in the right half-plane, so that
for $z$ on these contours we have that $\frac{3}{2\pi i} \int_{z_3}^z (\psi_Q-\psi_R)(s) ds$
is real. This establishes part (d) of Proposition \ref{prop:GammaPR}
and this proposition is now proved completely.

Large parts of Proposition \ref{prop:GammaP*R*} have also been proved. All statements
regarding $\Gamma_P$, $\Gamma_R$, $\Gamma_{E,k}$, $k=1,2,3,4$ are proved.
We also know that $\Gamma_P^*$ and $\Gamma_R^*$ meet at the point $iy^*$ on
the positive imaginary axis. What remains is to investigate how $\Gamma_P^*$
and $\Gamma_R^*$ continue beyond $iy^*$.
The following lemma will be useful.

\begin{lemma} \label{lem:Revarphi}
The real part of $\varphi_P$ is strictly positive on $\Gamma_R$, $\Gamma_{E,3}$,
and $\Gamma_{E,4}$.
The real part of $\varphi_R$ is strictly positive on $\Gamma_P$, $\Gamma_{E,1}$,
and $\Gamma_{E,2}$.
\end{lemma}

We delay the proof of this lemma to Section \ref{subsec54},
since its proof is easier once we have the functions $g_P$ and $g_R$ and know
some of their properties. The reader may verify that the proofs that follow
in Section 4 do not depend on the curves $\Gamma_P^*$ and $\Gamma_R^*$.

Assuming the lemma we can complete the proof of Proposition \ref{prop:GammaP*R*}
as follows. Consider the quadratic differential $-(\psi_Q-\psi_P)(z)^2 dz^2$
in the domain $G_1$ bounded by $\Gamma_{E,4}$,  the curve $\Gamma_R$ in the upper
half-plane and the real interval $(-\infty, x^*]$ where $x^*$ is the point
where $\Gamma_R$ intersects the real axis. The critical trajectory $\Gamma_P^*$
has to end in $G_1$ in one of the following ways
\begin{enumerate}
\item[(1)] on the real axis,
\item[(2)] on $\Gamma_R$ or $\Gamma_{E,4}$,
\item[(3)] at infinity.
\end{enumerate}
It cannot end on $\Gamma_R$ or $\Gamma_{E,4}$ because of Lemma \ref{lem:Revarphi},
and it cannot extend to infinity for the same reason that it could not
extend to infinity in the left half-plane. So it has to end on the real axis,
but it cannot end on the negative real axis, again for the same reason as before.
So it has to end on the positive real axis, somewhere in the open interval $(0,x^*)$.
Then by symmetry $\Gamma_P^*$ continues in the lower half-plane to the branch point
$z_2$ and $\Gamma_P^*$ does not intersect with $\Gamma_R$.
This proves part (c) of Proposition \ref{prop:GammaP*R*}.

By symmetry, also part (d) holds, and this completes the proof of
Proposition \ref{prop:GammaP*R*}.

\section{Measures and functions associated with the Riemann surface}

In this section we study the measures $\mu_P$, $\mu_Q$, $\mu_R$, and $\mu_E$,
and the functions $\varphi_P$, $\varphi_R$, $g_P$, $g_Q$, and $g_R$.
These measures and functions that are associated with the Riemann surface
satisfy many relations that will be used in the transformations
of the Riemann--Hilbert problem that follow in later sections.
We also prove Theorem \ref{mes-proba} and Lemmas \ref{lemma:phiPR}
and \ref{lem:Revarphi} in this section.

\subsection{Properties of the measures $\mu_P$, $\mu_R$, and $\mu_E$}

We start with a lemma.
\begin{lemma}  \label{lem:int1}
We have
\begin{equation} \label{eq:intGammaP}
    \frac{3}{2\pi i} \int_{\Gamma_P} (\psi_Q - \psi_P)_+(s) \, ds = 1,
\end{equation}
and
\begin{equation} \label{eq:intGammaR}
   \frac{3}{2\pi i} \int_{\Gamma_R} (\psi_Q - \psi_R)_+(s) \, ds = 1.
\end{equation}
\end{lemma}
\begin{proof}
Let $\gamma$ be a closed contour on the sheet $\mathcal{R}_P$ going around $\Gamma_P$ once
in the positive direction. Then the residue theorem for the exterior of $\gamma$ gives
\[   \frac{3}{2\pi i} \int_{\gamma} \psi_P(s)\, ds = 1, \]
because $\psi_P$ is analytic outside $\gamma$ and we have (\ref{eq:Pinf}).
If we shrink $\gamma$ to $\Gamma_P$, then the integral becomes
\[  \frac{3}{2\pi i} \int_{\Gamma_P} ( (\psi_P)_-(s) - (\psi_P)_+(s) )\, ds = 1. \]
Taking into account that $(\psi_P)_- = (\psi_Q)_+$, we obtain (\ref{eq:intGammaP}).
The reasoning is similar for the proof of (\ref{eq:intGammaR}), where we use a closed
contour going around $\Gamma_R$ on $\mathcal R_R$ and the behavior (\ref{eq:Rinf})
of $\psi_R$ at infinity.
\end{proof}

Now we can prove the parts of Theorem \ref{mes-proba} dealing
with $\mu_P$, $\mu_R$ and $\mu_E$. We give the proof of the part
dealing with $\mu_Q$ only after we have a better understanding of
the curve $\Gamma_P^*$ and $\Gamma_R^*$, see Proposition \ref{prop:muQ}.

\begin{proposition} \label{prop:muPmuR}
$\mu_P$ is a probability measure on $\Gamma_P$, $\mu_R$ is
a probability measure on $\Gamma_R$ and $\mu_E$ is a positive measure
on $\Gamma_E$.
\end{proposition}
\begin{proof}
The curve $\Gamma_P$ is such that for $z \in \Gamma_P$, the integral
$\frac{3}{2\pi i} \int_{z_1}^z (\psi_Q-\psi_P)_+(s) \,ds$
is real, see Proposition \ref{prop:GammaPR}. For $z=z_1$, it has the value $0$,
and for $z = z_2$ it has the value $1$ by (\ref{eq:intGammaP}).
Let $z = z(t)$, $t \in [0,T]$, be arclength parametrization of
$\Gamma_P$. The derivative of
\begin{equation} \label{eq:parametr}
    t \mapsto \frac{3}{2\pi i} \int_{z_1}^{z(t)} (\psi_Q-\psi_P)_+(s) \,ds
\end{equation}
is equal to $\frac{3}{2\pi i} \left(\psi_Q(z(t)) - \psi_P(z(t))\right) z'(t)$
and this is different from $0$ for $t \in (0,T)$.
Thus (\ref{eq:parametr}) is strictly increasing from $0$ for $t=0$
to $1$ for $t=T$.
This immediately implies that $\mu_P$ defined
by (\ref{eq:muP}) is a probability measure on $\Gamma_P$.
Similarly $\mu_R$ defined in (\ref{eq:muR}) is a probability measure on $\Gamma_R$.

For $\mu_E$ we observe that (\ref{eq:muE}) defines a real measure on $\Gamma_E$,
since, by Proposition \ref{prop:GammaPR}, $\frac{3}{2\pi i} \int_{z_1}^z
(\psi_Q - \psi_P)(s) ds$ is real for $z \in \Gamma_{E,1} \cup \Gamma_{E,2}$,
and $\frac{3}{2\pi i} \int_{z_3}^z (\psi_Q - \psi_R)(s) ds$ is
real for $z \in \Gamma_{E,3} \cup \Gamma_{E,4}$.
Using an argument based on arclength parametrization, similar to the
one above, we find that on each part $\Gamma_{E,j}$, $j=1,2,3,4$, the
measure $\mu_E$ is either positive or negative.
Since $\psi_Q(s)-\psi_P(s)=1 + \mathcal{O}(1/s^2)$ as $s \to \infty$ (see
(\ref{eq:Pinf})--(\ref{eq:Qinf})), we have for $z \in \Gamma_{E,1}$,
\[  \frac{3}{2\pi i} \int_{z_1}^{z} (\psi_Q-\psi_P)(s)\, ds = \frac{3}{2\pi i} (z-z_1)
   + \mathcal{O}(1) \qquad \mbox{ as } z \to \infty, \, z \in \Gamma_{E,1}. \]
Since $\Im(z-z_1) \to +\infty$ as $z \to \infty$ along $\Gamma_{E,1}$,
the integral is positive  as $z \to \infty$ along $\Gamma_{E,1}$. As the
measure $\mu_E$ is of constant sign on $\Gamma_{E,1}$, we may thus
deduce that it is positive everywhere on $\Gamma_{E,1}$. The reasoning
for $\Gamma_{E,2}$, $\Gamma_{E,3}$, and $\Gamma_{E,4}$ is similar.

This  completes the proof of Proposition \ref{prop:muPmuR}.
\end{proof}

\subsection{Properties of $g_P$ and $g_R$}

The functions $g_P$ and $g_R$ were defined in (\ref{eq:gP}) and (\ref{eq:gR}).
These are multi-valued functions, depending on the choice of the branch of
the logarithm $\log (z-s)$, which we assume depends on $s \in \Gamma$ in
a continuous way. Since $\mu_P$  and $\mu_R$ are probability measures,
the $g$-functions are defined modulo $2\pi i$.

\begin{lemma}  \label{lem:gder}
For the derivatives of the functions $g_P$ and $g_R$ we have
\begin{equation} \label{eq:gPder}
   g_P'(z)  =  3 \psi_P(z)+3,   \qquad z \in \mathbb C \setminus \Gamma_P,
\end{equation}
and
\begin{equation} \label{eq:gRder}
   g_R'(z)  =  3 \psi_R(z)-3,   \qquad z \in \mathbb C \setminus \Gamma_R.
\end{equation}
\end{lemma}
\begin{proof}
The derivative of $g_P$ is easily obtained as
\[  g_P'(z) = \frac{3}{2\pi i} \int_{\Gamma_P} \frac{1}{z-s} (\psi_Q-\psi_P)_+(s) \,ds. \]
If $\gamma$ is a closed contour going around $\Gamma_P$ on
$\mathcal{R}_P$ in the positive direction but with $z$ outside
$\gamma$, then  (since $(\psi_Q)_+=(\psi_P)_-$)
\[  g_P'(z) = \frac{3}{2\pi i} \oint_{\gamma} \frac{\psi_P(s)}{z-s} \, ds. \]
The integral over $\gamma$ can be calculated with the residue
theorem for the exterior of $\gamma$, for which there is a residue
at $z$ and at $\infty$ given by (\ref{eq:Pinf}). This proves
(\ref{eq:gPder}). The proof of (\ref{eq:gRder}) is similar.
\end{proof}

\begin{corollary}
For $z \in \mathbb C \setminus (\Gamma_P \cup \Gamma_R)$, we have
\begin{equation} \label{eq:gPRder1}
   2 g_P'(z)+g_R'(z) = \frac{3}{z} + 3 - 3[\psi_Q(z) - \psi_P(z)],
\end{equation}
and
\begin{equation} \label{eq:gPRder2}
   2 g_R'(z)+g_P'(z) = \frac{3}{z} - 3 - 3[\psi_Q(z) - \psi_R(z)].
\end{equation}
\end{corollary}
\begin{proof}
We recall that $\psi_P(z)$, $\psi_Q(z)$, and $\psi_R(z)$ are roots of the
cubic equation (\ref{eq:cubic}), so that
\begin{equation} \label{sum-PQR}
    \psi_P(z)+ \psi_Q(z)+\psi_R(z)= \frac{1}{z}.
\end{equation}
Hence by (\ref{eq:gPder}) and (\ref{eq:gRder}),
\begin{eqnarray*}
   2g_P'(z)+g_R'(z) & = & 6\psi_P(z) + 6 + 3 \psi_R(z) - 3 \\
                    & = & 3[\psi_P(z)+\psi_Q(z)+\psi_R(z)] -3[\psi_Q(z)-\psi_P(z)] + 3 \\
                    & = & \frac{3}{z} + 3 - 3[\psi_Q(z)-\psi_P(z)].
\end{eqnarray*}
A similar computation gives (\ref{eq:gPRder2}).
\end{proof}

A useful explicit expression of $g_P$ and $g_R$ in terms of the
mapping functions $\psi_P$ and $\psi_R$ is given in the next lemma.
\begin{lemma} \label{explicit}
The functions $g_P$ and $g_R$ have the following explicit
expressions in terms of the mapping functions $\psi_P$  and $\psi_R$,
\begin{equation} \label{eq:gPexplicit}
    g_P(z) = 3z(\psi_P(z)+1) - \log[\psi_P(z)(\psi_P^2(z)-1)] - 1 + \log \left(\frac{2}{3} \right),
         \mbox{ for } z \in \mathbb C \setminus \Gamma_P,
\end{equation}
and
\begin{equation} \label{eq:gRexplicit}
    g_R(z) = 3z(\psi_R(z)-1) - \log[\psi_R(z)(\psi^2_R(z)-1)] - 1 + \log \left(\frac{2}{3}\right),
        \mbox{ for } z \in \mathbb C \setminus \Gamma_R.
\end{equation}
\end{lemma}
\begin{proof}
We let $z \in \mathbb C \setminus \Gamma_P$ and put $w = \psi_P(z)$.
Taking a derivative of (\ref{eq:zw}), we find, since $w'(z) = \frac{1}{z'(w)}$,
\[ 3 = - \left(\frac{1}{w^2} + \frac{1}{(w-1)^2} + \frac{1}{(w+1)^2} \right) w'. \]
Thus
\begin{eqnarray*}
   3\psi_P(z) + 3 & = & - (w+1) \left(\frac{1}{w^2} + \frac{1}{(w-1)^2} + \frac{1}{(w+1)^2} \right) w' \\
     & = & - \left(\frac{1}{w} + \frac{1}{w^2} + \frac{1}{w-1} + \frac{2}{(w-1)^2} + \frac{1}{w+1} \right) w' \\
     & = & \frac{d}{dz} \left[ - \log w(w^2-1) + \frac{1}{w} + \frac{2}{w-1} \right].
     \end{eqnarray*}
By (\ref{eq:gPder}), we then see that
\[ g_P(z) = - \log w(w^2-1) + \frac{1}{w} + \frac{2}{w-1} + C \]
for some constant $C$. The constant can be determined from the behavior for $z \to \infty$,
since $g_P(z) = \log z + \O(1/z)$, and $w = \psi_P(z) = -1 + \frac{1}{3z} + \O(1/z^2)$.
The result is that $C = 2 + \log \frac{2}{3}$.
Using (\ref{eq:zw}) we then find (\ref{eq:gPexplicit}).

The proof of (\ref{eq:gRexplicit}) follows along similar lines.
\end{proof}

\subsection{Jump properties of $g_P$ and $g_R$}

We recall that the functions $\varphi_P$ and $\varphi_R$ were
introduced in  (\ref{eq:varphiP}) and (\ref{eq:varphiR}).
The next lemma connects these functions with $g_P$ and $g_R$.
It will be frequently used in what follows.
Throughout the rest of the paper we use $\ell$ to denote the constant
\begin{equation} \label{eq:ell}
    \ell = \log 2 - \pi i.
\end{equation}

\begin{lemma}  \label{lem:gPR}
For $z \in \mathbb C \setminus \left( \Gamma_P \cup \Gamma_R \right)$, we have
\begin{eqnarray}
    2g_P(z)+g_R(z) & = & 3 \log z + 3z - 2 \varphi_P(z) + \ell,  \label{eq:gPR1} \\[10pt]
    2g_R(z)+g_P(z) & = & 3 \log z - 3z - 2 \varphi_R(z) + \ell.  \label{eq:gPR2}
\end{eqnarray}
On the contours we have
\begin{eqnarray}
     g_{P+}(z)+g_{P-}(z)+g_R(z) & = & 3\log z + 3z +\ell,  \qquad z \in \Gamma_P, \label{eq:gPRG1} \\[10pt]
     g_P(z)+g_{R+}(z)+g_{R-}(z) & = & 3\log z - 3z +\ell,  \qquad z \in \Gamma_R, \label{eq:gPRG2}
\end{eqnarray}
and
\begin{eqnarray}
    g_{P+}(z) - g_{P-}(z) & = & -\varphi_{P+}(z) + \varphi_{P-}(z) = -2\varphi_{P+}(z) = 2 \varphi_{P-}(z),
        \qquad z \in \Gamma_P, \label{eq:gPRG3} \\[10pt]
    g_{R+}(z) - g_{R-}(z) & = & -\varphi_{R+}(z) + \varphi_{R-}(z) = -2\varphi_{R+}(z) = 2 \varphi_{R-}(z),
        \qquad z \in \Gamma_R \label{eq:gPRG4}.
\end{eqnarray}
\end{lemma}
\begin{proof}
Integrating (\ref{eq:gPRder1}) from $z_1$ to $z$ over some path in
$\mathbb C \setminus (\Gamma_P \cup \Gamma_R \cup \{ 0 \})$, we get
\[  2g_P(z) + g_R(z) - 2g_P(z_1) - g_R(z_1) =  3 \log z - 3 \log z_1 + 3(z-z_1) - 2\varphi_P(z), \]
so that (\ref{eq:gPR1}) holds with constant
\[ \ell = 2 g_P(z_1) + g_R(z_1) - 3\log z_1 - 3z_1. \]
Using the explicit expressions (\ref{eq:gPexplicit}) and (\ref{eq:gRexplicit}) for $g_P$
and $g_R$ we are able to show that $\ell$ is equal to (\ref{eq:ell}).
To obtain (\ref{eq:gPR2}) we integrate (\ref{eq:gPRder2}) from $z_3$ to $z$. It turns out
that the constant of integration is again given by (\ref{eq:ell}).

Next, we use Lemma \ref{lem:gder} to find
\[   g_{P+}(z) + g_{P-}(z) + g_R(z) = 3 \int_{z_1}^z (\psi_{P+}+\psi_{P-}+\psi_R)(s)\,ds + 3(z-z_1),
   \qquad z\in \Gamma_P. \]
On $\Gamma_P$ we have $\psi_{P-}(s) = \psi_{Q+}(s)$, so that
\[   g_{P+}(z) + g_{P-}(z) + g_R(z) = 3 \int_{z_1}^z (\psi_P+\psi_Q+\psi_R)_+(s)\,ds + 3z - 3z_1,
\qquad z\in \Gamma_P. \]
Since $\psi_P(s)+\psi_Q(s)+\psi_R(s)=\frac1s$, we obtain
\[   g_{P+}(z) + g_{P-}(z) + g_R(z) = 3 \log z - 3 \log z_1 + 3z - 3z_1 = 3 \log z - 3z + \ell,
\qquad z\in \Gamma_P. \]
This proves (\ref{eq:gPRG1}).  A similar computation leads to (\ref{eq:gPRG2}).

Finally, if we take (\ref{eq:gPR1}) on the $+$ and $-$sides of $\Gamma_P$ and subtract,
we get
\[ 2 g_{P+}(z) - 2g_{P-}(z) = - 2 \varphi_{P+}(z) + 2 \varphi_{P-}(z), \qquad z \in \Gamma_P. \]
This gives (\ref{eq:gPRG3}), since $\varphi_{P+}(z) = - \varphi_{P-}(z)$. Similarly
we find (\ref{eq:gPRG4}).
\end{proof}

Since $g_P(z)=\log z + \O(1/z)$ and $g_R(z)=\log z + \O(1/z)$
as $z \to \infty$, we get from (\ref{eq:gPR1})--(\ref{eq:gPR2}) that
\begin{equation} \label{eq:varphiPasymp}
   2 \varphi_P(z) = 3z + \ell + \O\left( \frac1z \right)
\qquad \mbox{as } z \to \infty,
\end{equation}
and
\begin{equation} \label{eq:varphiRasymp}
   2 \varphi_R(z) = -3z + \ell + \O\left( \frac1z \right)
 \qquad \mbox{as }z \to \infty.
\end{equation}
We also see from (\ref{eq:gPR1}) and (\ref{eq:gPR2}) that $\varphi_P$
and $\varphi_R$ are multivalued functions, which are defined modulo $\pi i$,
since $g_P$ and $g_R$ are defined modulo $2\pi i$.

\subsection{Proof of Lemma \ref{lem:Revarphi}} \label{subsec54}

Now we are in a position to prove Lemma \ref{lem:Revarphi}
which was needed in Section 3 to complete the discussion on the contours
$\Gamma_P^*$ and $\Gamma_R^*$.

\begin{proof}
Let $\gamma$ be the contour which consists of $\Gamma_R$,
$\Gamma_{E,3}$ and $\Gamma_{E,4}$.
We define a function $G_R$ on the right half-plane $\{ \Re z > 0 \}$ by
\[ G_R(z) = \begin{cases}
       g_R(z) & \textrm{to the left of } \gamma, \\
       g_R(z) + 2\varphi_R(z) & \textrm{to the right of } \gamma.
      \end{cases}  \]
This function is analytic across $\Gamma_R$ because of the jump
relation (\ref{eq:gPRG4}) of $g_R$. The function
$u(z) = \Re G_R(z)$ is therefore harmonic in $\{ \Re z > 0 \} \setminus
(\Gamma_{E,3} \cup \Gamma_{E,4})$. On $\Gamma_{E,3}$ and $\Gamma_{E,4}$
we have $\Re \varphi_R = 0$, and thus $u$ is continuous in the full
right half-plane. Furthermore $\Re \varphi_R < 0$ on the right of
$\gamma$. Since $g_R$ is analytic across $\Gamma_{E,3}$ and $\Gamma_{E,4}$
it then easily follows that $u$ is superharmonic in $\{ \Re z > 0 \}$.
Near infinity we have
\[   u(z) = \begin{cases}
       \Re \log z + \O(1/z) & \textrm{as $z \to \infty$ to the left of } \gamma, \\
       \Re (-3z + \log z + \ell) + \O(1/z) &
    \textrm{as $z\to\infty$ to the right of } \gamma.
          \end{cases} \]
By symmetry $\Re g_P(z) = \Re g_R(z)$ on the imaginary axis.
Thus $\Re g_P(z) - 3\Re z$ is harmonic in the right half plane,
$\Re g_P(z) - 3\Re z = u(z)$ on the imaginary axis, and
$u(z)-\Re g_P(z) + 3 \Re z$ is bounded from below as $z \to \infty$.
Then it follows from the minimum principle for superharmonic functions
that $u(z) - \Re g_P(z) + 3\Re z \geq 0$ in the right half-plane. In
particular this gives $\Re( g_R(z) - g_P(z) + 3 z) \geq 0$ to the
left of $\gamma$ and $\Re (g_R(z) - g_P(z) + 2\varphi_R(z)) \geq 0$
to the right of $\gamma$.
Taking the difference between (\ref{eq:gPR2}) and (\ref{eq:gPR1})
gives $g_R(z)-g_P(z)=-6z-2\varphi_R(z)+2\varphi_P(z)$.
Then we see that
$\Re(-3z -2\varphi_R(z)+2\varphi_P(z)) \geq 0$ to the left of $\gamma$
and $\Re(-3z+2\varphi_P(z)) \geq 0$ to the right of $\gamma$.
Letting now $z \to \gamma$ from either side, we get, since
$\Re \varphi_R(z) = 0$ on $\gamma$, that
$2 \Re \varphi_P(z) \geq 3 \Re z > 0$ on $\Gamma_R$,
$\Gamma_{E,3}$, and $\Gamma_{E,4}$.

The proof that $\Re \varphi_R > 0$ on $\Gamma_P$, $\Gamma_{E,1}$,
and $\Gamma_{E,2}$ is similar, and also follows because of symmetry.
\end{proof}

\subsection{Proof of Lemma \ref{lemma:phiPR}} \label{subsec55}

\begin{proof}
By Proposition \ref{prop:GammaPR} and the definition of $\varphi_P$,
we have that $\Re \varphi_P = 0$ on $\Gamma_P$, $\Gamma_P^*$, $\Gamma_{E,1}$,
and $\Gamma_{E,2}$.

We know that $\Re \varphi_P$ is a harmonic function in
$\mathbb C \setminus (\Gamma_P \cup \Gamma_R \cup \{0\})$.
Since $\psi_Q(z) \sim 1/z$ as $z \to 0$, it easily follows from
(\ref{eq:varphiP}) that $\Re \varphi_P(z) \to -\infty$ as $z \to 0$.
Then by the maximum principle for harmonic functions we get
that $\Re \varphi_P < 0$ on $D_P^*$.

As $z \to \infty$, we have (\ref{eq:varphiPasymp}).
On the unbounded curves $\Gamma_{E,1}$ and $\Gamma_{E,2}$
we have $\Re \varphi_P = 0$. From (\ref{eq:varphiPasymp})
it then follows that these curves are asymptotic to
the vertical line $\Re (3z + \ell) = 0$, which is
$\Re z = - \frac{1}{3} \log 2$.
As $z \to \infty$ in the unbounded domain $D_{\infty,P}$ we then
have $\limsup \Re (3z + \ell) \leq 0$, so that
$\limsup \Re \varphi_P(z) \leq 0$ by (\ref{eq:varphiPasymp}). Again it follows
by the maximum principle for harmonic functions that
$\Re \varphi_P < 0$ on $D_{\infty,P}$.

For the remaining domain $G = \mathbb C \setminus (\overline{D_P^* \cup
D_{\infty, P}} )$, we have that $\Re \varphi_P$ is harmonic
on $G \setminus \Gamma_R$, with $\liminf \varphi_P(z)\geq 0$ as $z \to \infty$ with
$z \in G$ by (\ref{eq:varphiPasymp}), and $\Re \varphi_P > 0$ on
$\Gamma_R$ by Lemma \ref{lem:Revarphi}.
Thus again by the maximum principle, $\Re \varphi_P > 0$ on $G$.
This completes the proof of part (a) of Lemma \ref{lemma:phiPR}.

The proof of part (b) is similar and also follows because of symmetry.

For part (c) we note that $\Re \varphi_P = \Re \varphi_R$ on
the imaginary axis because of symmetry. As $z \to \infty$, we have by
(\ref{eq:varphiPasymp}) and (\ref{eq:varphiRasymp}) that
$\liminf \Re (\varphi_R(z) - \varphi_P(z)) \geq 0$ as
$z \to\infty$ in the left half-plane. Since $\Re \varphi_R$
and $\Re \varphi_P$ are harmonic in $\{ \Re z < 0 \} \setminus \Gamma_P$
and $\Re \varphi_R > \Re \varphi_P$ on $\Gamma_P$, we find
that
$\Re \varphi_P < \Re \varphi_R$ in the left half-plane.
Similarly, we have $\Re \varphi_P > \Re \varphi_R$ in the right half-plane.

This completes the proof of Lemma \ref{lemma:phiPR}.
\end{proof}

\subsection{The measure $\mu_Q$} \label{subsec56}

Now that we have the full curves $\Gamma_P^*$ and $\Gamma_R^*$, we can
prove the part of Theorem \ref{mes-proba} that deals with $\mu_Q$. We start with
the analogue of Lemma \ref{lem:int1}.

\begin{lemma}  \label{lem:int2}
We have
\begin{equation} \label{eq:intGammaP*}
    \frac{3}{2\pi i} \int_{\Gamma_P^*} (\psi_Q - \psi_P)(s) \, ds = 2,
\end{equation}
and
\begin{equation} \label{eq:intGammaR*}
    \frac{3}{2\pi i} \int_{\Gamma_R^*} (\psi_Q - \psi_R)(s) \, ds = 2.
\end{equation}
\end{lemma}
\begin{proof}
If we deform the integral (\ref{eq:intGammaP*}) over $\Gamma_P^*$ to an integral
over $\Gamma_P$, then we pick up a residue contribution at $0$ due to the pole of
$\psi_Q$ at the origin. So by (\ref{eq:poleQ}) and (\ref{eq:intGammaP}) we have
\[ \frac{3}{2\pi i} \int_{\Gamma_P^*} (\psi_Q - \psi_P)(s) \, ds =
    3 - \frac{3}{2\pi i} \int_{\Gamma_P} (\psi_Q - \psi_P)_+(s) \, ds
    = 2.
\]
This gives (\ref{eq:intGammaP*}) and in a similar way we obtain (\ref{eq:intGammaR*}).
\end{proof}

Now we can prove that $\mu_Q$ is a probability measure.
\begin{proposition} \label{prop:muQ}
$\mu_Q$ is a probability measure on $\Gamma_Q$.
\end{proposition}
\begin{proof}
We prove first that (\ref{eq:muQ}) defines a positive measure.
On $\Gamma_P^*$ we have that
$\frac{3}{2\pi i} \int_{z_2}^z (\psi_Q-\psi_P)(s) ds$ takes on real values,
and has the value $0$ at $z=z_2$ and the value $2$ for $z=z_1$ by (\ref{eq:intGammaP*}).
Since its derivative does not vanish between $z_2$ and $z_1$, we find that
$\frac{3}{2\pi i} (\psi_Q-\psi_P)(s) ds$ is a positive measure on $\Gamma_P^*$
of total mass $2$. Its restriction
to $\Gamma_P^* \cap \{ \Re z < 0 \}$ is part of the measure $\mu_Q$, and so
this part is a positive measure. Similarly, $\frac{3}{2\pi i} (\psi_R - \psi_Q)(s) ds$
is a positive measure on $\Gamma_R^*$ (due to (\ref{eq:intGammaR*})),
and so also on $\Gamma_R^* \cap \{ \Re z > 0\}$.
On the imaginary axis, we have, because of symmetry, that $\psi_P(s) = - \overline{\psi_R(s)}$,
which means that $\psi_R(s) - \psi_P(s)$ is real for $s \in i \mathbb R$. Then
$\frac{3}{2\pi i} (\psi_R - \psi_P)(s) ds$ is a real measure on the imaginary axis.
As $s \to \infty$, we have $\psi_R(s) - \psi_P(s) = 2 + \O(\frac{1}{s})$.
Since $\psi_R - \psi_P$ never vanishes in ${\overline\C}$, this
implies that
$\frac{3}{2\pi i} (\psi_R-\psi_P)(s) ds$ is a positive measure on the imaginary
axis, and so in particular on the interval $[-iy^*, iy^*]$. Hence all three parts
of $\mu_Q$ in (\ref{eq:muQ}) are positive, so that $\mu_Q$ is a positive measure.

Now we want to calculate the total mass of $\mu_Q$. If we separate in
$\int d\mu_Q$ the contributions from $\psi_P$, $\psi_Q$, and $\psi_R$,
we get
\[ \int d\mu_Q =
    \frac{3}{2\pi i} \int_{z_1}^{z_2} \psi_P(s) ds
    + \frac{3}{2\pi i} \int_{z_2}^{z_3} \psi_Q(s) ds
    + \frac{3}{2\pi i} \int_{z_3}^{z_4} \psi_R(s) ds
    + \frac{3}{2\pi i} \int_{z_4}^{z_1} \psi_Q(s) ds.
\]
Here $\int_{z_j}^{z_k}$ stands for integration along the path from $z_j$ to $z_k$
within $\Gamma_Q$. We deform the contour from $z_1$ to $z_2$ to
the contour $\Gamma_P$. On $\Gamma_P$ we have $(\psi_P)_+ = (\psi_Q)_-$.
Thus
\[ \frac{3}{2\pi i} \int_{z_1}^{z_2} \psi_P(s) ds =
   \frac{3}{2\pi i} \int_{\Gamma_P} (\psi_Q)_-(s) ds. \]
Similarly, we have
\[ \frac{3}{2\pi i} \int_{z_3}^{z_4} \psi_R(s) ds =
   \frac{3}{2\pi i} \int_{\Gamma_R} (\psi_Q)_-(s) ds. \]
Therefore,
\[ \int d\mu_Q = \frac{3}{2\pi i} \oint_{\gamma} \psi_Q(s) ds \]
where $\gamma$ is the closed contour, consisting of the $-$sides of $\Gamma_P$
and $\Gamma_R$ and of the paths from $z_2$ to $z_3$ and from $z_4$ to $z_1$
within $\Gamma_Q$. Now we can use the residue theorem for the exterior
domain of $\gamma$. The only contribution comes from infinity, and
since $\psi_Q(s) = \frac{1}{3s} + \O(\frac{1}{s^2})$ as $s \to \infty$,
see (\ref{eq:Qinf}), we find $\int_{\Gamma_Q} d \mu_Q = 1$.
Hence $\mu_Q$ is a probability measure on $\Gamma_Q$.
\end{proof}

\subsection{The function $g_Q$} \label{subsec57}

For $g_Q$ we have properties that are analogous to the properties of
$g_P$ and $g_R$ obtained in Lemmas \ref{lem:gder} and \ref{explicit}.

\begin{lemma}  \label{lem:gQder}
For the derivative of  $g_Q$ we have
\begin{equation} \label{eq:gQder}
    g_Q'(z) =\left\{ \begin{array}{l}
        3 \psi_Q(z), \qquad z \in D_{\infty}, \\
        3 \psi_P(z), \qquad z \in D_P, \\
        3 \psi_R(z), \qquad z \in D_R.
        \end{array} \right.
\end{equation}
\end{lemma}
\begin{proof}
We have
\[ g_Q'(z) = \int \frac{1}{z-s} d\mu_Q(s). \]
As we did in the proof of Proposition \ref{prop:muQ} to evaluate
$\int d\mu_Q(s)$, we can write the integral as an integral
over the closed contour $\gamma$ that consists of the $-$sides of $\Gamma_P$
and $\Gamma_R$ and of the paths from $z_2$ to $z_3$ and from $z_4$ to $z_1$
within $\Gamma_Q$. We get for $z \in D_{\infty}$,
\[ g_Q'(z) = \frac{3}{2\pi i} \oint_{\gamma} \frac{\psi_Q(s)}{z-s} ds. \]
Since $\psi_Q$ is analytic in $D_{\infty}$ and vanishes at infinity, the
only contribution comes from the pole at $s=z$ and the result is
that $g_Q'(z) = 3 \psi_Q(z)$ if $z \in D_{\infty}$.
This gives (\ref{eq:gQder}) for $z \in D_{\infty}$.
The  expressions in the other two domains follow by analytic continuation.
\end{proof}

\begin{lemma} \label{lem:gQexplicit}
$g_Q$ has the explicit representation
\begin{equation} \label{eq:gQexplicit}
    g_Q(z) = 3zw - \log[w(w^2-1)] - 1 + \log \left(-\frac{1}{3}\right),
        \mbox{ with}
        \left\{ \begin{array}{ll}
         w = \psi_Q(z) & \mbox{ if } z \in D_{\infty}, \\
         w = \psi_P(z) & \mbox{ if } z \in D_P, \\
         w = \psi_R(z) & \mbox{ if } z \in D_R.
        \end{array} \right.
\end{equation}
\end{lemma}
\begin{proof}
This is analogous to the proof of Lemma \ref{explicit}.
\end{proof}

We conclude this section with  connections between $g_P$, $g_Q$ and
$g_R$ that easily follow from the explicit expressions
(\ref{eq:gPexplicit}), (\ref{eq:gRexplicit}), and (\ref{eq:gQexplicit}).

\begin{corollary} \label{corexplicit}
\begin{enumerate}
\item[\rm (a)]
For $z \in D_{\infty}$, we have
\begin{equation} \label{eq:sumgs}
    g_{P}(z) + g_Q(z) + g_R(z) = 3 \log z.
\end{equation}
\item[\rm (b)]
For $z\in D_P$, we have
\begin{equation} \label{relgPgQ}
    g_P(z)-g_Q(z)=3z+\ell.
\end{equation}
\item[\rm (c)]
For $z\in D_R$, we have
\begin{equation} \label{relgRgQ}
    g_R(z)-g_Q(z)=-3z+\ell.
\end{equation}
\end{enumerate}
\end{corollary}
\begin{proof}
Adding the formulas (\ref{eq:gPexplicit}), (\ref{eq:gRexplicit}), and (\ref{eq:gQexplicit})
with $z \in D_{\infty}$, we get
\begin{eqnarray} \nonumber
    \lefteqn{g_P(z) + g_Q(z) + g_R(z) =} \\
    && \label{eq:sumgshulp}
    3z(w_1 + w_2 + w_3) - \log\left[w_1(w_1^2-1) w_2(w_2^2-1) w_3(w_3^2-1) \right]
    -3 + \log\left(-\frac{4}{27}\right),
\end{eqnarray}
where $w_1 = \psi_P(z)$, $w_2= \psi_Q(z)$, and $w_3=\psi_R(z)$.
Since $w_1$, $w_2$, and $w_3$ are the solutions of the cubic
equation (\ref{eq:cubic}) we have
\[ z(w-w_1)(w-w_2)(w-w_3) = zw^3 - w^2 - zw + \frac{1}{3} \]
for every $w \in \mathbb C$. Comparing coefficients of $w^2$, we get
$z(w_1 + w_2 + w_3) = 1$, and taking $w=0$, $w=1$, and $w= -1$, we get
$z w_1 w_2 w_3 = - \frac{1}{3}$, $z (w_1-1)(w_2-1)(w_3-1) = \frac{2}{3}$,
and $z(w_1+1)(w_2+1)(w_3+1) = \frac{2}{3}$, respectively.
Using this in (\ref{eq:sumgshulp}), we arrive at (\ref{eq:sumgs}).

For parts (b) and (c), we start with the equation (\ref{eq:sumgs})
in $D_\infty$. Making use
of the jump relation (\ref{eq:gPRG1}) and  (\ref{eq:gPRG2}) to go
into $D_P$ and $D_R$, we obtain (\ref{relgPgQ}) and (\ref{relgRgQ}).
\end{proof}

\section{The Riemann--Hilbert problem and the first two transformations}

Our asymptotic analysis is based on the Riemann--Hilbert problem for
$Y$ formulated in the introduction, see (\ref{eq:Yjump}), (\ref{eq:Yasym}).
In this section we prove that the Riemann--Hilbert problem
has a unique solution and that the solution is given in terms
of the polynomials $P_n$, $Q_n$, and $R_n$, and the remainder $E_n$.
We also do the first two transformation of the Riemann--Hilbert problem,
which consist of a normalization of the problem at infinity, and
a deformation of contours.

\subsection{The Riemann--Hilbert problem}

We show that the Riemann--Hilbert problem for $Y$ has a solution
in terms of the Hermite-Pad\'e polynomials.

\begin{theorem}
Let $P_n$, $Q_n$, $R_n$, and $E_n$ be as above. Then the solution of the
Riemann--Hilbert problem for $Y$ (see the introduction)
is unique and is given by
\begin{equation} \label{eq:Yout}
 Y(z) =
         \begin{pmatrix}
         p_{n+1,n-1,n}(3nz) & z^{-3n-2} q_{n+1,n-1,n}(3nz) & r_{n+1,n-1,n}(3nz) \\[10pt]
         P_{n}(z) & z^{-3n-2} Q_{n}(z) & R_{n}(z) \\[10pt]
         p_{n,n-1,n+1}(3nz) & z^{-3n-2} q_{n,n-1,n+1}(3nz) & r_{n,n-1,n+1}(3nz)
         \end{pmatrix},
\end{equation}
for $z$ outside  $\Gamma$, and
\begin{equation} \label{eq:Yin}
 Y(z) =  \begin{pmatrix}
         p_{n+1,n-1,n}(3nz) &  z^{-3n-2}e_{n+1,n-1,n}(3nz) & r_{n+1,n-1,n}(3nz) \\[10pt]
         P_{n}(z) &  z^{-3n-2} E_{n}(z) & R_{n}(z) \\[10pt]
         p_{n,n-1,n+1}(3nz) &  z^{-3n-2}e_{n,n-1,n+1}(3nz) & r_{n,n-1,n+1}(3nz)
         \end{pmatrix},
\end{equation}
for  $z$  inside $\Gamma$.
In the first rows of (\ref{eq:Yout}) and (\ref{eq:Yin}) we use the
Hermite--Pad\'e polynomials of indices
$n+1,n-1,n$ normalized so that $p_{n+1,n-1,n}(3nz)$ is a monic polynomial,
and in the third row we use the Hermite--Pad\'e polynomials of indices
$n,n-1,n+1$ normalized so that $r_{n,n-1,n+1}(3nz)$ is monic.
\end{theorem}

\begin{proof}
The given $Y$ is analytic inside and outside the contour $\Gamma$. This is
clear from (\ref{eq:Yout}) and (\ref{eq:Yin}), except perhaps for the second column
of (\ref{eq:Yin}) which has a possible singularity at the origin. However,
the singularity is removable since $e_{n_1,n_2,n_3}(z) = \O(z^{n_1+n_2+n_3+2})$
as $z\to 0$.
The asymptotic condition (\ref{eq:Yasym}) is satisfied because of the
fact that $p_{n+1,n-1,n}(3nz)$, $Q_{n}(z)$,
and $r_{n,n-1,n+1}(3nz)$ are monic polynomials.

The jump condition can easily be checked. For the first and third entries in the
first row it reads
\[   (Y_{11})_+ = (Y_{11})_-, \qquad (Y_{13})_+ = (Y_{13})_- , \]
which is indeed so since $Y_{11}$ and $Y_{13}$ are both polynomials. We also have
\begin{eqnarray*}
    (Y_{12})_+(z) & = & z^{-3n-2} e_{n+1,n-1,n}(3nz) \\
    & = & z^{-3n-2} \left(p_{n+1,n-1,n}(3nz) e^{-3nz} + q_{n+1,n-1,n}(3nz) + r_{n+1,n-1,n}(3nz)e^{3nz}\right) \\
    & = & z^{-3n-2} e^{-3nz} (Y_{11})_-(z) + (Y_{12})_-(z) + z^{-3n-2} e^{3nz} (Y_{13})_-(z),
\end{eqnarray*}
and this is the jump condition (\ref{eq:Yjump}) for the second entry in the first row.
The second and third rows are handled in the same way.

To prove uniqueness, we assume that $\widetilde{Y}$ is another solution
of the Riemann--Hilbert problem. First observe that $\det Y$ is a
scalar function which is analytic in $\mathbb{C} \setminus \Gamma$.
Because of (\ref{eq:Yjump}), we have that $(\det Y)_+(z) = (\det Y)_-(z)$
for $z \in \Gamma$, so that
$\det Y$ has no jump, making $\det Y$ an entire function. For
large $z$ we have $\det Y(z) = 1 + \O(1/z)$ by (\ref{eq:Yasym}), hence by
Liouville's theorem $\det Y = 1$ everywhere. We can therefore consider
$\widetilde{Y} Y^{-1}$, which is analytic in $\mathbb{C} \setminus \Gamma$.
There is no jump on $\Gamma$ since $(\widetilde{Y} Y^{-1})_+(z)
= (\widetilde{Y} Y^{-1})_-(z)$ for every $z \in \Gamma$, hence
$\widetilde{Y} Y^{-1}$ is entire (i.e., each entry is an entire
function). For large $z$ we have $\widetilde{Y} Y^{-1}(z) =
I + \O(1/z)$, hence Liouville's theorem implies
that $\widetilde{Y} Y^{-1}(z) = I$ for every $z$, and hence
$\widetilde{Y}(z) = Y(z)$.
\end{proof}

\begin{remark}
A Riemann--Hilbert characterization for general Hermite--Pad\'e polynomials
was given in \cite{VAGeKu}.
\end{remark}

\subsection{First transformation}

We will use the functions $g_P$ and $g_R$, and the constant $\ell = \log 2 - \pi i$
from Section 4.3
to transform the Riemann--Hilbert problem for $Y$ to a Riemann--Hilbert problem for $U$, given by
\begin{equation}  \label{eq:U}
   U(z) = L^{-n-1} Y(z) \begin{pmatrix}
                      e^{-(n+1)g_P(z)} & 0 & 0 \\
                      0 & e^{(n+1)[g_P(z)+g_R(z)]} & 0 \\
                      0 & 0 & e^{-(n+1)g_R(z)}
                      \end{pmatrix}    L^{n+1},
\end{equation}
where $L$ is the constant diagonal matrix
\begin{equation}   \label{eq:L}
 L = \begin{pmatrix}
            e^{\ell/3} & 0 & 0 \\
             0 & e^{-2\ell/3} & 0 \\
            0 & 0 & e^{\ell/3}
          \end{pmatrix}.
\end{equation}
For the contour $\Gamma$ we take $\Gamma = \Gamma_P \cup \Gamma_R \cup \Gamma_U \cup \Gamma_L$, where
$\Gamma_U$ is a contour  connecting  $z_4$ to $z_1$ and lying in $D_{\infty,U}$, and
$\Gamma_L$ is a contour connecting $z_2$ to $z_3$ and lying in $D_{\infty,L}$. Then
$\Re \varphi_P > 0$ and $\Re \varphi_R > 0$ on $\Gamma_U \cup \Gamma_L$
by Lemma \ref{lemma:phiPR}.

We note that $U$ is analytic on $\mathbb C \setminus \Gamma$, since
$e^{g_P(z)}$ and $e^{g_R(z)}$ are analytic and single-valued on
$\mathbb{C} \setminus (\Gamma_P \cup \Gamma_R)$ and $\Gamma_P \cup \Gamma_R \subset \Gamma$.

Since $g_P(z) = \log z + \O(1/z)$ as $z \to \infty$,
we have $e^{(n+1)g_P(z)} = z^{n+1} [1 + \O(1/z)]$
as $z \to \infty$. Similarly we also have $e^{(n+1)g_R(z)} = z^{n+1} [1 + \O(1/z)]$
as $ z \to \infty$. Hence
\begin{equation}  \label{eq:Uasym}
    U(z) = I + \O\left(\frac1z\right), \qquad z \to \infty.
\end{equation}
So $U$ is normalized at infinity.

The jump relation for $U$ needs to be worked out on the four pieces of
the contour $\Gamma = \Gamma_P \cup \Gamma_R \cup \Gamma_U \cup \Gamma_L$.

For $z \in \Gamma_P$ we have
\begin{equation}  \label{eq:UjumpP1}
    U_+(z) = U_-(z) \begin{pmatrix}
      e^{-(n+1)[g_{P+}(z)-g_{P-}(z)]} &
        z e^{3z} e^{(n+1)[-3\log z - 3z +g_{P+}(z)+g_{P-}(z)+g_R(z)-\ell]} & 0 \\
      0 &  e^{(n+1) [g_{P+}(z)-g_{P-}(z)]} & 0 \\
      0 & z e^{-3z} e^{(n+1)[-3\log z +3z + g_{P+}(z)+2g_R(z)-\ell]} & 1
    \end{pmatrix}.
\end{equation}
Taking into account Lemma \ref{lem:gPR}, we can simplify the jump (\ref{eq:UjumpP1})
to
\[ U_+(z) = U_-(z) \begin{pmatrix}
                   e^{2(n+1)\varphi_{P+}(z)} & z e^{3z} & 0 \\
                   0 & e^{2(n+1)\varphi_{P-}(z)} & 0 \\
                   0 & ze^{-3z} e^{-2(n+1)\varphi_R(z)} & 1
                   \end{pmatrix}, \qquad z \in \Gamma_P.
\]
Similarly, the jump on $\Gamma_R$ is
\[  U_+(z) = U_-(z) \begin{pmatrix}
                   1 & ze^{3z} e^{-2(n+1)\varphi_P(z)} &  0 \\
                   0 & e^{2(n+1)\varphi_{R-}(z)} & 0 \\
                   0 & ze^{-3z} & e^{2(n+1)\varphi_{R+}(z)}
                   \end{pmatrix}, \qquad z \in \Gamma_R.
\]
On the parts $\Gamma_U$ and $\Gamma_L$ we have
\begin{equation}    \label{eq:UjumpUL1}
    U_+(z) = U_-(z) \begin{pmatrix}
        1 & ze^{3z} e^{(n+1)[-3\log z -3z + 2g_P(z)+g_R(z)-\ell]} & 0 \\
        0 & 1 & 0 \\
        0 & ze^{-3z} e^{(n+1)[-3\log z +3z + g_P(z)+2g_R(z)-\ell]} & 1
        \end{pmatrix}.
 \end{equation}
If we use Lemma \ref{lem:gPR} then (\ref{eq:UjumpUL1})
can be re-written as
\[ U_+(z) = U_-(z) \begin{pmatrix}
                   1 & ze^{3z} e^{-2(n+1)\varphi_P(z)} & 0 \\
                   0 & 1 & 0 \\
                   0 & ze^{-3z} e^{-2(n+1)\varphi_R(z)} & 1
                   \end{pmatrix}, \qquad z \in \Gamma_U \cup \Gamma_L.
\]

Summarizing, we have  the
following Riemann--Hilbert problem for $U$
\begin{enumerate}
\item $U$ is analytic on $\mathbb C \setminus \Gamma$.
\item $U$ satisfies the following jump relations
\begin{equation}  \label{eq:UjumpP}
  U_+(z) = U_-(z) \begin{pmatrix}
                   e^{2(n+1)\varphi_{P+}(z)} & z e^{3z} & 0 \\
                   0 & e^{2(n+1)\varphi_{P-}(z)} & 0 \\
                   0 & ze^{-3z} e^{-2(n+1)\varphi_R(z)} & 1
                   \end{pmatrix}, \qquad z \in \Gamma_P,
\end{equation}
\begin{equation}  \label{eq:UjumpR}
  U_+(z) = U_-(z) \begin{pmatrix}
                   1 & ze^{3z} e^{-2(n+1)\varphi_P(z)} &  0 \\
                   0 & e^{2(n+1)\varphi_{R-}(z)} & 0 \\
                   0 & ze^{-3z} & e^{2(n+1)\varphi_{R+}(z)}
                   \end{pmatrix}, \qquad z \in \Gamma_R,
\end{equation}
\begin{equation}  \label{eq:UjumpUL}
   U_+(z) = U_-(z) \begin{pmatrix}
                   1 & ze^{3z} e^{-2(n+1)\varphi_P(z)} & 0 \\
                   0 & 1 & 0 \\
                   0 & ze^{-3z} e^{-2(n+1)\varphi_R(z)} & 1
                   \end{pmatrix}, \qquad z \in \Gamma_U \cup \Gamma_L.
\end{equation}
\item $U(z) = I + \O\left(\frac1z\right)$ as $z \to \infty$.
\end{enumerate}

In Figure \ref{fig:Re} the shaded region is where $\Re \varphi_P$ is negative. For $\Re \varphi_R$ there
is a similar picture, but reflected along the imaginary axis. The contours $\Gamma_U$ and $\Gamma_L$
are in the region where both $\Re \varphi_P$ and $\Re \varphi_R$ are positive.
The jump matrix in (\ref{eq:UjumpUL}) for $U$ on the contours $\Gamma_U$ and $\Gamma_L$
is then the identity matrix $I$ plus a matrix
with entries that tend to zero exponentially fast as $n \to \infty$.

The factor $e^{-2(n+1)\varphi_R(z)}$ in the last row of the jump matrix in  (\ref{eq:UjumpP})
tends exponentially fast to $0$ because $\Re \varphi_R >0$ on $\Gamma_P$, see
Lemma \ref{lem:Revarphi}. In a similar way the factor $e^{-2(n+1)\varphi_P(z)}$
in the first row of the jump matrix in (\ref{eq:UjumpR}) on $\Gamma_R$ tends to $0$ exponentially fast
since $\Re \varphi_P > 0$ on $\Gamma_R$ (see Lemma \ref{lem:Revarphi} and  Figure \ref{fig:Re}).
Furthermore $\varphi_{P+} = -\varphi_{P-}$ is purely imaginary on $\Gamma_P$
because of (\ref{eq:GammaP})
and $\varphi_{R+}=-\varphi_{R-}$ is purely imaginary on $\Gamma_R$ because of (\ref{eq:GammaR}),
so that the diagonal elements of the jump matrices on $\Gamma_P$ and $\Gamma_R$ are oscillatory.

\subsection{Deformation of contours}
\label{deform}
The jump matrix in (\ref{eq:UjumpP}) can be written as a product of four matrices
\begin{multline}   \label{eq:productP}
\begin{pmatrix}
     e^{2(n+1)\varphi_{P+}(z)} & ze^{3z} & 0 \\
     0 & e^{2(n+1)\varphi_{P-}(z)} & 0 \\
     0 & ze^{-3z} e^{-2(n+1)\varphi_R(z)} & 1
    \end{pmatrix} \\ =
    \begin{pmatrix}
     1 & 0 & 0 \\
     z^{-1} e^{-3z} e^{2(n+1)\varphi_{P-}(z)} & 1 & 0 \\
     0 & 0 & 1
     \end{pmatrix}
     \begin{pmatrix}
     0 & ze^{3z} & 0 \\
     -z^{-1} e^{-3z} & 0 & 0 \\
     0 & 0 & 1
     \end{pmatrix} \\
     \begin{pmatrix}
     1 & 0 & 0 \\
     z^{-1} e^{-3z} e^{2(n+1)\varphi_{P+}(z)} & 1 & 0 \\
     0 & 0 & 1
     \end{pmatrix}
     \begin{pmatrix}
      1 & 0 & 0 \\
      0 & 1 & 0 \\
      0 & z e^{-3z} e^{-2(n+1)\varphi_R(z)} & 1
     \end{pmatrix}.
\end{multline}
Instead of jumping over $\Gamma_P$ in one jump, we will make four smaller jumps, and rather
than jumping over one contour, we jump over four contours, and each contour deals with
one of the matrices in the product (\ref{eq:productP}).
We will open up a lens around $\Gamma_P$ and introduce an extra contour
for the last factor of the matrix. The lens consists of two contours $\Gamma_{P^-} \cup \Gamma_{P^+}$
connecting $z_1$ and $z_2$, such that $\Gamma_{P^-}$ is on the minus side of $\Gamma_P$ and
$\Gamma_{P^+}$ is on the plus side of $\Gamma_P$, but still inside the region where $\Re \varphi_R > 0$
and $\Re \varphi_P < 0$. The contour $\Gamma_{P^{++}}$ stays away from $z_1$ and $z_2$.
It starts at a point on the upper contour $\Gamma_U$ and connects with a point on
the lower contour $\Gamma_L$, so that it lies inside the region
where $\Re \varphi_R > 0$, (so it stays to the left of $\Gamma_R^*$),
and does not intersect $\Gamma_{P^+}$. The subarc of $\Gamma_U$ between $z_1$
and the starting point of $\Gamma_{P^{++}}$ we call $\Gamma_{U^p}$, and the
subarc of $\Gamma_L$ from $z_2$ to the endpoint of $\Gamma_{P^{++}}$ is
called $\Gamma_{L^p}$.

\begin{figure}[ht]
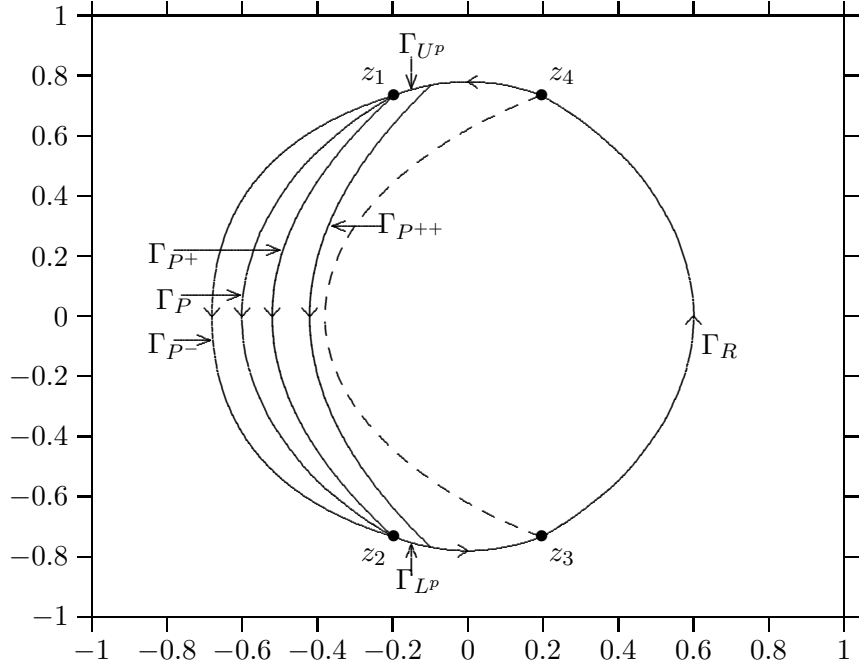

\hfil
\input traject2
\hfil
\caption{Deformation of contours around $\Gamma_P$. The dotted line is $\Gamma_R^*$}
\label{fig:deform}
\end{figure}

These contours are drawn in Figure \ref{fig:deform}.

In a similar way the jump (\ref{eq:UjumpR}) over $\Gamma_R$ has a factorization
\begin{multline}   \label{eq:productR}
\begin{pmatrix}
     1 & ze^{3z}e^{-2(n+1)\varphi_P(z)}  & 0 \\
     0 & e^{2(n+1)\varphi_{R-}(z)} & 0 \\
     0 & ze^{3z} & e^{2(n+1)\varphi_{R+}(z)}
    \end{pmatrix} \\ =
    \begin{pmatrix}
     1 & 0 & 0 \\
     0 & 1 & z^{-1}e^{3z}e^{2(n+1)\varphi_{R-}(z)} \\
     0 & 0 & 1
     \end{pmatrix}
     \begin{pmatrix}
     1 & 0 &  0 \\
     0 & 0 & -z^{-1}e^{3z}   \\
     0 & ze^{-3z} & 0
     \end{pmatrix} \\
     \begin{pmatrix}
     1 & 0 & 0 \\
     0 & 1 & z^{-1}e^{3z}e^{2(n+1)\varphi_{R+}(z)} \\
     0 & 0 & 1
     \end{pmatrix}
     \begin{pmatrix}
      1 & ze^{3z}e^{-2(n+1)\varphi_P(z)} & 0 \\
      0 & 1 & 0 \\
      0 & 0 & 1
     \end{pmatrix}.
\end{multline}
We now open up contours $\Gamma_{R^-}$, $\Gamma_{R^+}$, and $\Gamma_{R^{++}}$ around $\Gamma_R$,
which are mirror images of the contours around $\Gamma_P$. The mirror images
of $\Gamma_{U^p}$ and $\Gamma_{L^p}$ are denoted by $\Gamma_{U^r}$ and $\Gamma_{L^r}$, respectively.
The remaining middle parts of $\Gamma_U$ and $\Gamma_L$ are denoted by $\Gamma_{U^m}$ and $\Gamma_{L^m}$,
respectively.

\begin{figure}[ht]
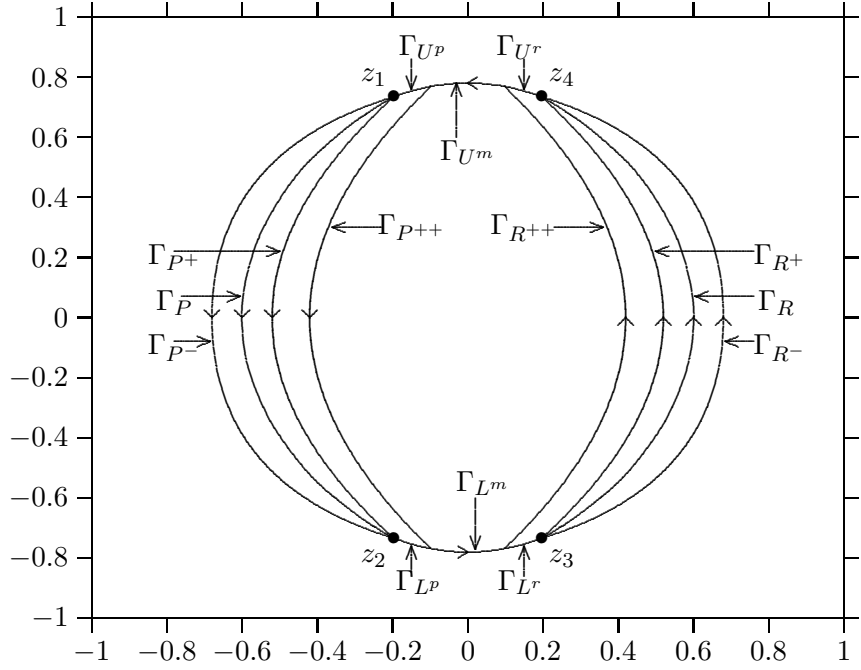

\hfil
\input traject3
\hfil
\caption{Contours of the RHP for $T$}  \label{fig:T}
\end{figure}

All these contours are shown in Figure \ref{fig:T}. All together there are 14 contours,
and they determine 8 regions in the plane.
The second transformation $U \mapsto T$ will be defined in each of these regions separately.
We define $T$ as follows.

We take
\begin{equation} \label{defT}
    T(z) = U(z),
\end{equation}
for $z$ in the unbounded region, and in the middle region bounded by
$\Gamma_{P^{++}}$, $\Gamma_{R^{++}}$, $\Gamma_{U^m}$, and $\Gamma_{L^m}$.
In the three regions near $\Gamma_P$ we put
\begin{equation} \label{defTnearGP}
   T(z) =  \begin{cases}
    U(z) V_{P^-}(z), &
    \mbox{for $z$ in the region bounded by $\Gamma_{P^-}$ and $\Gamma_P$}, \\
    U(z) V_{P^{++}}^{-1}(z) V_{P^+}^{-1}(z), &
    \mbox{for $z$ in the region bounded by $\Gamma_P$ and $\Gamma_{P^+}$}, \\
    U(z) V_{P^{++}}^{-1}(z), &
    \mbox{for $z$ in the region bounded by $\Gamma_{P^+}$, $\Gamma_{P^{++}}$,
    $\Gamma_{U^p}$, and $\Gamma_{L^p}$},
    \end{cases}
\end{equation}
where
\begin{eqnarray}
     V_{P^-}(z) & = & \begin{pmatrix}
                  1 & 0 & 0 \\
                  z^{-1}e^{-3z}e^{2(n+1)\varphi_P(z)} & 1 & 0 \\
                  0 & 0 & 1
                  \end{pmatrix}, \label{VP-}\\
     V_{P^+}(z) & = & \begin{pmatrix}
                   1 & 0 & 0 \\
                   z^{-1}e^{-3z}e^{2(n+1)\varphi_P(z)} & 1 & 0 \\
                   0 & 0 & 1
                   \end{pmatrix}, \\
     V_{P^{++}}(z) & = & \begin{pmatrix}
                   1 & 0 & 0 \\
                   0 & 1 & 0 \\
                   0 & ze^{-3z}e^{-2(n+1)\varphi_R(z)} & 1
                  \end{pmatrix}.
\end{eqnarray}
In the regions near $\Gamma_R$ we put similarly
\begin{equation} \label{defTnearGR}
    T(z) = \begin{cases}
    U(z) V_{R^-}(z), &
    \mbox{for $z$ in the region bounded by $\Gamma_{R^-}$ and $\Gamma_R$}, \\
    U(z) V_{R^{++}}^{-1}(z) V_{R^+}^{-1}(z), &
    \mbox{for $z$ in the region bounded by $\Gamma_R$ and $\Gamma_{R^+}$}, \\
    U(z) V_{R^{++}}^{-1}(z), &
    \mbox{for $z$ in the region bounded by $\Gamma_{R^+}$, $\Gamma_{R^{++}}$,
    $\Gamma_{U^r}$, and $\Gamma_{L^r}$},
    \end{cases}
\end{equation}
where
\begin{eqnarray}
     V_{R^-}(z) & = & \begin{pmatrix}
                  1 & 0 & 0 \\
                  0 & 1 & z^{-1}e^{3z}e^{2(n+1)\varphi_R(z)} \\
                  0 & 0 & 1
                  \end{pmatrix}, \\
     V_{R^+}(z) & = & \begin{pmatrix}
                   1 & 0 & 0 \\
                   0 & 1 & z^{-1}e^{3z}e^{2(n+1)\varphi_R(z)} \\
                   0 & 0 & 1
                   \end{pmatrix}, \\
     V_{R^{++}}(z) & = & \begin{pmatrix}
                   1 & ze^{3z}e^{-2(n+1)\varphi_P(z)} & 0 \\
                   0 & 1 & 0 \\
                   0 & 0 & 1
                  \end{pmatrix}.
\end{eqnarray}

Then we have the following Riemann--Hilbert problem for $T$.
\begin{enumerate}
\item $T$ is analytic in each of the 8 regions,
\item $T$ has a jump on each of the 14 contours
\[ T_+(z) = T_-(z) V_s(z), \qquad z \in \Gamma_s, \]
where $s$ stand for any of the symbols $P$, $P^-$, $P^+$, $P^{++}$, $R$, $R^-$, $R^+$, $R^{++}$,
$U^p$, $U^m$, $U^r$, $L^p$, $L^m$, or $L^r$.
The matrices $V_{P^-}$, $V_{P^+}$, $V_{P^{++}}$, $V_{R^-}$, $V_{R^+}$, and $V_{R^{++}}$ have
already been defined above. The other jump matrices are
\begin{eqnarray}
    V_P(z) & = & \begin{pmatrix}
            0 & ze^{3z} & 0 \\
            -z^{-1}e^{-3z} & 0 & 0 \\
            0 & 0 & 1
            \end{pmatrix}, \\
    V_R(z) & = & \begin{pmatrix}
            1 & 0 & 0 \\
            0 & 0 & -z^{-1}e^{3z} \\
            0 & ze^{-3z} & 0
            \end{pmatrix}, \\
    V_{U^p}(z) \, = \, V_{L^p}(z) & = & \begin{pmatrix}
                       1 & ze^{3z}e^{-2(n+1)\varphi_P(z)} & 0 \\
                       0 & 1 & 0 \\
                       0 & 0 & 1
                       \end{pmatrix}, \\
 V_{U^r}(z) \, = \, V_{L^r}(z) & = & \begin{pmatrix}
                       1 & 0 & 0  \\
                       0 & 1 & 0 \\
                       0 & ze^{-3z}e^{-2(n+1)\varphi_R(z)} & 1
                       \end{pmatrix}, \\
    V_{U^m}(z) \, = \, V_{L^m}(z) & = & \begin{pmatrix}
                       1 & ze^{3z}e^{-2(n+1)\varphi_P(z)} & 0  \\
                       0 & 1 & 0 \\
                       0 & ze^{-3z}e^{-2(n+1)\varphi_R(z)} & 1
                       \end{pmatrix}.
\end{eqnarray}
\item $T(z) = I + \O\left(\frac1z\right)$ as $z \to \infty$.
\end{enumerate}

Observe that all jumps, except for the jumps on $\Gamma_P$ and $\Gamma_R$, tend
to the identity matrix exponentially fast as $n \to \infty$.
Hence we expect that the dominating contributions are the jumps $V_P$ on $\Gamma_P$
and $V_R$ on $\Gamma_R$.

\section{Construction of parametrices and final transformation}

\subsection{Parametrix for the exterior region}
\begin{figure}[ht]
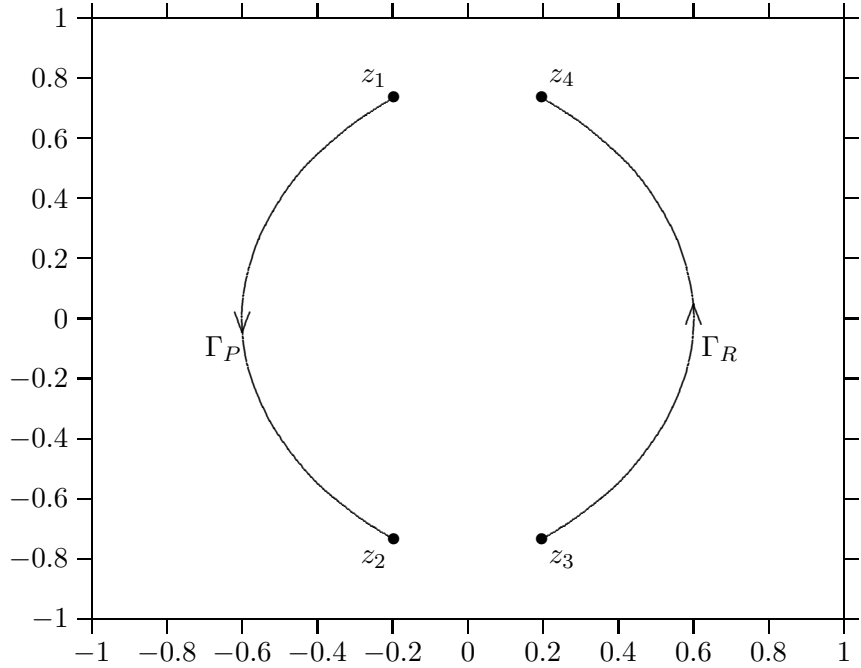

\hfil
\input traject4
\hfil
\caption{Contours of the RHP for $N$}  \label{fig:N}
\end{figure}

We will now solve a Riemann--Hilbert problem for a matrix valued function $N$
 on the contours $\Gamma_P \cup \Gamma_R$ (see Figure \ref{fig:N}) which,
in view of what was said at the end of the previous section,
is expected to describe the main contribution
of the Riemann--Hilbert problem of $T$.

We look for $N: \mathbb{C}\setminus (\Gamma_P \cup \Gamma_R)
 \to \mathbb{C}^{3 \times 3}$ satisfying
\begin{enumerate}
\item $N$ is analytic in $\mathbb{C} \setminus (\Gamma_P \cup \Gamma_R)$.
\item $N$ has jumps on $\Gamma_P$ and $\Gamma_R$ given by
\begin{equation}   \label{eq:NjumpP}
   N_+(z) = N_-(z) \begin{pmatrix}
                    0 & ze^{3z} & 0 \\
                    -z^{-1}e^{-3z} & 0 & 0 \\
                    0 & 0 & 1
                    \end{pmatrix}, \qquad z \in \Gamma_P,
\end{equation}
and
\begin{equation}   \label{eq:NjumpR}
   N_+(z) = N_-(z) \begin{pmatrix}
                    1 & 0 & 0 \\
                    0 & 0 & -z^{-1}e^{3z} \\
                    0 & ze^{-3z} & 0
                    \end{pmatrix}, \qquad z \in \Gamma_R.
\end{equation}
\item $N(z) = I + \O\left(\frac1z\right)$ as $z \to \infty$.
\end{enumerate}

\begin{proposition}  \label{lem:N}
A solution of the Riemann--Hilbert problem for $N$ is given by
\begin{equation}  \label{eq:NF}
  N(z) = \begin{pmatrix}
          F_1(\psi_P(z)) & F_1(\psi_Q(z)) & F_1(\psi_R(z)) \\
          F_2(\psi_P(z)) & F_2(\psi_Q(z)) & F_2(\psi_R(z)) \\
          F_3(\psi_P(z)) & F_3(\psi_Q(z)) & F_3(\psi_R(z))
          \end{pmatrix},
\end{equation}
where
\begin{eqnarray}
   F_1(w) & = &  \frac{-w(w-1)G(w)}{\sqrt{3 w^4 + 1}}, \label{eq:F1} \\
   F_2(w) & = & \frac{3(w^2-1)G(w)}{\sqrt{3 w^4 + 1}}, \label{eq:F2} \\
   F_3(w) & = & \frac{w(w+1)G(w)}{\sqrt{3 w^4 + 1}}, \label{eq:F3}
\end{eqnarray}
with $\sqrt{3w^4 + 1}$ defined and analytic in
$\mathbb C \setminus (\psi_{P+}(\Gamma_P) \cup \psi_{R+}(\Gamma_R))$, and such that
it is positive for large positive $w$.
The function $G$ is defined by
\begin{equation} \label{eq:G}
    G(w)= \left\{ \begin{array}{ll}
    we^{-\frac{(w+1)(2w-1)}{w(w-1)}} & \mbox{ for } w\in\psi(\mathcal{R}_P), \\[10pt]
    \left(\frac{w^2-1/3}{w^2-1}\right)e^{-\frac{2w^2}{w^2-1}} & \mbox{ for } w \in \psi(\mathcal{R}_Q), \\[10pt]
    we^{-\frac{(w-1)(2w+1)}{w(w+1)}} & \mbox{ for } w\in\psi(\mathcal{R}_R).
    \end{array} \right.
\end{equation}
\end{proposition}
\begin{proof}
Let us consider the first row $(N_{11}, N_{12}, N_{13})$ of $N$.
From (\ref{eq:NjumpP}) we get the following jumps on $\Gamma_P$
\begin{equation}  \label{eq:NP}
   \left\{ \begin{array}{lll}
   (N_{11})_+(z) & = & -z^{-1} e^{-3z} (N_{12})_-(z), \\
   (N_{12})_+(z) & = & ze^{3z} (N_{11})_-(z), \\
   (N_{13})_+(z) & = &  (N_{13})_-(z),
   \end{array} \right. \quad z \in \Gamma_P,
\end{equation}
and from (\ref{eq:NjumpR}) the following jumps on $\Gamma_R$
\begin{equation}  \label{eq:NR}
    \left\{ \begin{array}{lll}
     (N_{11})_+(z) & = & (N_{11})_-(z), \\
     (N_{12})_+(z) & = & z e^{-3z} (N_{13})_-(z), \\
     (N_{13})_+(z) & = &  -z^{-1} e^{3z} (N_{12})_-(z),
     \end{array} \right. \quad z \in \Gamma_R.
\end{equation}
Clearly $N_{11}$ is analytic on $\Gamma_R$ and $N_{13}$ is analytic on $\Gamma_P$.
So we can (and do) see $N_{11}$ as a function on the sheet $\mathcal R_P$ of
the Riemann surface $\mathcal R$, $N_{12}$ as a function on $\mathcal R_Q$
and $N_{13}$ as a function on $\mathcal R_R$.
Then we transform the problem from  $\mathcal{R}$ with the variable $z$,
to the complex $w$-plane, via the mapping $\psi: \mathcal{R} \to \mathbb{C}$.
The variables $z$ and $w$ are connected by (\ref{eq:zw}). The images
of the three sheets, the images of the branch points $z_1$, $z_2$, $z_3$, $z_4$,
and the images of the cuts $\Gamma_P$ and $\Gamma_R$ are shown in Figure \ref{fig:psi}.

Note that the images of $\Gamma_P$ under the mappings $\psi_{P+}$ and $\psi_{P-}$
(positive and negative boundary values of $\psi_P$ on $\Gamma_P$) give two arcs
from $w_1$ to $w_2$. They are oriented as shown in Figure \ref{fig:psi}.
The orientation corresponds to the orientation of $\Gamma_P$.
Together the arcs make up a simple closed loop around $-1$.

Similar remarks hold for the two images $\psi_{R+}(\Gamma_R)$ and $\psi_{R-}(\Gamma_R)$
of $\Gamma_R$. They make up a simple closed loop around $1$.

Now we transplant the (as yet unknown) functions $N_{11}$, $N_{12}$, and $N_{13}$ from
the Riemann surface to the $w$-plane, by defining $F_1$ as follows.
\begin{equation} \label{eq:defF1}
     F_1(w)  = \begin{cases}
   N_{11}\left( \frac{w^2-\frac13}{w(w^2-1)} \right), & \quad w \in \psi(\mathcal{R}_P), \\
   N_{12}\left( \frac{w^2-\frac13}{w(w^2-1)} \right), & \quad w \in \psi(\mathcal{R}_Q), \\
   N_{13}\left( \frac{w^2-\frac13}{w(w^2-1)} \right), & \quad w \in \psi(\mathcal{R}_R).
      \end{cases}
\end{equation}
Then $F_1$ is analytic in $\mathbb{C} \setminus (\psi_{P\pm}(\Gamma_P) \cup \psi_{R\pm}(\Gamma_R))$.
The jumps that $F_1$ should satisfy can be determined from (\ref{eq:NP})--(\ref{eq:NR}) and
are given by
\begin{equation} \label{eq:F1jumps}
    \begin{cases}
    F_{1+}(w) = z e^{3z} F_{1-}(w), & \quad w \in \psi_{P-}(\Gamma_P), \\
    F_{1+}(w) = -z^{-1} e^{-3z} F_{1-}(w), & \quad w \in \psi_{P+}(\Gamma_P), \\
    F_{1+}(w) = -z^{-1} e^{3z} F_{1-}(w), & \quad w \in \psi_{R+}(\Gamma_R), \\
    F_{1+}(w) = z e^{-3z} F_{1-}(w), & \quad w \in \psi_{R-}(\Gamma_R),
    \end{cases}
\end{equation}
where $z = z(w) = \frac{w^2-\frac{1}{3}}{w(w^2-1)}$.

The asymptotic condition on $N$ implies that $N_{11}(z) \to 1$, $N_{12}(z) \to 0$,
$N_{13}(z) \to 0$ as $z \to \infty$. For $F_1$, this means that
\begin{equation} \label{eq:F1asym}
    F_1(-1) = 1, \qquad F_1(0) = 0, \qquad F_1(1) = 0.
\end{equation}
We also want $F_1(w)$ to have a finite limit as $w \to \infty$, since $w = \infty$ corresponds to $z = 0$
on the $Q$-sheet.

We now seek $F_1$ in the form
\begin{equation} \label{eq:F1form}
    F_1(w)= \frac{-w(w-1) G(w)}{\sqrt{3w^4+1}}.
\end{equation}
Then $G$ should be analytic in $\mathbb{C} \setminus (\psi_{P\pm}(\Gamma_P) \cup \psi_{R\pm}(\Gamma_R))$
with jumps
\begin{equation} \label{eq:jumpG}
    \begin{cases}
    G_+(w) = z e^{3z} G_-(w), & \quad w \in \psi_{P-}(\Gamma_P), \\
    G_+(w) = z^{-1} e^{-3z} G_-(w), & \quad w \in \psi_{P+}(\Gamma_P), \\
    G_+(w) = z^{-1} e^{3z} G_-(w), & \quad w \in \psi_{R+}(\Gamma_R), \\
    G_+(w) = z e^{-3z} G_-(w), & \quad w \in \psi_{R-}(\Gamma_R),
    \end{cases}
\end{equation}
with $z = z(w)$.
The normalization for $G$ is
\begin{equation} \label{eq:Gin-1}
    G(-1) = -1.
\end{equation}
It is straightforward to check that $G$ given by (\ref{eq:G}) indeed satisfies
(\ref{eq:jumpG}) and (\ref{eq:Gin-1}). Then by (\ref{eq:F1form}) it follows that
$F_1$ has the correct jumps (\ref{eq:F1jumps}) and normalization (\ref{eq:F1asym}).
Then from (\ref{eq:defF1}) we recover $N_{11}$, $N_{12}$, and $N_{13}$ in
terms of $F_1$ by
\[ N_{11}(z) = F_1(\psi_P(z)), \quad N_{12}(z) = F_1(\psi_Q(z)),
    \quad N_{13}(z) = F_1(\psi_R(z)). \]
Then the jumps (\ref{eq:NP}) and (\ref{eq:NR}) are satisfied, and in addition
the normalization at infinity is correct. So we have found the first row of $N$.

The proof for the second and third rows is similar. The only difference
is that we have a different normalization at infinity, which leads to the
construction of functions $F_2$ and $F_3$ that satisfy the same jumps (\ref{eq:F1jumps})
as $F_1$, but are normalized by
\[ F_2(-1) = 0, \qquad F_2(0) = 1, \qquad F_2(1) = 0, \]
and
\[ F_3(-1) = 0, \qquad F_3(0) = 0, \qquad F_3(1) = 1. \]
Similar calculations then lead to the formulas (\ref{eq:F2}) and (\ref{eq:F3})
with the same function $G$.
\end{proof}

We remark that the entries of $N$ have fourth root singularities at the
branch points $z_k$ ($k=1,\ldots,4$). More precisely,
\begin{lemma} \label{lemmaOconditions}
The entries of $N$ behave as follows near the branch points.
As $z \to z_j$ with $j=1,2$, we have
\begin{equation} \label{Nnearz1z2}
    \begin{cases}
    N_{k1}(z) & = \O\left( |z-z_j|^{-1/4} \right), \\
    N_{k2}(z) & = \O\left( |z-z_j|^{-1/4} \right), \\
    N_{k3}(z) & = \O\left( 1 \right),
   \end{cases}
   \qquad k=1,2,3,
\end{equation}
and as $z \to z_j$ with $j=3,4$, we have
\begin{equation} \label{Nnearz3z4}
    \begin{cases}
    N_{k1}(z) & = \O\left( 1 \right), \\
    N_{k2}(z) & = \O\left( |z-z_j|^{-1/4} \right), \\
    N_{k3}(z) & = \O\left( |z-z_j|^{-1/4} \right),
   \end{cases}
   \qquad k=1,2,3.
\end{equation}
\end{lemma}
\begin{proof}
Since  $w_1$ is a non-degenerate critical point of the mapping $z = z(w)$,
we have for the inverse $w = \psi_P(z)$  as $z \to z_1 = z(w_1)$,
\begin{equation} \label{psiProndz1}
     \psi_P(z) =  w_1 + c (z-z_1)^{1/2} + \O\left(z-z_1\right)
\end{equation}
where $c$ is a non-zero constant.
Since $w_1$ is a simple root of $3 w^4 + 1$, it then follows that
\begin{equation}
    \sqrt{3\psi_P^4(z) + 1}  =  c_2 (z-z_1)^{1/4} + \O\left(|z-z_1|^{3/4}\right), \qquad z \to z_1,
\end{equation}
with $c_2 \neq 0$.
Since the  numerators of $F_1$, $F_2$, $F_3$ as given by (\ref{eq:F1})--(\ref{eq:F3}),
do not vanish for $w = w_1$, we find
\[ N_{k1}(z) = F_k(\psi_P(z)) = \O\left(|z-z_1|^{-1/4}\right), \qquad
    k=1,2,3,
\]
as $z \to z_1$.
In a similar way we find that $N_{k2}(z) = \O\left(|z-z_1|^{-1/4}\right)$
as $z \to z_1$.
The entries in the third column of $N$ are  analytic at $z_1$, since
$\psi_R$ is analytic at $z_1$ and the functions $F_1$, $F_2$, $F_3$ are
analytic at $\psi_R(z_1)$. This proves (\ref{Nnearz1z2}) for $j=1$.

The behavior near the other branch points follows in a similar way.
\end{proof}
\begin{remark}
It will be useful to have another representation for the entries in the second
row of $N$. They are
\begin{equation} \label{eq:N21N22N23}
    N_{21}(z) = \frac{2 e^{-g_P(z)}}{\sqrt{3 \psi_P^4(z) + 1}}, \quad
   N_{22}(z) = - \frac{e^{g_P(z) + g_R(z)}}{z^2 \sqrt{3 \psi_Q^4(z) + 1}}, \quad
   N_{23}(z) = \frac{2 e^{-g_R(z)}}{\sqrt{3\psi_R^4(z) + 1}}.
\end{equation}
It may be checked directly that these functions have the right asymptotics
as $z \to \infty$, and satisfy the correct jump relations on $\Gamma_P$
and $\Gamma_R$. They also satisfy the $\cal O$-conditions of
Lemma \ref{lemmaOconditions}.
\end{remark}

The Riemann--Hilbert problem for $T$ is now very close to the Riemann--Hilbert problem
for $N$ because the jumps for $T$ and $N$ on the contours $\Gamma_P$ and $\Gamma_R$
are the same and the jumps for $T$ on the other contours tend to the identity matrix
as $n \to \infty$, uniformly away from the branch points.
So we expect that $T$ behaves like $N$ as $n \to \infty$ away from the branch
points. However, in order to justify this we need to analyze the Riemann--Hilbert problem
in more detail near the branch points.

\subsection{Parametrices near the branch points}
\label{sec-branch}
\begin{figure}[ht]
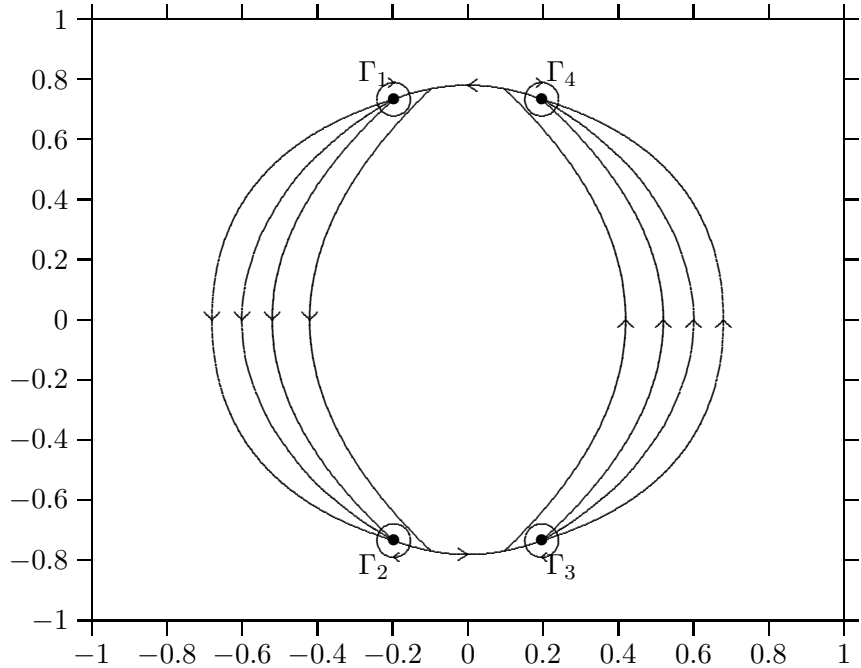

\hfil
\input traject5
\hfil
\caption{Small circles $\Gamma_1$, $\Gamma_2$, $\Gamma_3$, $\Gamma_4$ around the branch points}  \label{fig:t5}
\end{figure}

We introduce four new contours $\Gamma_1$, $\Gamma_2$, $\Gamma_3$, $\Gamma_4$ which are
small circles of radius $\delta$ centered at the branch points, as indicated in Figure \ref{fig:t5}.
We choose $\delta$ small enough so that the circles do not intersect the curves
$\Gamma_{P^{++}}$ and $\Gamma_{R^{++}}$.
The circles are oriented clockwise.
Inside each of these contours we will solve the Riemann--Hilbert
problem for $T$ exactly. The analysis is similar for the four branch points, and we will
only work out the analysis near $z_1$ in detail.

\begin{figure}[ht]
\centering
\unitlength=1.00mm
\begin{picture}(70.00,75.00)(20,70)
\bezier{224}(90.00,110.00)(88.00,82.00)(60.00,80.00)
\bezier{224}(60.00,80.00)(32.00,82.00)(30.00,110.00)
\bezier{224}(60.00,140.00)(88.00,138.00)(90.00,110.00)
\bezier{224}(30.00,110.00)(32.00,138.00)(60.00,140.00)
\bezier{120}(60.00,110.00)(47.00,111.00)(30.00,109.00)
\bezier{128}(60.00,110.00)(47.00,97.00)(44.00,84.00)
\bezier{256}(34.00,95.00)(54.00,115.00)(90.00,113.00)
\put(60.00,110.00){\circle*{2.00}}
\put(44.00,112.00){\line(-2,-1){4.00}}
\put(40.00,110.00){\line(2,-1){4.00}}
\put(40.00,100.00){\line(1,2){2.00}}
\put(40.00,100.00){\line(1,0){4.00}}
\put(48.00,94.00){\line(0,1){4.00}}
\put(48.00,94.00){\line(3,2){4.00}}
\put(71.00,112.50){\line(4,3){4.00}}
\put(71.00,112.50){\line(2,-1){4.00}}
\put(60.00,115.00){\makebox(0,0)[cc]{$z_1$}}
\put(58.00,142.00){\line(2,-1){4.00}}
\put(62.00,140.00){\line(-2,-1){4.00}}
\put(35.00,136.00){\makebox(0,0)[cc]{$\Gamma_1$}}
\put(83.00,110.00){\makebox(0,0)[cc]{$\Gamma_{U^p}$}}
\put(48.00,88.00){\makebox(0,0)[cc]{$\Gamma_{P^+}$}}
\put(38.00,95.00){\makebox(0,0)[cc]{$\Gamma_P$}}
\put(34.00,113.00){\makebox(0,0)[cc]{$\Gamma_{P^-}$}}
\end{picture}
\caption{The contours near the branch point $z_1$} \label{fig:branch}
\end{figure}
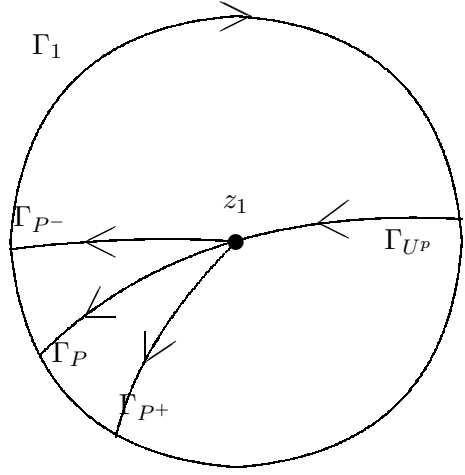

Zooming in near the branch point $z_1$ gives a Riemann--Hilbert problem with
five contours $\Gamma_{P_1}, \Gamma_{P}, \Gamma_{P^+}, \Gamma_U, \Gamma_1$, as
shown in Figure \ref{fig:branch}.

The jumps on these contours are
\begin{eqnarray*}
     V_{P^-}(z) \, = \, V_{P^+}(z)& = & \begin{pmatrix}
                  1 & 0 & 0 \\
                  z^{-1}e^{-3z} e^{2(n+1)\varphi_P(z)} & 1 & 0 \\
                  0 & 0 & 1
                  \end{pmatrix} \\
     V_P & = & \begin{pmatrix}
                 0 & z e^{3z} & 0 \\
                -z^{-1} e^{-3z} & 0 & 0 \\
                 0 & 0 & 1
                 \end{pmatrix} \\
     V_{U^p}(z) & = & \begin{pmatrix}
                       1 & ze^{3z} e^{-2(n+1)\varphi_P(z)} & 0 \\
                       0 & 1 & 0 \\
                       0 & 0 & 1
                       \end{pmatrix}.
\end{eqnarray*}

We look for a $3\times 3$ matrix valued function $M^{(1)}$ defined within the disk $\Delta_1$
surrounded by $\Gamma_1$, such that
\begin{enumerate}
\item $M^{(1)}$ is analytic in $\Delta_1 \setminus (\Gamma_P \cup \Gamma_{P^-} \cup \Gamma_{P^+} \cup \Gamma_{U^p})$,
\item $M^{(1)}$ has the jumps
    \[ M^{(1)}_+(z) = M^{(1)}_-(z) V_s(z), \qquad z \in \Gamma_s, \]
    where $s$ stands for any of the symbols $P$, $P^-$, $P^+$, and $U^p$.
\item On $\Gamma_1$ we have that $M^{(1)}$ matches $N$ in the sense that
    \[ M^{(1)}(z) = \left(I + \O\left(\frac{1}{n}\right)\right) N(z) \]
    uniformly for $z \in \Gamma_1$.
\end{enumerate}

The construction of $M^{(1)}$ is basically a $2\times 2$ problem, since
the jump matrices are non-trivial only in the $2\times 2$ left upper block.
The jumps have a standard form and a local parametrix can be built out
of Airy functions. To be more precise, we will find $M^{(1)}$ in the form
\begin{equation} \label{M1form}
    M^{(1)}(z) = E^{(1)}(z) \Psi^{(1)}((n+1)^{2/3} f_1(z))
    \begin{pmatrix}
     z^{-1/2} e^{-3z/2}  e^{(n+1)\varphi_P(z)} & 0 & 0 \\
                     0 & z^{1/2} e^{3z/2} e^{-(n+1)\varphi_P(z)} & 0 \\
                     0 & 0 & 1
                    \end{pmatrix}.
\end{equation}
Here we take the branch of $\varphi_P$ which is zero at $z_1$.
So for $z \in \Delta_1$, we take
\[ \varphi_P(z) = \frac{3}{2} \int_{z_1}^z (\psi_Q - \psi_P)(s) ds \]
with integration from $z_1$ to $z$ in $\Delta_1 \setminus \Gamma_P$.
Recall that $\varphi_P$ is defined modulo $\pi i$, and that so far
the precise branch did not matter, since we only had expressions like
$e^{2(n+1)\varphi_P(z)}$. Since
\[ \psi_P(s) = w_1 + c(s-z_1)^{1/2} + \O\left(s-z_1\right),
   \quad
   \psi_Q(s) = w_1 - c(s-z_1)^{1/2} + \O\left(s-z_1\right),
\]
as $s \to z_1$, with the same non-zero constant $c_1$,
see also (\ref{psiProndz1}), we then have from this definition of $\varphi_P$
that
\[ \varphi_P(z) = (z - z_1)^{3/2} h(z),     \qquad z \in \Delta_1 \setminus \Gamma_P, \]
with $h$ analytic and without zeros in $\Delta_1$. The function $(z-z_1)^{3/2}$ is
defined with a branch cut along $\Gamma_P$.

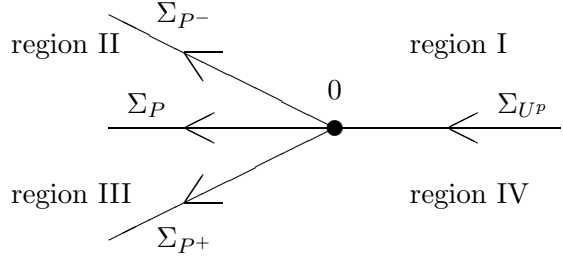
\begin{figure}[ht]
\centering
\unitlength=1.00mm
\begin{picture}(70.00,40.00)(20,90)
\put(30,110){\line(1,0){60}}
\put(60,110){\line(-2,-1){30}}
\put(60,110){\line(-2,1){30}}
\put(60.00,110.00){\circle*{2.00}}
\put(75.00,110.00){\line(2,1){4.00}}
\put(75.00,110.00){\line(2,-1){4.00}}
\put(40.00,110.00){\line(2,-1){4.00}}
\put(40.00,110.00){\line(2,1){4.00}}
\put(40.00,120.00){\line(1,0){5.00}}
\put(40.00,120.00){\line(2,-3){2.50}}
\put(40.00,100.00){\line(1,0){5.00}}
\put(40.00,100.00){\line(2,3){2.50}}
\put(60.00,115.00){\makebox(0,0)[cc]{$0$}}
\put(85.00,113.00){\makebox(0,0)[cc]{$\Sigma_{U^p}$}}
\put(40.00,125.00){\makebox(0,0)[cc]{$\Sigma_{P^-}$}}
\put(35.00,113.00){\makebox(0,0)[cc]{$\Sigma_P$}}
\put(40.00,95.00){\makebox(0,0)[cc]{$\Sigma_{P^+}$}}
\put(70,120){region I}
\put(17,120){region II}
\put(17,100){region III}
\put(70,100){region IV}
\end{picture}
\caption{The contours for $\Psi^{(1)}$}  \label{fig:Ai}
\end{figure}

The function $f_1(z)$ is defined such that
\[ \varphi_P(z) = \frac{2}{3} \left[ f_1(z) \right]^{3/2}. \]
It is a conformal map from $\Delta_1$ onto a convex neighborhood of $0$.
(We may have to shrink $\Delta_1$, if necessary.) It maps $\Gamma_P$ onto
a part of the negative real axis. We still have some freedom in the choice
of $\Gamma_{P^-}$, $\Gamma_{P^+}$ and $\Gamma_{U^p}$. We take $\Gamma_{U^p}$
so that $f_1$ maps $\Gamma_{U^p}$ to a part of the positive real line.
This means that $\Gamma_{U^p}$ is an analytic continuation of $\Gamma_P$.
$\Gamma_{P^-}$ and $\Gamma_{P^+}$ are chosen so that they are mapped
by $f_1$ onto rays in the complex $s$-plane.
We denote the images of $\Gamma_P$, $\Gamma_{P^-}$, $\Gamma_{P^+}$,
and $\Gamma_{U^p}$, by $\Sigma_{P}$, $\Sigma_{P^-}$, $\Sigma_{P^+}$
and $\Sigma_{U^p}$. These contours are shown in Figure \ref{fig:Ai}.
On these contours we use the constant jump matrices
\[ \widehat{V}_{P^-} = \widehat{V}_{P^+} =
    \begin{pmatrix} 1 & 0 & 0 \\ 1 & 1 & 0 \\ 0 & 0 & 1 \end{pmatrix} \]
on $\Sigma_{P^-}$ and $\Sigma_{P^+}$,
\[ \widehat{V}_P = \begin{pmatrix} 0 & 1 & 0 \\ -1 & 0 & 0 \\ 0 & 0 & 1 \end{pmatrix} \]
on $\Sigma_P$, and
\[ \widehat{V}_{U^p} = \begin{pmatrix} 1 & 1 & 0 \\ 0 & 1 & 0 \\ 0 & 0 & 1 \end{pmatrix} \]
on $\Sigma_{U^p}$. This Riemann--Hilbert problem is well-known when we
deal with $2 \times 2$ matrices (see, e.g., \cite[\S 7.6]{deift},
\cite{deiftal1}, \cite[p.~1522]{deiftal2}). Note however that  the arrows are pointing
in the opposite direction compared to Figures 7.18--7.20 in \cite{deift}, so that we
need to modify the solution a little.
The  solution $\Psi^{(1)}$ uses the Airy function $\Ai(z)$, which is the unique
solution of $y''(z)=zy(z)$ with asymptotics as $z \to \infty$ given by
\begin{eqnarray}
   \Ai(z) & = & \frac{1}{2\sqrt{\pi}} z^{-1/4} e^{-\frac23 z^{3/2}}
   \left( 1 + \O\left(\frac{1}{z^{3/2}} \right) \right), \label{eq:Aiasym} \\
   \Ai'(z) & = & \frac{-1}{2\sqrt{\pi}} z^{1/4} e^{-\frac23 z^{3/2}}
   \left( 1 + \O\left(\frac{1}{z^{3/2}} \right) \right), \label{eq:Aiderasym}
\end{eqnarray}
which holds for $|\arg z| < \pi$, where $z^{1/4}$ and $z^{3/2}$ are defined with
principal branch (i.e., with a cut along the negative real axis).
We will define $\Psi^{(1)}$ by
\begin{eqnarray*}
 \Psi^{(1)}(s) & = & \begin{pmatrix}
             \Ai(s) & -\Ai(\omega_3^2 s) & 0 \\
             \Ai'(s) & -\omega_3^2\Ai'(\omega_3^2 s) & 0 \\
              0 & 0 & 1
             \end{pmatrix}
             \begin{pmatrix}
             e^{-i\pi/6} & 0 & 0 \\
             0 & e^{i\pi/6} & 0 \\
             0 & 0 & 1
             \end{pmatrix}   , \qquad s \in I, \\
  \Psi^{(1)}(s) & = & \begin{pmatrix}
             \Ai(s) & -\Ai(\omega_3^2 s) & 0 \\
             \Ai'(s) & -\omega_3^2 \Ai'(\omega_3^2 s) & 0 \\
              0 & 0 & 1
             \end{pmatrix}
             \begin{pmatrix}
             e^{-i\pi/6} & 0 & 0 \\
             0 & e^{i\pi/6} & 0 \\
             0 & 0 & 1
             \end{pmatrix} \widehat{V}_{P^-}   , \qquad s \in II, \\
   \Psi^{(1)}(s) & = &  \begin{pmatrix}
             \Ai(s) & \omega_3^2\Ai(\omega_3 s) & 0 \\
             \Ai'(s) &  \Ai'(\omega_3 s) & 0 \\
              0 & 0 & 1
             \end{pmatrix}
             \begin{pmatrix}
             e^{-i\pi/6} & 0 & 0 \\
             0 & e^{i\pi/6} & 0 \\
             0 & 0 & 1
             \end{pmatrix}
             \widehat{V}_{P^+}^{-1}, \qquad s \in III, \\
   \Psi^{(1)}(s) & = & \begin{pmatrix}
             \Ai(s) & \omega_3^2\Ai(\omega_3 s) & 0 \\
             \Ai'(s) & \Ai'(\omega_3 s) & 0 \\
              0 & 0 & 1
             \end{pmatrix}
             \begin{pmatrix}
             e^{-i\pi/6} & 0 & 0 \\
             0 & e^{i\pi/6} & 0 \\
             0 & 0 & 1
             \end{pmatrix}, \qquad s \in IV,
\end{eqnarray*}
where $\omega_3 = e^{2\pi i/3}$ is a primitive third root of unity.

With the above definitions of $\Psi^{(1)}$ and $f_1$ it may then
be shown that for any analytic prefactor $E^{(1)}$ the matrix $M^{(1)}$
defined by (\ref{M1form}) satisfies the jump conditions on $\Gamma_s$,
where $s$ is any of the symbols
$P$, $P^{-}$,
$P^{+}$, and $U^p$. The extra factor $E^{(1)}$ has to
be chosen in such a way that $M^{(1)}$ satisfies the matching condition
on $\Gamma_1$ as well.

On the part of $\Gamma_1$ that lies in $\left(f_1\right)^{-1}(I)$
(the arc between $\Gamma_{P^-}$ and $\Gamma_{U^p}$) the asymptotic expansions for
the Airy functions (\ref{eq:Aiasym})--(\ref{eq:Aiderasym}) give
\begin{eqnarray*}
    \Ai((n+1)^{2/3}f_1 (z)) & = &
    \frac{1}{2\sqrt{\pi}} (n+1)^{-1/6} \left(f_1(z)\right)^{-\frac{1}{4}} e^{-(n+1)\varphi_P(z)}
   \left( 1 + \O\left(\frac1n\right) \right), \\
    \Ai((n+1)^{2/3}\omega_3^2 f_1(z)) & = &
    \frac{1}{2\sqrt{\pi}} (n+1)^{-1/6} \left(f_1(z)\right)^{-\frac{1}{4}} e^{i\pi/6} e^{(n+1)\varphi_P(z)}
   \left( 1 + \O\left(\frac1n\right) \right), \\
   \Ai'((n+1)^{2/3} f_1(z)) & = &
   \frac{-1}{2\sqrt{\pi}} (n+1)^{1/6} \left(f_1(z)\right)^{\frac{1}{4}} e^{-(n+1)\varphi_P(z)}
   \left( 1 + \O\left(\frac1n\right) \right), \\
   \Ai'((n+1)^{2/3}\omega_3^2 f_1(z)) & = &
   \frac{-1}{2\sqrt{\pi}} (n+1)^{1/6} \left(f_1(z)\right)^{\frac{1}{4}} e^{-i\pi/6} e^{(n+1)\varphi_P(z)}
   \left( 1 + \O\left(\frac1n\right) \right).
\end{eqnarray*}
Here the fourth root in  $\left(f_1(z) \right)^{\frac{1}{4}}$ is defined with
a cut along $\Gamma_P$.
This means that on this part of $\Gamma_1$ we have
\begin{multline}  \label{eq:Mpasym}
   M^{(1)}(z) = \frac{1}{2\sqrt{\pi}} E^{(1)}(z)
   \begin{pmatrix}
   ((n+1)^{2/3} f_1(z))^{-\frac14} & 0 & 0 \\
   0 & ((n+1)^{2/3} f_1(z))^{\frac14} & 0 \\
   0 & 0 & 1 \end{pmatrix} \\
   \begin{pmatrix}
   e^{-i\pi/6} & -e^{i\pi/3} & 0 \\
   -e^{-i\pi/6} & -e^{i\pi/3} & 0 \\
   0 & 0 & 1
   \end{pmatrix}
   \begin{pmatrix} z^{-1/2} e^{-3z/2} & 0 & 0 \\
   0 & z^{1/2} e^{3z/2} & 0 \\
   0 & 0 & 1 \end{pmatrix}
   \left( I + \O\left( \frac1n \right) \right).
\end{multline}
In order to have the matching condition we therefore choose
\begin{multline}
   E^{(1)}(z)  =  2\sqrt{\pi} N(z)
   \begin{pmatrix} z^{1/2} e^{3z/2} & 0 & 0 \\
   0 & z^{-1/2} e^{-3z/2} & 0 \\
   0 & 0 & 1 \end{pmatrix} \\
   \begin{pmatrix}
   e^{-i\pi/6} & -e^{i\pi/3} & 0 \\
   -e^{-i\pi/6} & -e^{i\pi/3} & 0 \\
   0 & 0 & 1
   \end{pmatrix}^{-1}
   \begin{pmatrix}
   ((n+1)^{2/3} f_1(z))^{\frac14} & 0 & 0 \\
   0 & ((n+1)^{2/3} f_1(z))^{-\frac14} & 0 \\
   0 & 0 & 1 \end{pmatrix} \\
    = \sqrt{\pi} e^{i\pi/6} N(z)
   \begin{pmatrix} z^{1/2} e^{3z/2} & 0 & 0 \\
   0 & z^{-1/2} e^{-3z/2} & 0 \\
   0 & 0 & 1 \end{pmatrix}\\
   \begin{pmatrix} 1 & - 1 & 0 \\
                  i & i & 0 \\
                  0 & 0 & 2 e^{-i\pi/6} \end{pmatrix}
   \begin{pmatrix} \\
   ((n+1)^{2/3} f_1(z))^{\frac14} & 0 & 0 \\
   0 & ((n+1)^{2/3} f_1(z))^{-\frac14} & 0 \\
   0 & 0 & 1 \end{pmatrix}.
      \label{eq:E}
\end{multline}
With this choice of $E^{(1)}$ it is then easy to see that the matching condition
holds on the part of $\Gamma_1$ in the region bounded by $\Gamma_{P^-}$ and
$\Gamma_{U^p}$. A similar analysis for the other parts of $\Gamma_1$ shows
that the same $E^{(1)}$ also works in the other regions.

From (\ref{eq:E}) it is easy to see that $E^{(1)}$ is analytic in $\Delta_1 \setminus \Gamma_P$.
On $\Gamma_P$, both $N$ and $\left(f_1\right)^{\frac{1}{4}}$ have a jump. $N$ has the
jump (\ref{eq:NjumpP}) and $\left(f_1\right)^{\frac{1}{4}}$ satisfies
$\left(f_1\right)^{\frac{1}{4}}_+ = -i \left(f_1\right)^{\frac{1}{4}}_-$.
Straightforward calculations then show that $E^{(1)}_+ = E^{(1)}_-$ on $\Gamma_P$,
so that $E^{(1)}$ is analytic across $\Gamma_P$. From (\ref{eq:E}) and
the fact that the entries of $N$ have at most fourth root singularities
at $z_1$, see (\ref{Nnearz1z2}), we see that the entries of $E^{(1)}$ have at most a square root
singularity at $z_1$. Since $E^{(1)}$ is analytic in $\Delta_1 \setminus \{z_1\}$,
the singularity at $z_1$ is removable, and this proves that $E^{(1)}$ is analytic
in the full $\Delta_1$. This completes the construction of the parametrix
$M^{(1)}$ in the neighborhood $\Delta_1$ of $z_1$.

From the definitions (\ref{M1form}) and (\ref{eq:E}) it is clear that
\[ M^{(1)} = N \begin{pmatrix} * & * & 0 \\ * & * & 0 \\ 0 & 0 & 1 \end{pmatrix}, \]
where $*$ denotes an unspecified entry. This means that the
third column of $M^{(1)}$ agrees with the third column of $N$.
We will use this in what follows.

In a similar way, we can construct parametrices $M^{(2)}$, $M^{(3)}$,
and $M^{(4)}$ near the other branch points $z_2$, $z_3$, $z_4$.
The third column of $M^{(2)}$ agrees with the third column of $N$,
and the first columns of $M^{(3)}$ and $M^{(4)}$ agree with the
first column of $N$.

\subsection{Third transformation}

\begin{figure}[ht]
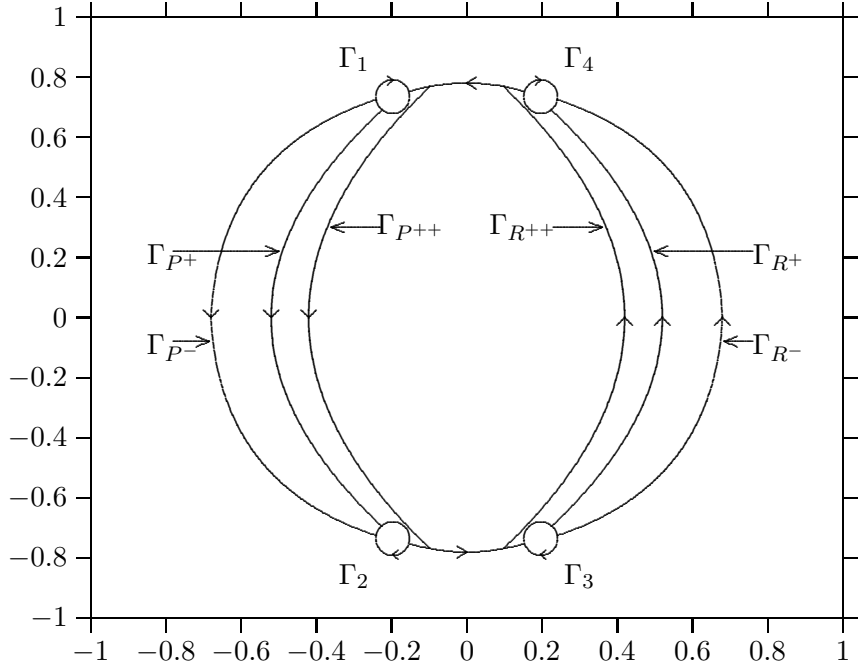

\hfil
\input traject6
\hfil
\caption{Contours of the RHP for $S$}  \label{fig:t6}
\end{figure}

We now introduce the final matrix
\begin{equation} \label{defS}
   S(z) = \begin{cases}
   T(z) \left(N(z)\right)^{-1},  &\qquad z \textrm{ outside } \Gamma_1, \Gamma_2, \Gamma_3, \Gamma_4, \\[10pt]
   T(z) \left(M^{(j)}(z) \right)^{-1}, & \qquad z \textrm{ inside } \Gamma_j,\ j=1,2,3,4.
   \end{cases}
\end{equation}
Inside each $\Gamma_j$ the matrices $T$ and $M^{(j)}$ have the same jumps, hence
$S$ has no jumps inside $\Gamma_j$.
Outside the $\Gamma_j$ the matrices $T$ and $N$ have the same jump matrices on
$\Gamma_P$ and $\Gamma_R$. Hence $S$ has no jump on $\Gamma_P$ and $\Gamma_R$.
This means that $S$ solves a Riemann--Hilbert problem on the system of curves shown in Figure \ref{fig:t6}.

$S$ is analytic outside the above system of contours and it is
normalized at infinity
\begin{equation}
\label{eq:Sasym}
S(z) = I + \O\left(\frac{1}{z}\right), \qquad z \to \infty.
\end{equation}
\begin{theorem}  \label{thm:T}
The matrix $S(z)$ has the behavior
\begin{equation}  \label{Sasym}
   S(z) = I + \O\left( \frac1n \right),
\end{equation}
uniformly on $\mathbb{C} \setminus \Sigma_S$, where
$\Sigma_S$ are the contours in Figure \ref{fig:t6}.
\end{theorem}

\begin{proof}
The jumps on all of the contours are uniformly of the form
$I +  \O(e^{-cn})$
with some fixed $c > 0$, except for the jumps on the circles $\Gamma_j$
where we have
\[ S_+(z) = S_-(z) M^{(j)}(z) N^{-1}(z), \qquad z \in \Gamma_j. \]
Because of the matching condition we have
\[ M^{(j)}(z) N^{-1}(z) = I + \O\left(\frac{1}{n}\right) \]
uniformly for $z \in \Gamma_j$.
Hence $S(z)$ solves a Riemann--Hilbert problem, normalized at $\infty$ with
jumps close to the identity matrix up to $\O(1/n)$, uniformly
on the contours $\Sigma_S$.
We can then use arguments as those leading to Theorem 7.171 in
\cite{deift} to obtain (\ref{Sasym}), see also \cite{deiftal1,deiftal2,Kuij}.
\end{proof}

Using the asymptotic expansion for the Airy function and its derivative, it is possible
to obtain a full asymptotic expansion for the jump matrix on the contours $\Gamma_j$
in powers of $1/n$. This in turn yields a full asymptotic expansion for $S$
\begin{equation} \label{Sasympexp}
     S(z) \sim I + \sum_{k=1}^{\infty} \frac{S_k(z)}{n^k}
\end{equation}
as in Theorem 7.10 of \cite{deiftal2}. We will not use the full expansion,
since (\ref{Sasym}) is enough for the asymptotic results
as stated in Section 2. It will be clear however that the use of
the full asymptotic expansion (\ref{Sasympexp}) will also lead to full
asymptotic expansions for the polynomials $P_n$, $Q_n$, and $R_n$, and for
the remainder $E_n$.

\section{Proofs of the asymptotic formulas}
We now know the asymptotic behavior (\ref{Sasym}) of $S$ as $n \to \infty$. We will trace back our steps
to the original Riemann--Hilbert problem for $Y$ to obtain asymptotics for the polynomials
in the Hermite--Pad\'e approximation problem for the exponential function.

\subsection{Proof of Theorem \ref{asymp}}
\begin{proof}
We start with the proof of the asymptotic formula  (\ref{asympP})
for $P_n$. Let $K$ be a compact subset of $\mathbb C \setminus \Gamma_P$.
We have the freedom to take the contours
$\Gamma_{P^-}$, $\Gamma_{P^+}$, $\Gamma_{P^{++}}$ near $\Gamma_P$, and
the circles $\Gamma_1$ and $\Gamma_2$ around $z_1$ and $z_2$ in such a way
that $K$ is in the exterior of these curves. Let $z \in K$.
Then we follow the transformations $Y \mapsto U \mapsto T \mapsto S$.
We see first from (\ref{eq:U}) that
\[ P_n(z) = Y_{21}(z) = U_{21}(z) e^{-(n+1) \ell} e^{(n+1)g_P(z)}. \]
Then from the definition of $T$ in (\ref{defT}) and (\ref{defTnearGR}),
we get that $U_{21}(z) = T_{21}(z)$.
We finally note that $T(z) = S(z) N(z)$, if $z$ is outside $\Gamma_3$ and $\Gamma_4$,
and $T(z) = S(z) M^{(j)}(z)$ if $z$ is inside $\Gamma_j$ with $j=3,4$.
For $T_{21}(z)$, we get in either case
\[ T_{21}(z) = S_{21}(z) N_{11}(z) + S_{22}(z) N_{21}(z) + S_{13}(z) N_{31}(z) \]
since the first columns of $M^{(3)}$ and $M^{(4)}$ agree with the first
column of $N$.

Since $S = I + \O\left(\frac{1}{n}\right)$, and since $N_{21}$ does not vanish on $K$,
we get
\[ P_n(z) = N_{21}(z) e^{-(n+1) \ell} e^{(n+1) g_P(z)}
    \left(1 + \O\left(\frac{1}{n}\right) \right)
\]
uniformly on $K$. Now we use the formula in (\ref{eq:N21N22N23}) for $N_{21}$
and we recall that $\ell = \log 2 - \pi i$, to obtain (\ref{asympP}).
\medskip

The proof of the asymptotic formula (\ref{asympR}) for $R_n$ is similar.
\medskip

For $Q_n$ we proceed differently, because of the way $Q_n$ appears in the
entries of $Y$. We take $z_0 \in \mathbb C \setminus \Gamma_Q$, and show that
there is a neighborhood $\Delta$ of $z_0$ such that (\ref{asympQ}) holds uniformly
for $z \in \Delta$. First we assume that $z_0$ belongs to the outside region
$D_{\infty}$. Then we can take the original contour $\Gamma$ so that
a neighborhood $\Delta$ of $z_0$ is in the outside region.
Then for $z \in \Delta$, we have by (\ref{eq:Yout})
\[ Q_n(z) = z^{3n+2} Y_{22}(z) = U_{22}(z) z^{3n+2} e^{-(n+1)(g_P(z) + g_R(z))}. \]
We can open the contours around $\Gamma_P$ and $\Gamma_R$ so that $\Delta$ is
in the exterior region to these contours. Then we have
$U = T = SN = \left(I + \O\left(\frac{1}{n}\right)\right)N$,
so that $U_{22}(z) = N_{22}(z) \left(1 + \O\left(\frac{1}{n}\right)\right)$
uniformly for $z \in \Delta$.
This leads to the first formula in (\ref{asympQ}), if
we use (\ref{eq:sumgs}) and the formula for $N_{22}$ in (\ref{eq:N21N22N23}).

If $z_0 \in D_P$, then we can also open up the lenses around $\Gamma_P$
and $\Gamma_R$ so that a neighborhood $\Delta$ of $z_0$ is not contained
in these lenses. Then we have for $z \in \Delta$, by (\ref{eq:Yin})
\[ Q_n(z) = E_n(z) - P_n(z) e^{-3nz} - R_n(z) e^{3nz} \]
and we need to find out, what is the dominant contribution as $n$ gets
large.
We already have asymptotic formulas for $P_n$ and $R_n$, from
which it follows that
\begin{equation} \label{Pnthroot}
    \frac{1}{n} \log |P_n(z)e^{-3nz}| \to \Re( g_P(z) - 3z - \ell),
\end{equation}
\begin{equation} \label{Rnthroot}
    \frac{1}{n} \log |R_n(z)e^{3nz}| \to \Re( g_R(z) + 3z - \ell).
\end{equation}
For $E_n$ it is easy to obtain in a similar way
\begin{equation} \label{Enthroot}
    \frac{1}{n} \log |E_n(z)| \to \Re (3 \log z - g_P(z) - g_R(z)).
\end{equation}
Now it turns out that for $z \in D_P$ the  term $P_n(z) e^{-3nz}$ dominates.
Indeed, we have by (\ref{eq:gPR1}) and (\ref{eq:gPR2})
\[ (g_P(z) - 3z - \ell) -  (g_R(z) + 3z - \ell) = (g_P - g_R - 6z)
    = 2(-\varphi_P(z) + \varphi_R(z)) \]
and we know that the real part of $-\varphi_P + \varphi_R$ is positive in
$D_P$. We also have by (\ref{eq:gPR1})
\[ (g_P(z) - 3z - \ell) - (3 \log z - g_P(z) - g_R(z)) = 2g_P(z) + g_R(z) - 3\log z- 3z - \ell
    = - 2\varphi_P(z) \]
and the real part of $-2\varphi_P$ is positive in $D_P$.
In view of (\ref{Pnthroot})--(\ref{Enthroot}) we then obtain the formula for $Q_n$
\[ Q_n(z) = - P_n(z) e^{-3nz} \left(1 + \O\left(\frac{1}{n}\right)\right) \]
uniformly for $z \in \Delta$, and from what we already know about $P_n$,
\[ Q_n(z) = - \frac{(-1)^{n+1} e^{n g_P(z)}}{2^n \sqrt{3\psi_P^4+1}} e^{-3nz}
     \left(1 + \O\left(\frac{1}{n}\right)\right) \]
uniformly for $z \in \Delta$. This leads to the second line in formula
(\ref{asympQ}) because
$g_P(z) - 3z + \log(-\frac{1}{2}) = g_Q(z)$ for $z \in D_P$, see
(\ref{eq:gPexplicit}) and (\ref{eq:gQexplicit}). The proof for the
case $z_0 \in D_R$ is similar.

We finally have to consider the case that $z_0$ is on
$\Gamma_P$ or on $\Gamma_R$ (but not one of the branch points).
We  consider only $z_0 \in \Gamma_P$, since $z_0 \in \Gamma_R$ will follow
in the same way. So let $z_0 \in \Gamma_P \setminus \{z_1,z_2\}$.
After opening up lenses around $\Gamma_P$ we have the contour
$\Gamma_{P^-}$ to the left of $z_0$, and the contour $\Gamma_{P^+}$
to the right. We can take a neighborhood $\Delta$ of $z_0$ that
is strictly contained in the domain bounded by $\Gamma_{P^-}$
and $\Gamma_{P^+}$. Then for $z \in \Delta \cap D_{\infty}$, we
have
\[ Q_n(z) = z^{3n+2} Y_{22}(z) = U_{22}(z) z^{3n+2} e^{-(n+1)(g_P(z) + g_R(z))} \]
and $U(z) = T(z) V^{-1}_{P^-}(z)$. Now we have
$U_{22}(z) = T_{22}(z)$, since the second column of
$V^{-1}_{P^-}$ is simply $\begin{pmatrix} 0 & 1 & 0
\end{pmatrix}^T$. We open the circles $\Gamma_1$ and $\Gamma_2$
around $z_1$ and $z_2$ so that $\Delta$ is in the exterior. Then
$T = SN = (I + \O(\frac{1}{n})) N$, and so
\[ Q_n(z) = z^{3n+2} N_{22}(z) z^{3n+2} e^{-(n+1)(g_P(z) + g_R(z))}
    \left(I + \O\left(\frac{1}{n}\right)\right) \]
uniformly for $z \in \Delta \cap D_{\infty}$. Using
(\ref{eq:sumgs}) and the formula (\ref{eq:N21N22N23}) for $N_{22}$, we
then find (\ref{asympQ}).

For $z \in \Delta \cap D_P$, we have
\begin{eqnarray*}
    Q_n(z) & = & E_n(z) - P_n(z) e^{-3nz} - R_n(z) e^{3nz} \\[10pt]
    & = & z^{3n+2} Y_{22}(z) - e^{-3nz} Y_{21}{z} - e^{3nz} Y_{13}(z) \\[10pt]
    & = & z^{3n+2} e^{-(n+1)(g_P(z) + g_R(z))} U_{22}(z)
        - e^{-3nz} e^{(n+1)(g_P(z) -\ell)} U_{21}(z) \\[10pt]
        & & - e^{3nz} e^{(n+1)(g_R(z) -\ell)}    U_{13}(z).
\end{eqnarray*}
Now we have that
\[ U(z) = T(z) V_{P^+}(z) V_{P^{++}}(z) \]
with
\[ V_{P^+}(z) V_{P^{++}}(z) = \begin{pmatrix} 1 & 0 & 0 \\
    z^{-1} e^{-3z} e^{2(n+1) \varphi_P(z)} & 1 & 0 \\
    0 & ze^{-3z} e^{-2(n+1) \varphi_R(z)} & 1 \end{pmatrix}. \]
Then
\begin{eqnarray*}
\lefteqn{ Q_n(z) = - e^{-3nz} e^{(n+1)(g_P(z) -\ell)} T_{21}(z) } \\[10pt]
   & & + \left[z^{3n+2} e^{-(n+1)(g_P(z) + g_R(z))}
    - e^{-3nz} e^{(n+1)(g_P(z) -\ell)} z^{-1} e^{-3z} e^{2(n+1) \varphi_P(z)}
        \right] T_{22}(z) \\[10pt]
   &&  + \left[z^{3n+2} e^{-(n+1)(g_P(z) + g_R(z))} ze^{-3z} e^{-2(n+1) \varphi_R(z)}
        - e^{3nz} e^{(n+1)(g_R(z) -\ell)}   \right]T_{13}(z).
\end{eqnarray*}
The funny thing is that the expressions multiplying $T_{22}(z)$ and $T_{13}(z)$
are exactly zero. This follows from (\ref{eq:gPR1}) and (\ref{eq:gPR2}).
Thus only the first term remains. This gives
\begin{eqnarray*}
     Q_n(z) & = & - e^{g_P(z) - \ell} e^{n g_Q(z)} T_{21}(z) \\[10pt]
     &  = & - e^{g_P(z) - \ell} N_{21}(z) e^{n g_Q(z)} \left(1 +
     \O\left(\frac{1}{n}\right)\right) \\[10pt]
     &  = & \frac{1}{\sqrt{\psi_P^4(z) + 1}} e^{n g_Q(z)} \left(1 + \O\left(\frac{1}{n}\right)\right)
\end{eqnarray*}
uniformly for $z \in \Delta \cap D_P$.
In the last step we used (\ref{eq:gPexplicit}) and the expression
for $N_{21}$. This proves the asymptotic formula for $z \in \Delta \cap D_P$.

The asymptotic formula (\ref{asympE}) for $E_n$ is established in
$D_P\cup D_R\cup \Gamma_Q$ again by following the transformations
$Y \mapsto U \mapsto T \mapsto S$. The asymptotics in the other
domains may be obtained by observing that the main contribution in the
sum
$$E_n(z)= P_n(z) e^{-3nz} + Q_n(z) + R_n(z) e^{3nz}$$
is given by $P_n(z) e^{-3nz}$, $Q_n(z)$ or $R_n(z) e^{3nz}$ if $z$
belongs
to the domains $D_{\infty,P}$, $D_{\infty,U}\cup D_{\infty,L}$ or
$D_{\infty,R}$ respectively.
This completes the proof of Theorem \ref{asymp}.
\end{proof}

\subsection{Proof of Theorem \ref{lim-zeros}}
\begin{proof}
The limits for the counting measures $\nu_{P_n}$, $\nu_{Q_n}$, and
$\nu_{R_n}$  follow from the strong convergence results established in
Theorem \ref{asymp}. The proof using the unicity theorem for logarithmic
potentials (see e.g.\ \cite[Theorem II.2.1]{SaffTotik})
is the same as the proof of \cite[Theorem 2.1]{stahl2}.

The proof for the limit of the measures $\nu_{E_n}$ is more difficult, since
these measures have unbounded support and infinite mass. Our proof is different
from the one given by Stahl \cite{stahl2}.

Since $E_n$ is an entire function of order $1$ with a zero of multiplicity
$3n+2$ at the origin, we have by the Hadamard factorization theorem
\[ E_n(z) = z^{3n+2} e^{A_n z + B_n} \prod_{j=1}^{\infty}
    \left(1 - \frac{z}{z_{j,n}}\right) e^{z/z_{j,n}}, \]
with $A_n, B_n \in \C$ and where $z_{j,n}$, $j = 1,2,\ldots$ are the
zeros of $E_n$ different from the origin. The zeros are counted according
to their multiplicities, and we assume that they are arranged in
non-decreasing absolute value, so that
$0 < |z_{1,n}| \leq |z_{2,n}| \leq \cdots$.
Taking logarithms, we get
\[ \frac{1}{n} \log \left(z^{-3n-2} E_n(z)\right)
    = \frac{1}{n} (A_n z + B_n)
        + \int \left(\log \left(1-\frac{z}{s}\right) + \frac{z}{s}\right) d\nu_{E_n}(s). \]
From the strong asymptotic formula (\ref{asympE}), we obtain from this that
\begin{equation} \label{limitgE}
    \lim_{n\to \infty} \left[ \frac{1}{n} (A_n z + B_n)
        + \int \left(\log \left(1-\frac{z}{s}\right) + \frac{z}{s}\right) d\nu_{E_n}(s) \right]
        = g_E(z)
\end{equation}
where
\begin{equation} \label{eq:gE}
    g_E(z) = \left\{ \begin{array}{ll}
       g_R(z) + 3z - 3 \log z - \ell & \mbox{ for } z \in D_{\infty, R}, \\[10pt]
       g_P(z) - 3z - 3 \log z - \ell & \mbox{ for } z \in D_{\infty, P}, \\[10pt]
       - g_P(z) - g_R(z) & \mbox{ for } z \in \C \setminus (\overline{D_{\infty,R} \cup D_{\infty,P}}).
    \end{array} \right.
\end{equation}
Note that $g_E$ is analytic across $\Gamma_P$ and $\Gamma_R$ because of the
relations (\ref{eq:gPRG1})--(\ref{eq:gPRG2}), so that $g_E$ is analytic in $\mathbb C \setminus \Gamma_E$.
The limit (\ref{limitgE}) holds uniformly for $z$ in compact subsets of $\C \setminus \Gamma_E$.

We first show that
\begin{equation} \label{formulagE}
    g_E(z) = \int \left(\log \left(1-\frac{z}{s}\right) + \frac{z}{s}\right) d\mu_{E}(s) - g_P(0) - g_R(0)
\end{equation}
for $z \in \mathbb C \setminus \Gamma_E$.
To prove (\ref{formulagE}) we assume $z \in \mathbb C \setminus (\overline{D_{\infty,R} \cup D_{\infty,P}})$.
We take the derivative of the right-hand side of (\ref{formulagE}), which is
\begin{eqnarray} \nonumber
    \int \left( \frac{1}{z-s} + \frac{1}{s} \right) d\mu_E(s)
    & = & \frac{3}{2\pi i} \int_{\Gamma_{E,1} \cup \Gamma_{E,2}} \left( \frac{1}{z-s} + \frac{1}{s} \right)
        (\psi_Q - \psi_P)(s) ds \\
    & & + \frac{3}{2\pi i} \int_{\Gamma_{E,3} \cup \Gamma_{E,4}} \left( \frac{1}{z-s} + \frac{1}{s} \right)
        (\psi_Q - \psi_R)(s) ds.
        \label{eq:derivrhs}
\end{eqnarray}
The integral over $\Gamma_{E,1} \cup \Gamma_{E,2}$ can be deformed to the integral
over the $-$side of $\Gamma_P$, which leads to
\begin{eqnarray*}
    \frac{3}{2\pi i} \int_{\Gamma_{E,1} \cup \Gamma_{E,2}} \left( \frac{1}{z-s} + \frac{1}{s} \right)
        (\psi_Q- \psi_P)(s) ds
    & = & \frac{3}{2\pi i} \int_{\Gamma_P} \left( \frac{1}{z-s} + \frac{1}{s} \right) (\psi_Q - \psi_P)_-(s) ds \\
    & = & \frac{3}{2\pi i} \oint_{\gamma} \left(\frac{1}{z-s} + \frac{1}{s} \right) \psi_P(s) ds
\end{eqnarray*}
where $\gamma$ is a closed contour going around $\Gamma_P$ in the negative direction and with
$z$ and $0$ outside $\gamma$, see also the proof of Lemma \ref{lem:gder}.
We compute the integral using the residue theorem for the exterior
of $\gamma$, for which there is a contribution from the poles at $z$ and $0$. The result is
that the first term on the right-hand side of (\ref{eq:derivrhs}) is equal to $- 3 \psi_P(z) + 3 \psi_P(0)$.
In the same way, we find that the second term is equal to $-3\psi_R(z) + 3 \psi_R(0)$.
Using (\ref{eq:gPder}), (\ref{eq:gRder}) and the fact that $\psi_R(0) =  \sqrt{1/3} = - \psi_P(0)$,
we see that the derivative of the right-hand side of (\ref{formulagE}) is equal to $-g_P'(z) - g_R'(z)$
for $z \in \mathbb C \setminus (\overline{D_{\infty,R} \cup D_{\infty,P}})$.
Because of the definition of $g_E$ in (\ref{eq:gE}) we then get (\ref{formulagE}) in case
$z \in \mathbb C \setminus (\overline{D_{\infty,R} \cup D_{\infty,P}})$.
The formula (\ref{formulagE}) in the other domains follows by analytical continuation.

Knowing (\ref{formulagE}) we get from (\ref{limitgE}) and (\ref{eq:gE}) that
$\lim\limits_{n\to\infty} \frac{1}{n} A_n = 0$,
$\lim\limits_{n\to\infty} \frac{1}{n} B_n = - g_P(0) - g_R(0)$, and
\begin{equation} \label{eq:convnuEns}
    \lim_{n\to\infty} \int \left(\log \left(1-\frac{z}{s}\right) + \frac{z}{s}\right) d\nu_{E_n}(s)
    = \int \left(\log \left(1-\frac{z}{s}\right) + \frac{z}{s}\right) d\mu_E(s)
\end{equation}
for $z \in \C \setminus \Gamma_E$.
\medskip

We next consider the measures $\rho_n$ defined by
$$
    d\rho_n(s) = \frac{1}{|s|^2} d\nu_{E_n}(s).
$$
It will be our goal to prove that
\begin{equation} \label{eq:rhonestimate}
    \sup_{n \in \mathbb N} \int d\rho_n < +\infty
\end{equation}
and that the sequence $(\rho_n)$ is tight, that is,
\begin{equation} \label{eq:rhontight}
    \forall \varepsilon > 0 : \exists R > 0 : \forall n\in\mathbb N :
    \int_{|s| \geq R} d\rho_n(s) \leq \varepsilon,
\end{equation}
see \cite[Chapter 6]{Billing}.
We obtain (\ref{eq:rhonestimate}) and (\ref{eq:rhontight}) from
lower bounds on the absolute values of the zeros $z_{j,n}$ of $E_n$.

Let $N_n(r)$ be the number of non-zero roots of $E_n$ with absolute
value $\leq r$. From the theory of entire functions, we have that
\begin{equation} \label{eq:Nnestimate}
    N_n(r) < \log \max_{|z| = er} |f_n(z)| - \log |f_n(0)|,
\end{equation}
see \cite[Section 2.5]{Levin}, where
\begin{equation} \label{deffnz}
    f_n(z) = z^{-3n-2} E_n(z).
\end{equation}
The integral representations (\ref{rep-intP})--(\ref{rep-intR})
and (\ref{eq:constant}) imply that
\begin{equation} \label{rep-intE}
    E_n(z) = \frac{(-1)^{n+1} n!}{(3n)^n} \frac{1}{2\pi i}
        \oint_{C_R} \frac{e^{3nzw}dw}{[w(w^2-1)]^{n+1}},
\end{equation}
where $C_R$ is a circle centered at the origin with radius $R > 1$.
From (\ref{rep-intE}) it immediately follows that
\[  E_n(z) = \frac{(-1)^{n+1} n! (3n)^{2n+2}}{(3n+2)!} z^{3n+2} + \O(z^{3n+3})
    \qquad \mbox{as } z \to 0,
\]
and so
\begin{equation} \label{fn0}
    f_n(0) = \frac{(-1)^{n+1} n! (3n)^{2n+2}}{(3n+2)!}.
\end{equation}
Using $E_n(z) = P_n(z) e^{-3nz} + Q_n(z) + R_n(z) e^{3nz}$
and the fact that the polynomials $P_n$, $Q_n$, and $R_n$, of leading
coefficients
$(-1/2)^{n+1}$, 1, and  $(-1/2)^{n+1}$ respectively, have
their zeros in a compact set (independent of $n$), we easily get
that $|E_n(z)| \leq e^{4n|z|}$ for every $n \in \mathbb N$ and for
every $|z| \geq R$ with $R$ sufficiently large, say $R \geq R_0 \geq 1$.
Then, from (\ref{eq:Nnestimate}), (\ref{deffnz}) and (\ref{fn0}) we see
that there exists a constant $C >0$ so that $N_n(r) < Cn r$ if $r > R_0$.
Hence, for every $j$ with $|z_{j,n}| > R_0$, we have $j \leq N_n(|z_{j,n}|)
< Cn |z_{j,n}|$, so that
\begin{equation} \label{zjnestimate}
    |z_{j,n}| > \frac{j}{Cn} \qquad \mbox{if } |z_{j,n}| > R_0.
\end{equation}
This is our required estimate on $|z_{j,n}|$.

To estimate
\[ \int d\rho_n = \frac{1}{n} \sum_{j=1}^{\infty} \frac{1}{|z_{j,n}|^2}, \]
we note that there are at most $CnR_0$ zeros with absolute value $\leq R_0$.
So for $j > CnR_0$, we have $|z_{j,n}| > R_0$, and we can use the
estimate (\ref{zjnestimate}). For $j \leq CnR_0$, we use the fact that
the zeros can only accumulate on $\Gamma_E$, which implies that there is a
$\delta > 0$ so that $|z_{j,n}| > \delta$ for every $j$ and $n$.
Then for every $n \in \mathbb N$,
\[   \int d\rho_n \leq
     \frac{1}{n} \sum_{j=1}^{[CnR_0]} \frac{1}{\delta^2}
        + \frac{1}{n} \sum_{j=[CnR_0]+1}^{\infty} \frac{(Cn)^2}{j^2}
     \leq \frac{CR_0}{\delta^2} + \frac{C^2 n}{[CnR_0]},
\]
and (\ref{eq:rhonestimate}) follows.

Next, we estimate for $R > R_0$,
\begin{eqnarray*}
    \int_{|s| \geq R} d\rho_n(s) & = &
        \frac{1}{n} \sum_{j=N_n(R)+1}^{\infty} \frac{1}{|z_{j,n}|^2}
    \,  \leq \, \frac{1}{n} \sum_{j=N_n(R)+1}^{\infty}
        \min\left(\frac{1}{R^2}, \frac{C^2n^2}{j^2}\right) \\
    & \leq & \frac{1}{n} \sum_{j=1}^{[CnR]} \frac{1}{R^2}
        + \frac{1}{n} \sum_{j=[CnR]+1}^{\infty} \frac{C^2n^2}{j^2}
    \, \leq \, \frac{C}{R} + \frac{C^2n}{[CnR]},
\end{eqnarray*}
and this shows that, given $\varepsilon >0$, we can choose $R > R_0$ sufficiently
large so that $\int_{|s| \geq R} d\rho_n(s) \leq \varepsilon$ for every $n$.
Hence (\ref{eq:rhontight}) holds and the sequence $(\rho_n)$ is tight.
\medskip

Now let $\rho^*$ be the limit of a weakly convergent subsequence of
$(\rho_n)$, where $\rho_n \to \rho$ weakly means that $\int f d\rho_n
\to \int f d\rho$
for every bounded continuous function $f$.
Because of (\ref{eq:rhonestimate}), we have that $\int d\rho^* < \infty$.
Let $d\nu^*(s) = |s|^2 d\rho^*(s)$. Then (\ref{eq:convnuEns}) implies
\begin{equation} \label{eq:intnuintmuE}
    \int \left(\log \left(1-\frac{z}{s}\right) + \frac{z}{s} \right) d\nu^*(s) =
    \int \left(\log \left(1-\frac{z}{s}\right) + \frac{z}{s} \right) d\mu_E(s)
\end{equation}
for every $z \in \mathbb C \setminus \Gamma_E$.
Let $R > 0$ and let $\nu^*_R$ and $\mu_{E,R}$ be the restrictions to $\{|z| < R\}$
of $\nu^*$ and $\mu_E$, respectively. Then we get by taking the real part
and rewriting (\ref{eq:intnuintmuE}) that
\[ \int \log |z-s| d\nu^*_R(s) = \int \log|z-s| d\mu_{E,R}(s) + h_R(z), \]
where $h_R$ is harmonic in $\{|z| < R \}$. Since $\Gamma_E$ has
two--dimensional Lebesgue measure 0, the unicity theorem \cite[Theorem
II.2.1]{SaffTotik} applies, showing
that $\nu^*_R$ and $\mu_{E,R}$ are equal. As $R$ is arbitrary, we have
that $\nu^* = \mu_E$ and so $d\rho^*(s) = \frac{1}{|s|^2} d\mu_E(s)$.
Since the sequence $(\rho_n)$ is tight, and $\frac{1}{|s|^2} d\mu_E(s)$
is the only possible limit of a weakly convergent subsequence, the full sequence
$(\rho_n)$ converges weakly to $\frac{1}{|s|^2} d\mu_E(s)$,
see \cite[Chapter 6]{Billing}.
Then it also follows that the measures $\nu_{E_n}$ converge to $\mu_E$
in the sense indicated in the theorem.

This completes the proof of the theorem.
\end{proof}

\subsection{Proof of Theorem \ref{asym-Q}}
\begin{proof}
Away from $\Gamma_Q$ we have that one of the terms in (\ref{asympQ2}) dominates
the others, and then (\ref{asympQ2}) reduces to (\ref{asympQ}).
So it suffices to prove that (\ref{asympQ2}) holds uniformly
for $z$ in a neighborhood $\Delta$ of a point $z_0 \in \Gamma_Q
\setminus \{z_1, z_2,z_3,z_4\}$.
So let $z_0 \in  \Gamma_Q \setminus \{z_1, z_2,z_3,z_4\}$. We take
the contours $\Gamma_{P^{++}}$ and $\Gamma_{R^{++}}$ so that
we can choose a disk  $\Delta$ centered at $z_0$ that lies
in the region bounded by $\Gamma_{P^{++}}$, $\Gamma_{R^{++}}$,
$\Gamma_{U^m}$ and $\Gamma_{L^m}$.

From (\ref{En}) we know that
\begin{equation}
\label{Qn}
Q_n(z)= E_n(z) - P_n(z) e^{-3nz} - R_n(z) e^{3nz}.
\end{equation}
We already know asymptotic expressions for the polynomials
$P_n$, $R_n$ and the remainder function $E_n$ that are valid in $\Delta$,
namely
(\ref{asympP}), (\ref{asympR}) and the last equation in
(\ref{asympE}).
Plugging these three formulas into
(\ref{Qn}), we deduce
$$
Q_n(z)  =  z^{3n}e^{-n[g_P(z)+g_R(z)]}
\left[\frac{e^{n[2g_P(z)+g_R(z)-\ell-3\log z-3z]}}{\sqrt{3\psi_P^4(z)+1}}
\left(1 +
     \O\left(\frac{1}{n}\right)\right)
-\frac{\left(1 +
     \O\left(\frac{1}{n}\right)\right)}{\sqrt{3\psi_Q^4(z)+1}}
\right.
$$
$$
\left.
 +\frac{e^{n[g_P(z)+2g_R(z)-\ell-3\log z+3z]}}{\sqrt{3\psi_R^4(z)+1}}
\left(1 +
     \O\left(\frac{1}{n}\right)\right)\right].
$$
Using (\ref{eq:gPR1}) and (\ref{eq:gPR2}) we are led to (\ref{asympQ2}),
uniformly for $z\in \Delta$.

The proof of (\ref{asympE2}) is similar. We use (\ref{asympP}),
(\ref{asympR}) and the first equation in
(\ref{asympQ}). \end{proof}

\subsection{Proof of Theorem \ref{asym-PR}}
\begin{proof}
Away from $\Gamma_P$ one of the terms in (\ref{asympP1}) is dominating
the other, and then (\ref{asympP1}) leads to (\ref{asympP}).
So it suffices to consider $z_0 \in \Gamma_P\setminus \{z_1,z_2\}$, and prove
that (\ref{asympP1}) holds uniformly for $z$ in a disk $\Delta$
centered at $z_0$.
We may choose this disk so that it is disjoint from $\Gamma_1$ and
$\Gamma_2$, and also so that it is contained in the region bounded
by $\Gamma_{P^-}$ and $\Gamma_{P^+}$, see Figure \ref{fig:t6}.

We consider $z \in \Delta$. Then we have as before
$P_{n}(z)=U_{21}(z) e^{(n+1)(g_{P}(z)-\ell)}$, but the relation between
$T$ and $U$ is $U(z)=T(z)V_{P^-}^{-1}(z)$, or $U(z) = T(z) V_{P^+}(z) V_{P^{++}}$
depending on whether $z$ lies in the region bounded by $\Gamma_P$ and
$\Gamma_{P^-}$ (the $-$side of $\Gamma_P$)
or in the region bounded by $\Gamma_P$ and $\Gamma_{P^+}$ (the $+$side).
This gives
\[ U_{21}(z) = T_{21}(z) \pm z^{-1} e^{-3z} e^{2(n+1)\varphi_P(z)} T_{22}(z), \]
for $z$ on the $\pm$side of $\Gamma_P$.

From the fact that $T = SN$ with $S = I + \O(\frac{1}{n})$, we get
\begin{multline}
 P_{n}(z) = e^{(n+1)(g_{P}(z)-\ell)}\left(N_{21}(z)
 \pm z^{-1}e^{2(n+1)\varphi_{P}(z)-3z} N_{22}(z)\right)
     \left( 1 + \O \left( \frac1n \right) \right)\\
+
\left(N_{11}(z)
 \pm z^{-1}e^{2(n+1)\varphi_{P}(z)3z} N_{12}(z)\right)
     \O \left( \frac1n \right)\\
+
\left(N_{31}(z)
 \pm z^{-1}e^{2(n+1)\varphi_{P-}(z)3z} N_{3,2}(z)\right)
     \O \left( \frac1n \right).
\end{multline}
Since all functions $N_{jk}$ do not vanish in $\Delta$, and since
$e^{\ell} = -2$, we obtain
\[ (-2)^{n+1} P_n(z) = e^{(n+1)g_P(z)}
    \left[ N_{21}(z) \left( 1 + \O \left( \frac1n \right) \right)
    \pm z^{-1}e^{2(n+1)\varphi_{P}(z)3z} N_{22}(z)
    \left(1+ \O \left( \frac1n \right) \right) \right]. \]
Now we use (\ref{eq:gPR1}) and the formulas in
(\ref{eq:N21N22N23}) for $N_{21}$ and $N_{22}$,
and (\ref{asympP1}) follows.

The proof of the asymptotic formula (\ref{asympR1}) for $R_n$ is similar.
\end{proof}

\subsection{Proof of Theorem \ref{asym-branch}}
\begin{proof}
If we unravel all the transformations for
$z \in \Delta_1$, we find that $Y(z)$ is a product of no less than
$10$ to $14$ matrices, the exact number
depending on the region where $z$ is. Recall that $z$ belongs
to one of the four regions $f^{-1}(I)$, $f^{-1}(II)$, $f^{-1}(III)$,
and $f^{-1}(IV)$, see Figure \ref{fig:Ai}.
For $z \in f_1^{-1}(III)$, we have the
following, where $s = (n+1)^{2/3} f_1(z)$,
\begin{multline} \label{fullproduct}
    Y(z) = e^{i\pi/6} \sqrt{\pi}
    \begin{pmatrix} e^{(n+1)\ell/3} & 0 & 0 \\ 0 & e^{-2(n+1)\ell/3} & 0 \\
            0 & 0 & e^{(n+1)\ell/3} \end{pmatrix}
    S(z) N(z) \\
    \begin{pmatrix} z^{1/2} e^{3z/2} & 0 & 0 \\ 0 & z^{-1/2} e^{-3z/2} & 0 \\ 0 & 0 & 1 \end{pmatrix}
    \begin{pmatrix} 1 & -1 & 0 \\ i & i & 0 \\ 0 & 0 & 2e^{-i\pi/6}
    \end{pmatrix}
    \begin{pmatrix} s^{1/4} & 0 & 0 \\ 0 & s^{-1/4} & 0 \\ 0 & 0 & 1 \end{pmatrix}
    \begin{pmatrix} \Ai(s) & -\Ai(\omega_3^2s) & 0 \\
    \Ai'(s) & -\omega_3^2\Ai'(\omega_3^2s) & 0 \\
    0 & 0 & 1             \end{pmatrix} \\
    \begin{pmatrix} e^{-i\pi/6} & 0 & 0 \\ 0 & e^{i\pi/6} & 0 \\ 0 & 0 & 1 \end{pmatrix}
    \begin{pmatrix} 1 & 0 & 0 \\ -1 & 1 & 0 \\ 0 & 0 & 1 \end{pmatrix}
    \begin{pmatrix} z^{-1/2} e^{-3z/2} e^{(n+1)\varphi_P(z)} & 0 & 0 \\
    0 & z^{1/2} e^{3z/2} e^{-(n+1)\varphi_P(z)} & 0 \\ 0 & 0 & 1 \end{pmatrix} \\
    \begin{pmatrix} 1 & 0 & 0 \\
    z^{-1}e^{-3z} e^{2(n+1)\varphi_P(z)} & 1 & 0 \\
    0 & 0 & 1 \end{pmatrix}
    \begin{pmatrix} 1 & 0 & 0 \\ 0 & 1 & 0 \\
    0 & ze^{-3z}e^{-2(n+1)\varphi_R(z)} & 1 \end{pmatrix} \\
    \begin{pmatrix} e^{-(n+1)\ell/3} & 0 & 0 \\
    0 & e^{2(n+1)\ell/3} & 0 \\ 0 & 0 & e^{-(n+1)\ell/3}
          \end{pmatrix}
\begin{pmatrix}
                      e^{(n+1)g_P(z)} & 0 & 0 \\
                      0 & e^{-(n+1)[g_P(z)+g_R(z)]} & 0 \\
                      0 & 0 & e^{(n+1)g_R(z)}
                      \end{pmatrix}.
    \end{multline}

Recall that $P_n(z)$ is the $(2,1)$ entry of $Y(z)$. We multiply
(\ref{fullproduct}) with the first unit vector
$\begin{pmatrix} 1 & 0 & 0 \end{pmatrix}^T$,
and we obtain after some simple calculations, where
$*$ denotes an unspecified unimportant entry
\begin{equation} \label{Pnnearz11}
 \begin{pmatrix} * \\ P_n(z) \\ * \end{pmatrix} =
    \frac{\sqrt{\pi}e^{(n+1)(g_P(z)+ \varphi_P(z))}}{(-2)^{n+1}}
    S(z) N(z) \begin{pmatrix} 1 & 0 & 0 \\ 0 & z^{-1} e^{-3z} & 0 \\ 0 & 0 & * \end{pmatrix}
    \begin{pmatrix} 1 & -1 & 0 \\ i & i & 0 \\ 0 & 0 & * \end{pmatrix}
    \begin{pmatrix} s^{1/4}\Ai(s) \\ s^{-1/4}\Ai'(s) \\ 0 \end{pmatrix}.
\end{equation}
Similar calculations for $z$ in the other regions near $z_1$ lead to exactly
the same expression (\ref{Pnnearz11}).
Then further calculations show that for $z \in \Delta_1$,
\begin{multline} \label{Pnnearz12}
    (-2)^{n+1} P_n(z) =
    \sqrt{\pi} e^{(n+1)(g_P(z)+ \varphi_P(z))} \\
    \left[ \sum_{k=1}^3 S_{2k}(z) \left(N_{k1}(z) + iz^{-1}e^{-3z} N_{k2}(z) \right)
    (n+1)^{1/6} f_1(z)^{1/4} \Ai((n+1)^{2/3} f_1(z))
    \right. \\ \left.
    + \sum_{k=1}^3 S_{2k}(z) \left(-N_{k1}(z) + iz^{-1}e^{-3z} N_{k2}(z) \right)
    (n+1)^{-1/6} f_1(z)^{-1/4} \Ai'((n+1)^{2/3} f_1(z))
    \right]. \end{multline}

From the jump condition (\ref{eq:NjumpP}) for $N$ on $\Gamma_P$
it easily follows that for each $k$,
\[ \left(N_{k1}(z) + iz^{-1}e^{-3z} N_{k2}(z) \right)_+
    = i \left(N_{k1}(z) + iz^{-1}e^{-3z} N_{k2}(z) \right)_-, \qquad
    z \in \Gamma_P. \]
The fourth root in $f_1(z)^{1/4}$ is defined with a cut along $\Gamma_P$,
and on $\Gamma_P$ there is a jump
\[ \left(f_1(z)^{1/4}\right)_+ = -i \left(f_1(z)^{1/4}\right)_-, \qquad z \in \Gamma_P. \]
Thus the products
\[ \left(N_{k1}(z) + iz^{-1} e^{-3z} N_{k2}(z) \right) f_1(z)^{1/4} \]
are analytic across $\Gamma_P$. Since the functions $N_{kj}$ have a fourth root
singularity at $z_1$, and $f_1(z)^{1/4}$ has a fourth root zero at $z_1$,
these functions are non-zero at $z_1$, and therefore also in a neighborhood
of $z_1$. By shrinking $\Delta_1$, if necessary, we may assume that these
functions are without zeros in $\Delta_1$.
Similarly, we have that
\[ \left(-N_{k1}(z) + iz^{-1}e^{-3z} N_{k2}(z) \right) f_1(z)^{-1/4} \]
is analytic and without zeros in $\Delta_1$.

Now we recall that $S(z) = I + \O\left(\frac{1}{n}\right)$.
So it follows that in the two sums in (\ref{Pnnearz12}) the terms with $k=2$
are dominating, and we get
\begin{multline} \label{Pnnearz13}
    (-2)^{n+1} P_n(z) =
    \sqrt{\pi} e^{(n+1)(g_P(z)+ \varphi_P(z))} \\
    \left[ n^{1/6} \left(N_{21}(z) + iz^{-1}e^{-3z} N_{22}(z) \right)
    f_1(z)^{1/4} \Ai((n+1)^{2/3} f_1(z)) \left(1 + \O\left(\frac{1}{n}\right)\right)
    \right. \\ \left.
    + n^{-1/6} \left(-N_{21}(z) + iz^{-1}e^{-3z} N_{22}(z) \right)
     f_1(z)^{-1/4} \Ai'((n+1)^{2/3} f_1(z)) \left(1 + \O\left(\frac{1}{n}\right)\right)
    \right].
\end{multline}
This proves (\ref{asympP3}). We have also shown above that the
functions $h_1$ and $h_2$ are analytic and without zeros in $\Delta_1$.
\medskip

For the asymptotic formula for $Q_n$ near $z_1$, we proceed in a similar way.
Actually, it is easier to do the calculations for $z$ in the region
$f_1^{-1}(II)$ or $f_1^{-1}(I)$. The result is that for $z \in \Delta_1$,
\begin{multline} \label{Qnnearz11}
 \begin{pmatrix} * \\ Q_n(z) \\ * \end{pmatrix}
     = \sqrt{\pi} e^{3z} e^{-(n+1) (g_P(z)+\varphi_P(z)+g_R(z)-3\log
     z)}  S(z) N(z) \\
     \begin{pmatrix} 1 & 0 & 0 \\ 0 & z^{-1} e^{-3z} & 0 \\ 0 & 0 & * \end{pmatrix}
        \begin{pmatrix} 1 & -1 & 0 \\ i & i & 0 \\ 0 & 0 & * \end{pmatrix}
        \begin{pmatrix}  s^{1/4} \omega_3^2 \Ai(\omega_3^2 s) \\ s^{-1/4} \omega_3 \Ai'(\omega_3^2 s) \\ 0
        \end{pmatrix},
\end{multline}
where again $s = (n+1)^{2/3} f_1(z)$.
We make use of (\ref{eq:gPR1}) to replace the exponential factor on the
right-hand side of (\ref{Qnnearz11}) by
$e^{(n+1)(g_P(z)+\varphi_P(z)-3z-l)}$. Then,
we obtain (\ref{asympQ3}) from (\ref{Qnnearz11}) in the same way
as we obtained (\ref{asympP3}) from (\ref{Pnnearz11}).

For the asymptotic formula for $E_n$ near $z_1$, we may again proceed
in a similar way, doing the calculations for $z$ in the region
$f_1^{-1}(III)$ or $f_1^{-1}(IV)$.
\end{proof}

\subsection{Proof of Corollary \ref{asym-zero}}
\begin{proof}
The behavior of the extreme zeros of $P_n$ near $z_1$ follows from
the asymptotic formula (\ref{asympP3}). Indeed, we consider
the function
\[ F_n(t) = \frac{(-2)^{n+1} e^{-(n+1)(g_P(z) + \varphi_P(z))}}{\sqrt{\pi} n^{1/6} h_1(z)}
    P_n(z),
    \qquad \mbox{where }z = z_1 + t n^{-2/3}. \]
Then $F_n$ has zeros $t_{\nu,n} = (z_{\nu,n}^P - z_1) n^{2/3}$,
$\nu=1,2,\ldots$, and these zeros are  ordered by increasing absolute value.
Because of (\ref{asympP3}) we have that
\begin{multline}
F_n(t) = \Ai\left( (n+1)^{2/3}
  f_1(z_1+tn^{-2/3})\right) \left(1 + \O\left(\frac{1}{n}\right)\right)
  \\
    + n^{-1/3} \frac{h_2(z_1+tn^{-2/3})}{h_1(z_1+tn^{-2/3})} \Ai'\left((n+1)^{2/3}
  f_1(z_1+tn^{-2/3})\right)
  \left(1+\O\left(\frac{1}{n} \right)\right).
\end{multline}
Since $f_1$ is an analytic function with a simple zero at
$z_1$ and $f_1'(z_1)=c_1$ (see (\ref{def-f}), (\ref{c_1})), and since $h_1$ and
$h_2$ are nonzero analytic functions near $z_1$, we get by expanding
these functions near $z_1$ that
\begin{equation}
\label{F_n}
F_n(t)= \Ai\left(tc_1+\O\left({n}^{-2/3} \right)\right)
+ n^{-1/3} \frac{h_2(z_1)}{h_1(z_1)} \Ai'\left(tc_1+\O\left({n}^{-2/3} \right)\right)
+\O\left({n}^{-2/3} \right).
\end{equation}
Now, expanding the Airy function $\Ai$ near $tc_1$, and observing that
the second term of the sum in the right-hand side of (\ref{F_n}) is of
order $n^{-1/3}$, we get simply
\begin{equation} \label{F_n2}
    F_n(t) = \Ai(tc_1)+\O\left({n}^{-1/3} \right).
\end{equation}
The $\O$-term holds uniformly on compact subsets of the complex $t$-plane.
From Hurwitz' theorem it follows that for every fixed $\nu \in \mathbb N$, we
have
\[ \lim_{n \to \infty} t_{\nu,n} = -\frac{\iota_{\nu}}{c_1}. \]
Using the fact that $-\iota_{\nu}$ is a simple zero of the Airy function,
we obtain from (\ref{F_n2}) that
\[ t_{\nu,n} = - \frac{\iota_{\nu}}{c_1} + \O\left({n}^{-1/3} \right). \]
This proves (\ref{zerosPn}), since $z_{\nu,n}^P = z_1 + t_{\nu,n} n^{-2/3}$.

The formulas (\ref{zerosQn}) and (\ref{zerosEn}) for the extreme zeros
of $Q_n$ and $E_n$ near $z_1$ are obtained in a similar way from the
asymptotics of
$Q_n$ and $E_n$ near $z_1$, given by (\ref{asympQ3}) and (\ref{asympE3}),
respectively.
\end{proof}

\subsection{Proof of Theorem \ref{algebra}}
\begin{proof}
We consider $z$ in a compact subset $K$ of $D_R$.
By the third equation in (\ref{asympQ}), we have
$$ \lim_{n \to \infty} \frac{1}{n}\log|Q_n^2(z)| = 2\Re g_Q(z),$$
and by (\ref{asympP}) and (\ref{asympR}),
$$ \lim_{n\to\infty} \frac{1}{n}\log|P_n(z)R_n(z)| = \Re (g_P(z)+g_R(z)-2 \ell).$$
In $D_R$ we may use (\ref{relgRgQ}), so that
$$ 2\Re g_Q(z) - \Re (g_P(z)+g_R(z)-2\ell) = \Re(g_R(z)-g_P(z)+6z).$$
By (\ref{eq:gPR1}) and (\ref{eq:gPR2}), this is
$2\Re(\varphi_P-\varphi_R)$ which is positive in $D_R$. Thus,
$4P_nR_n/Q_n^2$ is exponentially small, and there is a choice of the
square root so that
$$ \sqrt{Q_n^2(z)-4P_n(z)R_n(z)} = Q_n(z) \left(1+\O(e^{-cn})\right),$$
as $n\to\infty$, uniformly for $z \in K$, for some constant $c>0$.
Then we choose the $-$ sign in (\ref{formulaXn}) so that
$$ X_{n}(z)=\frac{-Q_n(z)-\sqrt{Q_n^2(z)-4P_n(z)R_n(z)}}{2R_n(z)}=
-\frac{Q_n(z)}{R_n(z)} \left(1+\O(e^{-cn})\right).$$
Now, we use (\ref{asympR}) and the third equation in (\ref{asympQ}) to
obtain
$$ X_{n}(z)=\frac{-(-2)^{n+1}e^{ng_Q(z)}}{2e^{ng_R(z)}}
    \left(1+\O\left(\frac{1}{n}\right)\right)=
    e^{3nz}\left(1+\O\left(\frac{1}{n}\right)\right),$$
where the last equality uses (\ref{relgRgQ}). This proves (\ref{conv-Xn})
for $z \in D_R$.

The case $z\in D_P$ is similar, but now we have to take
the $+$ sign in the formula (\ref{formulaXn}) for $X_n$.
This completes the proof of Theorem \ref{algebra}.
\end{proof}

\obeylines
\texttt{
A. B. J. Kuijlaars (arno@wis.kuleuven.ac.be)
W. Van Assche (walter@wis.kuleuven.ac.be)
Department of Mathematics
Katholieke Universiteit Leuven
Celestijnenlaan 200B
B-3001 Leuven, BELGIUM
\medskip

F. Wielonsky (wielonsk@ano.univ-lille1.fr)
Laboratoire de Math\'ematiques appliqu\'ees
FRE CNRS 2222 - Bat. M2
Universit\'e des Sciences et Technologies Lille 1
F-59655 Villeneuve d'Ascq Cedex, FRANCE
}

\end{document}